\documentclass[11pt]{amsart}
\usepackage{amssymb,amsmath,amsthm,graphicx,booktabs,xcolor}
\usepackage[margin=1.2in]{geometry}
\usepackage{tikz}
\newtheorem{theorem}{Theorem}[section]
\newtheorem{lemma}[theorem]{Lemma}
\newtheorem{proposition}[theorem]{Proposition}
\newtheorem{corollary}[theorem]{Corollary}
\theoremstyle{definition}\newtheorem{definition}[theorem]{Definition}
\theoremstyle{remark}\newtheorem{remark}[theorem]{Remark}
\numberwithin{equation}{section}
\newcommand{\Om}{\Omega}\newcommand{\C}{\mathbb{C}}\newcommand{\R}{\mathbb{R}}
\newcommand{\dz}{\partial_z}\newcommand{\dzb}{\partial_{\bar z}}
\renewcommand{\Re}{\operatorname{Re}}\renewcommand{\Im}{\operatorname{Im}}\newcommand{\dist}{\operatorname{dist}}
\newcommand{\kess}{\kappa_{\mathrm{ess}}}
\newcommand{\A}{A^2}\newcommand{\Az}{A^2_0}
\DeclareMathOperator{\Div}{div}\DeclareMathOperator{\rot}{rot}\DeclareMathOperator{\curl}{curl}

\title[Verified bounds for the inf-sup constant on polygons]{Verified two-sided bounds for the inf-sup constant\\ of the divergence on planar polygons}
\author{Xuefeng Liu}
\address{Department of Information and Sciences, Tokyo Woman's Christian University, 2-6-1 Zempukuji, Suginami-ku, Tokyo 167-8585, Japan}
\email{xfliu@lab.twcu.ac.jp}
\begin{document}
\begin{abstract}
The inf-sup constant $\beta(\Om)$ of the divergence governs the stability of the Stokes problem. For polygons, analytic lower bounds are far
from the true value, and numerical methods either bound $\beta$ from above or require spectral-gap and mesh-fineness hypotheses that cannot
presently be checked. We propose a method that produces verified lower bounds without any such hypothesis. It rests on the Friedrichs
formulation of $\beta$ in the Bergman space and on a multiplier identity that controls the non-compact corner contribution pointwise, with the
sharp corner constant; the remaining compact part is bounded by a finite matrix inequality that is checked in exact and ball arithmetic.
We prove that the bounds converge to $\beta$. Together with enclosed Rayleigh quotients of explicit trial functions (and, for the square, the
bottom $\frac12-\frac1\pi$ of the essential spectrum) this yields two-sided enclosures, e.g.\ $0.426245\le\beta\le0.426252$ for the square and
$0.303973\le\beta\le0.304078$ for the L-shaped domain. For a square with a square hole, which is neither simply connected nor star-shaped, we
obtain $0.154596\le\beta\le0.154781$. For the rectangles $(0,1)\times(0,r)$ we enclose $\beta$ for eight aspect ratios and prove for the five ratios
$r\in\{\frac14,\frac38,\frac12,\frac9{16},\frac58\}$ that $\beta^2$ is an isolated eigenvalue below the essential spectrum; the last two lie beyond the
range covered by the explicit bound of Costabel, Crouzeix, Dauge and Lafranche. As an application, the certified values give guaranteed a posteriori
bounds for Taylor--Hood approximations of the Stokes equations, including the doubly connected domain, for which, to our knowledge, no explicit lower
bound for $\beta$ had been available.
\end{abstract}
\maketitle

\section{Introduction}\label{sec:intro}
Throughout, $\Om\subset\R^2$ is a bounded Lipschitz domain, and the spaces in \eqref{beta} below are real:
$L^2_0(\Om):=\{q\in L^2(\Om;\R):\int_\Om q\,dA=0\}$, $H^1_0(\Om)^2$ consists of real vector fields, $(\cdot,\cdot)$ and $\|\cdot\|$ are the real $L^2$
inner product and norm, and $|v|_1:=\|\nabla v\|$. The inf-sup constant of the divergence, also called the LBB constant, is
\begin{equation}\label{beta}
\beta(\Om):=\inf_{0\ne q\in L^2_0(\Om)}\ \sup_{0\ne v\in H^1_0(\Om)^2}\frac{(q,\Div v)}{|v|_1\,\|q\|}.
\end{equation}
Positivity of $\beta$ is equivalent to the well-posedness of the Stokes system \cite{BabuskaAziz}, and $1/\beta$ is the norm of the optimal right
inverse of the divergence (Lemma~\ref{lem:zperp}). The value of $\beta$ enters quantitatively in a priori and a posteriori error constants for the Stokes equations. In computer-assisted existence proofs for the
Navier--Stokes equations, the a priori error estimate of Nakao, Yamamoto and Watanabe \cite{NakaoYamamotoWatanabe1998}, on which the verification of
\cite{WatanabeYamamotoNakao1999} for a square is built, rests on the constant of Korn's inequality, which is closely related to the inf-sup constant
\cite{HorganPayne}. Because explicit values of such constants were not available, Liu, Nakao, You and Oishi \cite{LiuNakaoYouOishi} and Liu, Nakao and Oishi
\cite{LiuNakaoOishi} avoid them altogether by working with exactly divergence-free (Scott--Vogelius) finite element approximations.

\subsection*{State of the art}
The exact value of $\beta$ is known only for a few shapes such as disks, ellipses and annuli, see the survey in \cite{CCDL}.
\emph{Analytic lower bounds} exist for star-shaped planar domains: the Horgan--Payne estimate \cite{HorganPayne}, shown by Costabel and
Dauge \cite{CostabelDauge} to hold for some classes of domains (among them rectangles) and to require a correction in general. For the square it gives
$\beta\ge\sin(\pi/8)=0.38268\ldots$\,. \emph{Upper bounds} follow from the essential spectrum of the Cosserat operator: on a polygon it is the
union of the intervals $[\frac12-\frac{|\sin\omega|}{2\omega},\frac12+\frac{|\sin\omega|}{2\omega}]$ over the corner openings $\omega$ \cite{CCDL}, so that
$\beta^2\le\frac12-\max_\omega\frac{|\sin\omega|}{2\omega}$ \cite{Stoyan1999,CCDL}. \emph{Numerical approximation} is delicate because the underlying eigenvalue problem
is not compact; many natural finite element discretisations do not converge \cite{BCDG}. Gallistl \cite{Gallistl19} obtained guaranteed
upper bounds by a Rayleigh--Ritz scheme, Storn \cite{Storn} a convergent least-squares scheme, and Gallistl \cite{Gallistl21} an a posteriori
analysis that yields a lower bound under two hypotheses: that $\beta^2$ is an eigenvalue separated from the rest of the spectrum by a known gap, and
that the mesh is sufficiently fine. Neither hypothesis can presently be verified by a computation. Guaranteed inf-sup bounds for other operators have been
obtained recently \cite{ChaumontFrelet,Grassle}, but they do not cover the divergence.
To our knowledge, no method has been available that delivers a rigorous lower bound for $\beta$ on a polygon which is unconditional and
converges to the exact value.

\subsection*{Contribution}
\begin{enumerate}
\item A \emph{multiplier method} for upper bounds of the Friedrichs constant $\kappa(\Om)$, defined in \eqref{kappa}, which satisfies
$\kappa=1-2\beta^2$ (Theorem~\ref{thm:beta}); Section~\ref{sec:mult}.
For holomorphic $f$ and any Lipschitz $m$ vanishing on $\partial\Om$ one has $\int f^2=\int(1-\dzb m)f^2$ (Lemma~\ref{lem:id}). With piecewise
linear $m$ the weight $w=1-\dzb m$ of \eqref{eq:W} can be made zero along the sides and bounded near the corners by a fan constant, which tends to the sharp corner value
$|\sin\omega|/\omega$ (Lemma~\ref{lem:fan}). No information on the spectrum is needed.
\item A \emph{computable certificate} for the remaining compact part (Section~\ref{sec:cert}): local Taylor data on cells, an explicit
remainder, and a Gram-matrix bound obtained from explicit representers. Reliability is Theorem~\ref{thm:rel}.
\item A \emph{convergence theorem} (Section~\ref{sec:conv}): on the triangles treated by the multiplier alone the method keeps the full quadratic form instead of a pointwise bound, so that
the certificate loses nothing there; every $\kappa_0$ above $\kappa(\Om)$ and above the corner constant of the mesh is certified by suitable discrete data, and the certified bounds converge to $\beta$.
\item A \emph{verified implementation} in exact integer, rational and ball arithmetic (Section~\ref{sec:impl}), and
\emph{two-sided enclosures} for several polygons, including a non-convex and a doubly connected one (Section~\ref{sec:num}).
\end{enumerate}
\emph{Results.} We obtain $0.426245\le\beta\le0.426252$ for the square, $0.387220\le\beta\le0.387263$ for the rectangle $(0,2)\times(0,1)$, and
$0.303973\le\beta\le0.304078$ for the L-shaped domain. For a square with a square hole, which is neither simply connected nor star-shaped, we obtain
$0.154596\le\beta\le0.154781$. For the rectangles $(0,1)\times(0,r)$ we enclose $\beta$ for eight aspect ratios and prove for the five ratios
$r\in\{\frac14,\frac38,\frac12,\frac9{16},\frac58\}$ that $\beta^2$ is an isolated eigenvalue below the essential spectrum; for
$r\in\{\frac14,\frac38,\frac12\}$ this already follows from \cite[Sect.~5.1]{CCDL}. 
As an application, the certified values give guaranteed a posteriori bounds for Taylor--Hood approximations of the Stokes equations
(Section~\ref{sec:app}), including the doubly connected domain, for which, to our knowledge, no explicit lower bound for $\beta$ had been available.

The method is restricted to two dimensions, because it uses holomorphic functions. The enclosures are computer-assisted: all numerical
inequalities are verified by machine, while the statements of Sections~\ref{sec:fried}--\ref{sec:conv} are proved in the paper.

\subsection*{Organisation}
Section~\ref{sec:fried} recalls the Friedrichs formulation of $\beta$ and the corner bound $\kappa\ge\kess$. Section~\ref{sec:mult} proves the
multiplier identity, constructs the corner fans and reduces the bound $\kappa\le\kappa_0$ to a compact quadratic form on the active triangles.
Section~\ref{sec:cert} turns this into a finite certificate, Section~\ref{sec:conv} proves convergence, Section~\ref{sec:impl} describes the
verified implementation, Section~\ref{sec:num} reports the enclosures and Section~\ref{sec:app} the application to the Stokes equations.
Appendix~\ref{app:complex} collects facts from complex analysis, Appendix~\ref{app:aux} auxiliary results, and Appendix~\ref{app:lit} self-contained
proofs of results from the literature. The main notation of Sections~\ref{sec:fried} and~\ref{sec:mult} is listed in Table~\ref{tab:notation}.

\begin{table}[ht]
\caption{Notation of Sections~\ref{sec:fried}--\ref{sec:mult}. The letter $\beta$ is reserved for $\beta(\Om)$ (and bounds for it), $\delta$ for the
exponent of the corner trial functions, and ``Friedrichs constant'' always means $\kappa(\Om)$.}\label{tab:notation}
\centering\small
\begin{tabular}{@{}lll@{}}
\toprule
Symbol & Meaning & Defined in\\
\midrule
$\beta(\Om)$ & inf-sup constant & \eqref{beta}\\
$L^2_0(\Om)$, $|v|_1$ & real mean-zero $L^2$ functions, $\|\nabla v\|$ & Section~\ref{sec:intro}\\
$D(a,R)$, $|\Om|$ & open disk with centre $a$ and radius $R$; area of $\Om$ & ---\\
$\langle\cdot,\cdot\rangle$, $\dz$, $\dzb$ & complex $L^2$ inner product, Wirtinger derivatives & Section~\ref{sec:fried}\\
$\A$, $\Az$ & Bergman space, its mean-zero subspace & \eqref{eq:bergman}\\
$\kappa(\Om)$ & Friedrichs constant & \eqref{kappa}\\
$u_q$ & Riesz representer of $q\in L^2_0(\Om)$ & \eqref{eq:uq}\\
$\omega_j$, $\kess$ & corner openings, corner constant & \eqref{eq:kess}\\
$c'$, $M$, $g_\delta$ & corner trial functions and their ingredients & \eqref{eq:gdelta}\\
$m$, $w$ & multiplier and its weight $1-\dzb m$ & \eqref{eq:W}\\
$m_\ell$ & reflection multiplier of a side & \eqref{eq:mell}\\
$S_l$, $E_l$, $a_l$, $b_l$, $w_l$ & fan sectors, points $e^{-2i\vartheta_l}$, coefficients, weights & Definition~\ref{def:fan}, \eqref{eq:fanm}\\
$\kappa_{n,\gamma}(\omega)$ & optimal fan constant & \eqref{eq:kng}\\
$\mathcal T=\mathcal P\cup\mathcal A$, $w_T$, $\kappa_{\mathcal P}$ & passive/active split, triangle weights, passive constant & Definition~\ref{def:split}\\
$\kappa_0$ & candidate upper bound for $\kappa(\Om)$, $\kappa_0>\kappa_{\mathcal P}$ & Definition~\ref{def:split}\\
$q_T$, $\|\cdot\|_q$, $\langle\cdot,\cdot\rangle_q$ & weighted forms, norm, inner product & \eqref{q}, \eqref{eq:qip}\\
$c$, $c_{\min}$ & pointwise lower bound for $q$ & \eqref{eq:cfun}\\
\bottomrule
\end{tabular}
\end{table}

\section{The Friedrichs formulation}\label{sec:fried}
We identify $\R^2$ with $\C$ and write $z=x+iy$. On $L^2:=L^2(\Om;\C)$ we use the inner product $\langle u,v\rangle:=\int_\Om u\bar v\,dA$, and the
Wirtinger derivatives $\dz:=\frac12(\partial_x-i\partial_y)$ and $\dzb:=\frac12(\partial_x+i\partial_y)$ (Appendix~\ref{app:complex}); a $C^1$
function $f$ is holomorphic if and only if $\dzb f=0$, which is the pair of Cauchy--Riemann equations for $\Re f$ and $\Im f$. The Bergman space and
its mean-zero subspace are ($\A=\A(\Om)$)
\begin{equation}\label{eq:bergman}
\A:=\{f\in L^2:f\text{ is holomorphic in }\Om\},\qquad \Az:=\Big\{f\in\A:\int_\Om f\,dA=0\Big\}.
\end{equation}
Functions in $\A$ are single valued; $\Om$ need not be simply connected. For real vector fields $v=(v_1,v_2)$ and real functions $\varphi$ we put
$\Div v:=\partial_xv_1+\partial_yv_2$, $\rot v:=\partial_xv_2-\partial_yv_1$ and $\curl\varphi:=(\partial_y\varphi,-\partial_x\varphi)$. By Poincar\'e's
inequality $|\cdot|_1$ is a norm on $H^1_0(\Om)^2$, which is a Hilbert space with the inner product $(\nabla\cdot,\nabla\cdot)$.

The \emph{Friedrichs constant} of $\Om$ is
\begin{equation}\label{kappa}
\kappa(\Om):=\sup_{f\in\Az\setminus\{0\}}\frac{|\int_\Om f^2\,dA|}{\int_\Om|f|^2\,dA}.
\end{equation}

\begin{lemma}[phase rotation]\label{lem:phase}
Let $f\in\Az$, and let $\phi\in\R$ satisfy $e^{2i\phi}\int_\Om f^2\,dA=|\int_\Om f^2\,dA|$ (such $\phi$ exists). Then $F:=e^{i\phi}f$ satisfies $F\in\Az$,
$\|F\|=\|f\|$ and $\int_\Om F^2\,dA=|\int_\Om f^2\,dA|$. Consequently
\begin{equation}\label{eq:kappa-re}
\kappa(\Om)=\sup_{f\in\Az\setminus\{0\}}\frac{\Re\int_\Om f^2\,dA}{\int_\Om|f|^2\,dA}.
\end{equation}
\end{lemma}
\begin{proof}
If $\int f^2=|\int f^2|e^{i\psi}$, then $\phi:=-\psi/2$ is admissible. $\Az$ is invariant under multiplication by the constant $e^{i\phi}$, $|F|=|f|$, and
$\int F^2=e^{2i\phi}\int f^2=|\int f^2|$. For \eqref{eq:kappa-re}: trivially $\Re\int f^2\le|\int f^2|$, so the right-hand side of \eqref{eq:kappa-re}
is at most $\kappa(\Om)$; conversely, for $f\ne0$ the function $F$ has $\Re\int F^2/\|F\|^2=|\int f^2|/\|f\|^2$.
\end{proof}

\begin{lemma}[div--rot identity]\label{lem:divrot}
For $u,v\in H^1_0(\Om)^2$,
\begin{equation}\label{eq:beta-a}
(\nabla u,\nabla v)=(\Div u,\Div v)+(\rot u,\rot v);
\end{equation}
in particular $|v|_1^2=\|\Div v\|^2+\|\rot v\|^2$, so that $\|\Div v\|\le|v|_1$. Moreover, $\int_\Om\Div v\,dA=0$ and $\int_\Om\rot v\,dA=0$.
\end{lemma}
\begin{proof}
For $v\in C_c^\infty(\Om)^2$, expanding the squares,
\[(\Div v)^2+(\rot v)^2=|\nabla v|^2+2\big(\partial_xv_1\,\partial_yv_2-\partial_yv_1\,\partial_xv_2\big).\]
Integrating by parts first in $x$ and then in $y$ (there are no boundary terms, $v$ having compact support),
\[\int_\Om\partial_xv_1\,\partial_yv_2\,dA=-\int_\Om v_1\,\partial_x\partial_yv_2\,dA=\int_\Om\partial_yv_1\,\partial_xv_2\,dA,\]
so the mixed term has integral zero. Both sides of $|v|_1^2=\|\Div v\|^2+\|\rot v\|^2$ are continuous in the $H^1$ norm and $C_c^\infty(\Om)^2$ is dense in
$H^1_0(\Om)^2$, so this identity holds on $H^1_0(\Om)^2$. Both sides are the quadratic forms of symmetric bilinear forms; polarisation gives
\eqref{eq:beta-a}. For $v\in C_c^\infty$ each of $\int\partial_xv_1$, $\int\partial_yv_2$, $\int\partial_xv_2$, $\int\partial_yv_1$ vanishes by the
fundamental theorem of calculus on lines parallel to the axes; these integrals are continuous in $H^1$ because $\Om$ is bounded.
\end{proof}
For smooth $u$ one more integration by parts turns \eqref{eq:beta-a} into the identity $\Delta u=\nabla\Div u-\curl\rot u$.

The following reformulation of \eqref{beta} is standard (cf.\ \cite[(1.4), (1.5)]{CCDL}). 
\begin{lemma}[Riesz representer]\label{lem:riesz}
For every $q\in L^2_0(\Om)$ there is a unique $u_q\in H^1_0(\Om)^2$, the \emph{Riesz representer} of $q$, with
\begin{equation}\label{eq:uq}
(\nabla u_q,\nabla v)=(q,\Div v)\qquad(v\in H^1_0(\Om)^2).
\end{equation}
It satisfies
\begin{equation}\label{eq:beta-b}
\sup_{0\ne v\in H^1_0(\Om)^2}\frac{(q,\Div v)}{|v|_1}=|u_q|_1,\qquad |u_q|_1^2=(q,\Div u_q),\qquad \beta(\Om)=\inf_{0\ne q\in L^2_0(\Om)}\frac{|u_q|_1}{\|q\|}.
\end{equation}
\end{lemma}
\begin{proof}
By Lemma~\ref{lem:divrot}, $|(q,\Div v)|\le\|q\|\,\|\Div v\|\le\|q\|\,|v|_1$, so the Riesz representation theorem gives a unique $u_q$ with \eqref{eq:uq}.
By the Cauchy--Schwarz inequality $(q,\Div v)=(\nabla u_q,\nabla v)\le|u_q|_1|v|_1$, with equality for $v=u_q$; this gives the first two relations in
\eqref{eq:beta-b}, and the third follows from the first and \eqref{beta}.
\end{proof}

The next theorem is a reformulation of a result of Costabel and Dauge \cite[Thm.~2.1]{CostabelDauge}, 
which asserts $C(\Om)=\Gamma(\Om)+1$ for every bounded planar domain, without regularity or simple-connectivity assumptions. Here $C(\Om)=\beta(\Om)^{-2}$
and $\Gamma(\Om):=\sup_{0\ne f\in\Az}\|\Re f\|^2/\|\Im f\|^2$; since $\Re\int f^2=\|\Re f\|^2-\|\Im f\|^2$ and $\|f\|^2=\|\Re f\|^2+\|\Im f\|^2$,
\eqref{eq:kappa-re} gives $\kappa=(\Gamma-1)/(\Gamma+1)$. For smooth simply connected domains the identity is due to Horgan and Payne
\cite[Sec.~5]{HorganPayne}; the inequality of Friedrichs goes back to \cite{Friedrichs}.
\begin{theorem}\label{thm:beta}
For every bounded Lipschitz domain $\Om\subset\R^2$, $\ \beta(\Om)^2=\dfrac{1-\kappa(\Om)}{2}$.
\end{theorem}
Both directions of Theorem~\ref{thm:beta} enter the verified bounds. A self-contained proof, included for completeness, is given in Appendix~\ref{app:lit-beta}.
The facts from complex analysis used here and below (Wirtinger derivatives, Weyl's lemma, Bergman spaces) are collected, with proofs or references, in
Appendix~\ref{app:complex}.

\begin{corollary}\label{cor:bounds}
\begin{enumerate}
\item[(a)] If $\kappa(\Om)\le\kappa_0$ for some $\kappa_0<1$, then $\beta(\Om)\ge\sqrt{(1-\kappa_0)/2}$.
\item[(b)] For every $f\in\Az\setminus\{0\}$, $\beta(\Om)^2\le\frac12\big(1-|\int_\Om f^2\,dA|/\|f\|^2\big)$.
\end{enumerate}
\end{corollary}
\begin{proof}
(a) follows from Theorem~\ref{thm:beta}. (b) follows from Theorem~\ref{thm:beta} and $|\int f^2|/\|f\|^2\le\kappa(\Om)$, by \eqref{kappa}.
\end{proof}

Trial functions need not have mean zero:
\begin{lemma}[mean correction]\label{lem:mean}
Let $g\in\A$ and $I:=\int_\Om g\,dA$. Then $f:=g-I/|\Om|$ belongs to $\Az$, and
\[\int_\Om f^2\,dA=\int_\Om g^2\,dA-\frac{I^2}{|\Om|},\qquad\|f\|^2=\|g\|^2-\frac{|I|^2}{|\Om|}.\]
\end{lemma}
\begin{proof}
Constants belong to $\A$ because $\Om$ is bounded, and $\int f=I-I=0$. With $\mu:=I/|\Om|$, expanding,
$\int f^2=\int g^2-2\mu I+\mu^2|\Om|=\int g^2-I^2/|\Om|$ and $\|f\|^2=\|g\|^2-2\Re(\bar\mu I)+|\mu|^2|\Om|=\|g\|^2-|I|^2/|\Om|$.
\end{proof}

\begin{definition}[polygonal domain]\label{def:polygon}
A \emph{polygonal domain} is a bounded, connected Lipschitz domain $\Om\subset\R^2$ whose boundary consists of finitely many straight sides; $\Om$ may be
multiply connected. Its \emph{corners} $c_1,\dots,c_J$ are the boundary vertices whose interior openings $\omega_j$ lie in $(0,2\pi)\setminus\{\pi\}$,
and its \emph{corner constant} is
\begin{equation}\label{eq:kess}
\kess:=\max_{1\le j\le J}\frac{|\sin\omega_j|}{\omega_j}.
\end{equation}
\end{definition}
From now on $\Om$ is a polygonal domain.

\emph{Corner trial functions.} Let $c$ be a corner of $\Om$ with opening $\omega$, and let $\theta_0\in\R$ and $R_0>0$ be such that
$\Om\cap D(c,R_0)=S:=\{c+re^{i\theta}:0<r<R_0,\ \theta_0<\theta<\theta_0+\omega\}$. Fix $t\in(0,R_0)$. For $0<\delta<1$ put
\begin{equation}\label{eq:gdelta}
c':=c-te^{i(\theta_0+\omega/2)},\qquad M(z):=\frac{z-c}{z-c'},\qquad g_\delta(z):=\exp\big((-1+\delta)\operatorname{Log}M(z)\big),
\end{equation}
where $\operatorname{Log}$ is the principal logarithm on $\C\setminus(-\infty,0]$ (Lemma~\ref{lem:holo}(e)); that is, $g_\delta=M^{-1+\delta}$ with the
principal power. The point $c'$ lies on the exterior bisector of the corner.
\begin{lemma}\label{lem:gdelta}
$[c,c']\cap\overline\Om=\{c\}$, and for every $\delta\in(0,1)$ the function $g_\delta$ of \eqref{eq:gdelta} is single valued and holomorphic in $\Om$,
$|g_\delta(z)|=|z-c|^{-1+\delta}|z-c'|^{1-\delta}$, and $g_\delta\in\A$.
\end{lemma}
\begin{proof}
The points $c+s(c'-c)$, $0<s\le1$, lie in $D(c,R_0)$ at the angle $\theta_0+\omega/2+\pi$, which is not in $[\theta_0,\theta_0+\omega]$ modulo $2\pi$
because $\omega<2\pi$; hence $[c,c']\cap\overline\Om=\{c\}$. The M\"obius map $M$ is a bijection of the Riemann sphere $\C\cup\{\infty\}$. For
$z=c+s(c'-c)$ one finds $M(z)=s/(s-1)$, so $M$ maps the segment $[c,c']$ onto $[-\infty,0]$ and hence $\C\setminus[c,c']$ onto $\C\setminus(-\infty,0]$.
There the principal power $\zeta^{-1+\delta}=\exp((-1+\delta)\operatorname{Log}\zeta)$ is holomorphic (Lemma~\ref{lem:holo}(e)). Since $\Om\cap[c,c']=\emptyset$,
$g_\delta=M^{-1+\delta}$ is holomorphic in $\Om$, and it is a genuine single-valued function there whatever the topology of $\Om$: its only branch cut
is the segment $[c,c']$, which lies outside $\Om$. Further $|g_\delta(z)|=|M(z)|^{-1+\delta}=|z-c|^{-1+\delta}|z-c'|^{1-\delta}$, which is bounded on
$\Om\setminus D(c,\rho)$ for every $\rho>0$, while
$\int_S|z-c|^{-2+2\delta}dA=\int_{\theta_0}^{\theta_0+\omega}\int_0^{R_0}r^{-1+2\delta}\,dr\,d\theta=\omega R_0^{2\delta}/(2\delta)<\infty$ and
$|z-c'|$ is bounded on $\Om$; hence $g_\delta\in\A$.
\end{proof}

By Theorem~\ref{thm:beta}, the next lemma is equivalent to $\beta(\Om)^2\le\frac12-\frac{\kess}2$, which is \cite[Thm.~3.3, eq.~(3.11)]{CCDL} 
(where the interval in \cite[(3.10)]{CCDL} is to be read with $|\sin\omega_c|$ at re-entrant corners); Costabel, Crouzeix, Dauge and Lafranche
attribute this bound to Stoyan \cite{Stoyan1999}.
\begin{lemma}\label{lem:kess}$\kappa(\Om)\ge\kess$.\end{lemma}
A self-contained proof, included for completeness, is given in Appendix~\ref{app:lit-kess}: suitably mean-corrected functions $g_\delta$ of
\eqref{eq:gdelta} have Rayleigh quotients $|\int f^2|/\|f\|^2$ tending to $|\sin\omega|/\omega$ as $\delta\to0$. These trial functions are used again
in Section~\ref{sec:impl}.

\begin{remark}\label{rem:noncompact}
The quadratic form $f\mapsto\int_\Om f^2$ is not compact on $\A$ when $\Om$ has a corner, which is the source of the difficulty of approximating $\kappa$
from above. Indeed, the functions $f_\delta$ of the proof of Lemma~\ref{lem:kess} (Appendix~\ref{app:lit-kess}) satisfy $\|f_\delta\|\to\infty$ and
$\|f_\delta\|_{L^2(\Om\setminus D(c,r))}=O(1)$ as $\delta\to0$, for every $r>0$. Hence $h_\delta:=f_\delta/\|f_\delta\|$ converges weakly to $0$ in $\A$
(split $\langle h,h_\delta\rangle$ into the integrals over $\Om\cap D(c,r)$ and $\Om\setminus D(c,r)$), while $|\int h_\delta^2|\to|\sin\omega|/\omega>0$; a compact form would give $\int h_\delta^2\to0$.
\end{remark}

\section{The multiplier identity}\label{sec:mult}
\emph{Standing notation.} In Sections~\ref{sec:mult}--\ref{sec:num} a \emph{multiplier} is a Lipschitz function $m\colon\overline\Om\to\C$ with $m=0$ on
$\partial\Om$, and the \emph{weight} of $m$ is
\begin{equation}\tag{W}\label{eq:W}
w:=1-\dzb m ,
\end{equation}
a bounded function defined almost everywhere (by Rademacher's theorem \cite{EvansGariepy}). The same definition is used for a Lipschitz function $m$
given only on a part of $\overline\Om$, such as the closed sector at a corner in Definition~\ref{def:fan}; then $w$ is defined on that part. Whenever a
multiplier is fixed, $w$ denotes its weight \eqref{eq:W}. If $m$ is affine on a set, $w$ is constant there (Appendix~\ref{app:complex}); the constant
values are $w_T$ on a triangle $T$ and $w_l$ on a sector $S_l$.

Weighted identities for $\int_\Om f^2$ obtained by integration by parts against a weight vanishing on $\partial\Om$ appear in Horgan and Payne's estimate
for star-shaped domains \cite[\S6]{HorganPayne} (see \cite[proof of Thm.~4.4, eq.~(4.10)]{CostabelDauge}); 
the general form below, for arbitrary Lipschitz multipliers and without regularity of $\Om$ or $f$ at the boundary, is elementary, and we include its proof.
\begin{lemma}[integration by parts]\label{lem:ibp}
Let $F\in L^1(\Om)$ be holomorphic in $\Om$, and let $u$ be Lipschitz on $\overline\Om$ with $u|_{\partial\Om}=0$. Then
\begin{equation}\tag{IBP}\label{eq:ibp}
\int_\Om F\,\dzb u\,dA=0 .
\end{equation}
\end{lemma}
\begin{proof}
A Lipschitz function is differentiable almost everywhere, with bounded derivatives that are its weak derivatives (Rademacher's theorem \cite{EvansGariepy});
if $L$ is a Lipschitz constant, $|\dzb u|\le\frac12(|\partial_xu|+|\partial_yu|)\le L$ a.e. Hence all integrals below are finite.

\emph{Step 1: compactly supported $u$.} Let $\phi$ be Lipschitz with compact support $K\subset\Om$, and $\epsilon_0:=\frac12\dist(K,\partial\Om)>0$. Let $j_\epsilon\ge0$ be a
smooth mollifier supported in $D(0,\epsilon)$ with $\int j_\epsilon=1$. For $\epsilon<\epsilon_0$, $\phi_\epsilon:=\phi*j_\epsilon$ is smooth and supported in the compact set
$K':=\{\dist(\cdot,K)\le\epsilon_0\}\subset\Om$. Since $F$ is $C^\infty$ with $\dzb F=0$ in $\Om$, classical integration by parts (Fubini and the fundamental theorem of
calculus in $x$ and in $y$) gives $\int F\dzb\phi_\epsilon=-\int\phi_\epsilon\dzb F=0$. Moreover $\dzb\phi_\epsilon=(\dzb\phi)*j_\epsilon\to\dzb\phi$ in $L^1$, and $F$ is
bounded on $K'$; hence $\int F\dzb\phi=\lim_{\epsilon\to0}\int F\dzb\phi_\epsilon=0$.

\emph{Step 2: cut-off functions.} For $k\in\mathbb N$ let $U_k:=\{z\in\Om:\dist(z,\partial\Om)>3/(2k)\}$ and $\eta_k:=\mathbf 1_{U_k}*j_{1/(2k)}$. Then $\eta_k$ is smooth and
$0\le\eta_k\le1$; $\eta_k(z)=1$ if $\dist(z,\partial\Om)>2/k$, because then $D(z,1/(2k))\subset U_k$; and $\eta_k(z)=0$ unless $z\in\Om$ and $\dist(z,\partial\Om)>1/k$, because
otherwise $D(z,1/(2k))\cap U_k=\emptyset$. So $\eta_k\in C_c^\infty(\Om)$, $\nabla\eta_k=0$ outside the boundary layer $\Sigma_k:=\{z\in\Om:\dist(z,\partial\Om)\le2/k\}$, and
$|\nabla\eta_k|\le\|\nabla j_{1/(2k)}\|_{L^1}=2k\|\nabla j_1\|_{L^1}=:Ck$.

\emph{Step 3: passage to the limit.} The function $u\eta_k$ is Lipschitz with compact support in $\Om$, so Step~1 and the product rule
$\dzb(u\eta_k)=\eta_k\dzb u+u\dzb\eta_k$ (a.e.) give
\begin{equation*}\tag{$*$}
\int_\Om F\,\eta_k\,\dzb u\,dA=-\int_\Om F\,u\,\dzb\eta_k\,dA .
\end{equation*}
On the left, $|F\eta_k\dzb u|\le L|F|\in L^1$ and $\eta_k\to1$ pointwise in $\Om$, so by dominated convergence the left side tends to $\int F\dzb u$.
On the right, for $z\in\Om$ let $z^*\in\partial\Om$ with $|z-z^*|=\dist(z,\partial\Om)$; since $u(z^*)=0$, $|u(z)|\le L\dist(z,\partial\Om)$. (For a continuous piecewise
linear $u$ it suffices that $L\ge\max|\nabla u|$, since the segment $[z,z^*)$ lies in $\Om$.) On $\Sigma_k$ this gives $|u|\le2L/k$, while
$|\dzb\eta_k|=\frac12|\nabla\eta_k|\le Ck/2$. Hence
\[\Big|\int_\Om Fu\,\dzb\eta_k\,dA\Big|\le\int_{\Sigma_k}|F|\cdot L\dist(\cdot,\partial\Om)\cdot\frac{Ck}{2}\,dA\le CL\int_{\Sigma_k}|F|\,dA\longrightarrow0\qquad(k\to\infty),\]
because $\mathbf 1_{\Sigma_k}\to0$ pointwise in $\Om$ and $|F|\in L^1$ (dominated convergence). The factor $k\dist(\cdot,\partial\Om)\le2$ is bounded on the layer: the
vanishing of $u$ on $\partial\Om$ exactly compensates the size $O(k)$ of $\nabla\eta_k$. No Hardy-type inequality, and no information on the boundary behaviour of $F$
beyond $F\in L^1$, is needed. Letting $k\to\infty$ in $(*)$ proves \eqref{eq:ibp}.
\end{proof}

\begin{lemma}[multiplier identity]\label{lem:id}
Let $m$ be a multiplier and $w$ its weight \eqref{eq:W}. Then $\int_\Om f^2\,dA=\int_\Om wf^2\,dA$ for all $f\in\A$.
\end{lemma}
\begin{proof}
$F:=f^2$ is holomorphic and $\int|F|=\|f\|^2<\infty$, so Lemma~\ref{lem:ibp} with $u=m$ gives $\int f^2\dzb m=0$. By \eqref{eq:W},
$\int wf^2=\int f^2-\int f^2\dzb m=\int f^2$.
\end{proof}
The lemma replaces the weight $1$ in $\int_\Om f^2$ by the weight $w$ of any multiplier. Since $\int_\Om w=|\Om|$ (take $f=1$ in Lemma~\ref{lem:id}), $w$ cannot be
small everywhere; the idea is to make it small \emph{near the boundary}, where functions of $\A$ are not controlled, and to accept large values in the interior,
where they are.

\begin{definition}[reflection multiplier]\label{def:refl}
Let $\ell$ be the line containing a side of $\Om$, $z_0\in\ell$, and $\nu$ the inward unit normal of the side, as a complex number. The \emph{signed
distance} to $\ell$ and the \emph{reflection multiplier} of the side are
\begin{equation}\label{eq:mell}
\dist(z,\ell):=\Re\big((z-z_0)\bar\nu\big),\qquad m_\ell(z):=2\dist(z,\ell)\,\bar\nu .
\end{equation}
\end{definition}
\begin{lemma}\label{lem:refl}
$\dist(z,\ell)$ does not depend on the choice of $z_0\in\ell$ and is positive on the side of $\ell$ that contains $\Om$ near the side. Moreover
$m_\ell(z)=(z-z_0)\bar\nu^2+(\bar z-\bar z_0)$; $m_\ell$ is affine, vanishes on $\ell$, and has $\dzb m_\ell=1$, i.e.\ weight $w=0$; and
$z-\overline{m_\ell(z)}=z-2\dist(z,\ell)\nu$ is the mirror image of $z$ in $\ell$.
\end{lemma}
\begin{proof}
For $z_0,z_1\in\ell$ the number $z_1-z_0$ is parallel to the side, hence orthogonal to $\nu$, so $\Re((z_1-z_0)\bar\nu)=0$; and
$\Re((z_0+s\nu-z_0)\bar\nu)=s>0$ for $s>0$. Writing $\dist(z,\ell)=\frac12\big((z-z_0)\bar\nu+(\bar z-\bar z_0)\nu\big)$ and using $|\nu|=1$ gives
the formula for $m_\ell$, which is affine and vanishes where $\dist(z,\ell)=0$. Since $\dzb z=0$ and $\dzb\bar z=1$, $\dzb m_\ell=1$. Finally
$\overline{m_\ell(z)}=2\dist(z,\ell)\nu$ because $\dist(z,\ell)$ is real.
\end{proof}
The reflection multiplier does not vanish on the other sides, so it is not itself a multiplier; it is used near its side, as a piece of multipliers
assembled from several formulas (Lemma~\ref{lem:fan}(c), Lemma~\ref{lem:goodmult}).

\emph{Corners.} Near a corner the reflection multipliers of the two sides meet, and they cannot be joined with $w=0$: by the fan identity
\eqref{fan} below, a weight that vanishes on the whole corner sector would force $\sin\omega=0$. The optimal piecewise linear construction is explicit;
Figure~\ref{fig:fan} shows its ingredients for a right angle.
\begin{definition}[fan multiplier]\label{def:fan}
Let the corner be at $0$ and $\Om\cap D(0,R_0)=\{re^{i\theta}:0<r<R_0,\ \theta_0<\theta<\theta_0+\omega\}$, and put
$\overline S:=\{re^{i\theta}:0\le r<R_0,\ \theta_0\le\theta\le\theta_0+\omega\}$. Let $\theta_0=\vartheta_0<\vartheta_1<\dots<\vartheta_N=\theta_0+\omega$ with
$0<\vartheta_l-\vartheta_{l-1}<\pi$ for every $l$, and put
\[S_l:=\{re^{i\theta}:0<r<R_0,\ \vartheta_{l-1}<\theta<\vartheta_l\}\quad(l=1,\dots,N),\qquad E_l:=e^{-2i\vartheta_l}\quad(l=0,\dots,N),\]
$\Delta E_l:=E_l-E_{l-1}$ and $\Delta E:=E_N-E_0$. A \emph{fan multiplier} for this partition is a continuous function $m\colon\overline S\to\C$ that
vanishes on both sides $\theta=\theta_0$ and $\theta=\theta_0+\omega$ and is $\R$-linear on each sector: for every $l$ there are $a_l,b_l\in\C$ with
$m(z)=a_lz+b_l\bar z$ on $S_l$. Since $\dzb(a_lz+b_l\bar z)=b_l$, the weight \eqref{eq:W} of a fan multiplier is constant on each sector, and
\begin{equation}\label{eq:fanm}
m(z)=a_lz+b_l\bar z\quad(z\in S_l),\qquad w_l:=w|_{S_l}=1-b_l\qquad(l=1,\dots,N).
\end{equation}
\end{definition}
The hypothesis $0<\vartheta_l-\vartheta_{l-1}<\pi$ is used in two places. It gives $|\Delta E_l|=2\sin(\vartheta_l-\vartheta_{l-1})>0$; for a sector
opening $\vartheta_l-\vartheta_{l-1}\in(\pi,2\pi)$ the sine is negative (e.g.\ $N=1$, $\omega=\frac{3\pi}2$), and for a sector opening $\pi$ one has $E_l=E_{l-1}$, so that the
sector does not enter \eqref{fan}. And the extremal weights in Lemma~\ref{lem:fan}(a) divide by $|\Delta E_l|$. The hypothesis holds for every partition if
$\omega<\pi$; at a re-entrant corner it requires at least two sectors.

\begin{lemma}[fan identity]\label{lem:fanid}
For a partition as in Definition~\ref{def:fan}, every fan multiplier satisfies
\begin{equation}\label{fan}\sum_{l=1}^Nw_l(E_l-E_{l-1})=E_N-E_0,\qquad |E_N-E_0|=2|\sin\omega|,\quad|E_l-E_{l-1}|=2\sin(\vartheta_l-\vartheta_{l-1}),\end{equation}
and conversely, for every $(w_1,\dots,w_N)\in\C^N$ satisfying the first identity in \eqref{fan} there is exactly one fan multiplier with these weights.
\end{lemma}
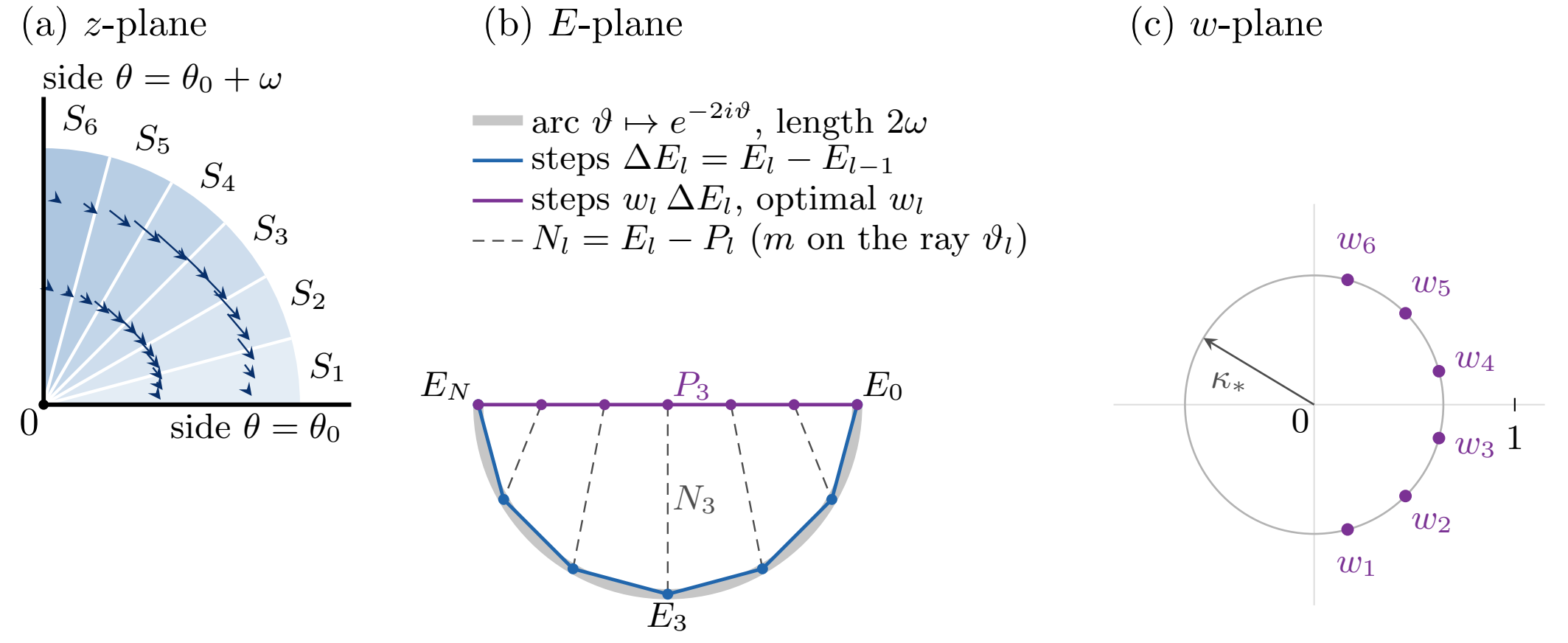
\begin{figure}[t]
\centering
\begin{tikzpicture}[>=stealth,font=\small]
\definecolor{fanblue}{RGB}{33,102,172}\definecolor{wpurple}{RGB}{123,50,148}\definecolor{ngrey}{RGB}{77,77,77}\definecolor{mblue}{RGB}{8,48,107}
\begin{scope}
\fill[fanblue!12] (0,0) -- (0.00:2.3) arc[start angle=0.00,end angle=15.00,radius=2.3] -- cycle;
\fill[fanblue!17] (0,0) -- (15.00:2.3) arc[start angle=15.00,end angle=30.00,radius=2.3] -- cycle;
\fill[fanblue!22] (0,0) -- (30.00:2.3) arc[start angle=30.00,end angle=45.00,radius=2.3] -- cycle;
\fill[fanblue!27] (0,0) -- (45.00:2.3) arc[start angle=45.00,end angle=60.00,radius=2.3] -- cycle;
\fill[fanblue!32] (0,0) -- (60.00:2.3) arc[start angle=60.00,end angle=75.00,radius=2.3] -- cycle;
\fill[fanblue!37] (0,0) -- (75.00:2.3) arc[start angle=75.00,end angle=90.00,radius=2.3] -- cycle;
\draw[white,line width=0.8pt] (0,0) -- (15.00:2.3);
\draw[white,line width=0.8pt] (0,0) -- (30.00:2.3);
\draw[white,line width=0.8pt] (0,0) -- (45.00:2.3);
\draw[white,line width=0.8pt] (0,0) -- (60.00:2.3);
\draw[white,line width=0.8pt] (0,0) -- (75.00:2.3);
\node[font=\footnotesize] at (7.50:2.576) {$S_1$};
\node[font=\footnotesize] at (22.50:2.576) {$S_2$};
\node[font=\footnotesize] at (37.50:2.576) {$S_3$};
\node[font=\footnotesize] at (52.50:2.576) {$S_4$};
\node[font=\footnotesize] at (67.50:2.576) {$S_5$};
\node[font=\footnotesize] at (82.50:2.576) {$S_6$};
\draw[->,mblue,line width=0.45pt] (1.0328,0.0677) -- (1.0505,0.0440);
\draw[->,mblue,line width=0.45pt] (1.0151,0.2019) -- (1.0679,0.1312);
\draw[->,mblue,line width=0.45pt] (0.9801,0.3327) -- (1.0645,0.2254);
\draw[->,mblue,line width=0.45pt] (0.9283,0.4578) -- (1.0405,0.3249);
\draw[->,mblue,line width=0.45pt] (0.8606,0.5750) -- (0.9940,0.4270);
\draw[->,mblue,line width=0.45pt] (0.7782,0.6824) -- (0.9257,0.5300);
\draw[->,mblue,line width=0.45pt] (0.6824,0.7782) -- (0.8349,0.6306);
\draw[->,mblue,line width=0.45pt] (0.5750,0.8606) -- (0.7231,0.7271);
\draw[->,mblue,line width=0.45pt] (0.4578,0.9283) -- (0.5907,0.8160);
\draw[->,mblue,line width=0.45pt] (0.3327,0.9801) -- (0.4400,0.8956);
\draw[->,mblue,line width=0.45pt] (0.2019,1.0151) -- (0.2726,0.9624);
\draw[->,mblue,line width=0.45pt] (0.0677,1.0328) -- (0.0914,1.0151);
\draw[->,mblue,line width=0.45pt] (1.8361,0.1203) -- (1.8675,0.0782);
\draw[->,mblue,line width=0.45pt] (1.8046,0.3590) -- (1.8984,0.2333);
\draw[->,mblue,line width=0.45pt] (1.7424,0.5914) -- (1.8925,0.4007);
\draw[->,mblue,line width=0.45pt] (1.6502,0.8138) -- (1.8498,0.5775);
\draw[->,mblue,line width=0.45pt] (1.5299,1.0222) -- (1.7671,0.7590);
\draw[->,mblue,line width=0.45pt] (1.3834,1.2132) -- (1.6457,0.9422);
\draw[->,mblue,line width=0.45pt] (1.2132,1.3834) -- (1.4842,1.1211);
\draw[->,mblue,line width=0.45pt] (1.0222,1.5299) -- (1.2855,1.2927);
\draw[->,mblue,line width=0.45pt] (0.8138,1.6502) -- (1.0501,1.4506);
\draw[->,mblue,line width=0.45pt] (0.5914,1.7424) -- (0.7822,1.5922);
\draw[->,mblue,line width=0.45pt] (0.3590,1.8046) -- (0.4846,1.7109);
\draw[->,mblue,line width=0.45pt] (0.1203,1.8361) -- (0.1625,1.8046);
\draw[line width=1.1pt] (2.760,0) -- (0,0) -- (0,2.760);
\fill (0,0) circle (1.3pt) node[below left,inner sep=1pt] {$0$};
\node[below left,font=\footnotesize,inner sep=2pt] at (2.760,0) {side $\theta=\theta_0$};
\node[anchor=south west,font=\footnotesize,inner sep=1pt] at (-0.05,2.76) {side $\theta=\theta_0+\omega$};
\node[anchor=west] at (-0.35,3.4) {(a) $z$-plane};
\end{scope}
\begin{scope}[xshift=5.6cm]
\draw[gray!45,line width=2.6pt,line cap=round] (1.7,0) arc[start angle=0,end angle=-180,radius=1.7];
\draw[ngrey,densely dashed,line width=0.45pt] (1.1333,-0.0000) -- (1.4722,-0.8500);
\draw[ngrey,densely dashed,line width=0.45pt] (0.5667,-0.0000) -- (0.8500,-1.4722);
\draw[ngrey,densely dashed,line width=0.45pt] (-0.0000,-0.0000) -- (0.0000,-1.7000);
\draw[ngrey,densely dashed,line width=0.45pt] (-0.5667,-0.0000) -- (-0.8500,-1.4722);
\draw[ngrey,densely dashed,line width=0.45pt] (-1.1333,-0.0000) -- (-1.4722,-0.8500);
\draw[fanblue,line width=0.9pt] (1.7000,-0.0000) -- (1.4722,-0.8500) -- (0.8500,-1.4722) -- (0.0000,-1.7000) -- (-0.8500,-1.4722) -- (-1.4722,-0.8500) -- (-1.7000,-0.0000);
\draw[wpurple,line width=0.9pt] (1.7000,0.0000) -- (1.1333,-0.0000) -- (0.5667,-0.0000) -- (-0.0000,-0.0000) -- (-0.5667,-0.0000) -- (-1.1333,-0.0000) -- (-1.7000,-0.0000);
\fill[fanblue] (1.7000,-0.0000) circle (1.4pt);
\fill[fanblue] (1.4722,-0.8500) circle (1.4pt);
\fill[fanblue] (0.8500,-1.4722) circle (1.4pt);
\fill[fanblue] (0.0000,-1.7000) circle (1.4pt);
\fill[fanblue] (-0.8500,-1.4722) circle (1.4pt);
\fill[fanblue] (-1.4722,-0.8500) circle (1.4pt);
\fill[fanblue] (-1.7000,-0.0000) circle (1.4pt);
\fill[wpurple] (1.7000,0.0000) circle (1.4pt);
\fill[wpurple] (1.1333,-0.0000) circle (1.4pt);
\fill[wpurple] (0.5667,-0.0000) circle (1.4pt);
\fill[wpurple] (-0.0000,-0.0000) circle (1.4pt);
\fill[wpurple] (-0.5667,-0.0000) circle (1.4pt);
\fill[wpurple] (-1.1333,-0.0000) circle (1.4pt);
\fill[wpurple] (-1.7000,-0.0000) circle (1.4pt);
\node[above right,inner sep=1pt,font=\footnotesize] at (1.7000,-0.0000) {$E_0$};
\node[above left,inner sep=1pt,font=\footnotesize] at (-1.7000,-0.0000) {$E_N$};
\node[below,inner sep=2pt,font=\footnotesize] at (0.0000,-1.7000) {$E_3$};
\node[above right,inner sep=1pt,wpurple,font=\footnotesize] at (-0.0000,-0.0000) {$P_3$};
\node[right,ngrey,inner sep=1pt,font=\footnotesize] at (-0.0000,-0.8500) {$N_3$};
\draw[gray!45,line width=2.6pt] (-1.750,2.550) -- (-1.300,2.550); \node[right,font=\footnotesize,inner sep=2pt] at (-1.300,2.550) {arc $\vartheta\mapsto e^{-2i\vartheta}$, length $2\omega$};
\draw[fanblue,line width=0.9pt] (-1.750,2.190) -- (-1.300,2.190); \node[right,font=\footnotesize,inner sep=2pt] at (-1.300,2.190) {steps $\Delta E_l=E_l-E_{l-1}$};
\draw[wpurple,line width=0.9pt] (-1.750,1.830) -- (-1.300,1.830); \node[right,font=\footnotesize,inner sep=2pt] at (-1.300,1.830) {steps $w_l\,\Delta E_l$, optimal $w_l$};
\draw[ngrey,densely dashed,line width=0.45pt] (-1.750,1.470) -- (-1.300,1.470); \node[right,font=\footnotesize,inner sep=2pt] at (-1.300,1.470) {$N_l=E_l-P_l$ ($m$ on the ray $\vartheta_l$)};
\node[anchor=west] at (-1.800,3.4) {(b) $E$-plane};
\end{scope}
\begin{scope}[xshift=11.4cm]
\draw[gray!25,line width=0.4pt] (-1.800,0) -- (2.070,0) (0,-1.800) -- (0,1.800);
\draw[gray!60,line width=0.5pt] (0,0) circle (1.1591);
\fill[wpurple] (0.3000,-1.1196) circle (1.6pt); \node[wpurple,font=\footnotesize] at (0.3885,-1.4500) {$w_1$};
\fill[wpurple] (0.8196,-0.8196) circle (1.6pt); \node[wpurple,font=\footnotesize] at (1.0614,-1.0614) {$w_2$};
\fill[wpurple] (1.1196,-0.3000) circle (1.6pt); \node[wpurple,font=\footnotesize] at (1.4500,-0.3885) {$w_3$};
\fill[wpurple] (1.1196,0.3000) circle (1.6pt); \node[wpurple,font=\footnotesize] at (1.4500,0.3885) {$w_4$};
\fill[wpurple] (0.8196,0.8196) circle (1.6pt); \node[wpurple,font=\footnotesize] at (1.0614,1.0614) {$w_5$};
\fill[wpurple] (0.3000,1.1196) circle (1.6pt); \node[wpurple,font=\footnotesize] at (0.3885,1.4500) {$w_6$};
\draw[->,ngrey,line width=0.5pt] (0,0) -- (-0.9932,0.5975) node[pos=0.55,below left,inner sep=1pt,font=\footnotesize] {$\kappa_*$};
\node[below left,inner sep=1pt,font=\footnotesize] at (0,0) {$0$};
\draw (1.8,0.06) -- (1.8,-0.06) node[below,font=\footnotesize] {$1$};
\node[anchor=west] at (-1.800,3.4) {(c) $w$-plane};
\end{scope}
\end{tikzpicture}
\caption{Definition~\ref{def:fan} and Lemmas~\ref{lem:fanid}, \ref{lem:fan} for a right angle ($\theta_0=0$, $\omega=\pi/2$) with $N=n=6$ equal sectors
and the optimal weights of Lemma~\ref{lem:fan}(a). (a) The sectors $S_l$ in the $z$-plane and the vectors $m(z)$ of the optimal fan multiplier, scaled by $0.21$; they
vanish on both sides and at the corner. (b) The points $E_l=e^{-2i\vartheta_l}$ lie on an arc of length $2\omega$ of the unit circle. The chord path
with steps $\Delta E_l$ has length $2n\sin(\omega/n)=3.1058$, the chord $\Delta E=E_N-E_0$ has length $2|\sin\omega|=2$. The path
$P_l:=E_0+\sum_{k\le l}w_k\Delta E_k$ with steps $w_l\Delta E_l$ also ends at $E_N$; this is the identity \eqref{fan}. For the optimal weights it runs
along the chord. The dashed segments $E_l-P_l$ are the ray values $N_l$ of $m$ (proof of Lemma~\ref{lem:fanid}). (c) The optimal weights all have
modulus $\kappa_*=|\sin\omega|/(n\sin(\omega/n))=0.6440$, larger than the corner value $|\sin\omega|/\omega=2/\pi=0.6366$.}\label{fig:fan}
\end{figure}
\begin{proof}
\emph{Values on the rays.} An $\R$-linear function $\mu(z)=az+b\bar z$ is positively homogeneous, $\mu(re^{i\vartheta})=r\mu(e^{i\vartheta})$ for $r>0$,
so on a ray it is determined by its value at any one point of that ray. Fix $\rho:=R_0/2$ and put
\[N_l:=\rho^{-1}e^{-i\vartheta_l}\,m(\rho e^{i\vartheta_l})\qquad(l=0,\dots,N),\]
the value of $m$ at the point $\rho e^{i\vartheta_l}$ of the ray $\theta=\vartheta_l$, rotated back and rescaled. On $S_l$,
$\rho^{-1}e^{-i\vartheta}m(\rho e^{i\vartheta})=a_l+b_le^{-2i\vartheta}$. Letting $\vartheta$ tend to the two bounding angles $\vartheta_{l-1}$, $\vartheta_l$
and using the continuity of $m$ on $\overline S$,
\[N_l=a_l+b_lE_l,\qquad N_{l-1}=a_l+b_lE_{l-1}\qquad(l=1,\dots,N).\]
Since $N_l$ is defined from the values of $m$ itself, the same number $N_l$ arises from sector $l$ and from sector $l+1$; this is where continuity enters.

\emph{Recursion and side conditions.} Subtracting the two relations and using $b_l=1-w_l$,
\begin{equation*}\tag{R}\label{eq:R}
N_l-N_{l-1}=b_l(E_l-E_{l-1})=(1-w_l)(E_l-E_{l-1}).
\end{equation*}
By homogeneity, $m=0$ on the side $\theta=\vartheta_0$ if and only if $N_0=0$, and $m=0$ on the side $\theta=\vartheta_N$ if and only if $N_N=0$.
(With $N_0=0$, \eqref{eq:R} gives $N_l=E_l-P_l$ for the path $P_l$ of Figure~\ref{fig:fan}(b), and $N_l=\rho^{-1}e^{-i\vartheta_l}m(\rho e^{i\vartheta_l})$
shows $m$ on the rays; $N_0=N_N=0$ are the side conditions.)

\emph{Sum identity.} Summing \eqref{eq:R} over $l=1,\dots,N$, the left side telescopes to $N_N-N_0=0$, and the right side equals
$(E_N-E_0)-\sum_lw_l(E_l-E_{l-1})$. This is the first identity in \eqref{fan}.

\emph{Moduli.} For real $\psi,\phi$, $|e^{-2i(\psi+\phi)}-e^{-2i\psi}|=|e^{-i\phi}-e^{i\phi}|=2|\sin\phi|$. With $\phi=\omega$ this gives $|E_N-E_0|=2|\sin\omega|$, and with
$\phi=\vartheta_l-\vartheta_{l-1}\in(0,\pi)$ it gives $|E_l-E_{l-1}|=2\sin(\vartheta_l-\vartheta_{l-1})$.

\emph{Converse.} Let $(w_l)$ satisfy the first identity in \eqref{fan}. Put $b_l:=1-w_l$, $N_0:=0$, $N_l:=N_{l-1}+b_l(E_l-E_{l-1})$ for $l=1,\dots,N$
(that is, \eqref{eq:R}), and $a_l:=N_l-b_lE_l$. Then $N_N=\sum_l(1-w_l)(E_l-E_{l-1})=(E_N-E_0)-\sum_lw_l(E_l-E_{l-1})=0$ and
$a_l+b_lE_{l-1}=N_l-b_l(E_l-E_{l-1})=N_{l-1}$. So the formula $a_lz+b_l\bar z$ of sector $l$ takes the value $\rho e^{i\vartheta_{l-1}}N_{l-1}$ at
$\rho e^{i\vartheta_{l-1}}$ and $\rho e^{i\vartheta_l}N_l$ at $\rho e^{i\vartheta_l}$. Adjacent formulas therefore agree at one point of their common ray, hence by
homogeneity on the whole ray; so they define a continuous function $m$ on $\overline S$, which vanishes on both sides because $N_0=N_N=0$. This $m$ is a fan
multiplier with weights $1-b_l=w_l$. Uniqueness: for any fan multiplier with weights $(w_l)$ necessarily $b_l=1-w_l$, the $N_l$ are determined by
$N_0=0$ and \eqref{eq:R}, and $a_l=N_l-b_lE_l$.
\end{proof}

\begin{lemma}[corner fan]\label{lem:fan}
Consider a corner and a partition as in Definition~\ref{def:fan}.
\begin{enumerate}
\item[(a)] Every fan multiplier satisfies $\max_l|w_l|\ge|\sin\omega|/\sum_l\sin(\vartheta_l-\vartheta_{l-1})>|\sin\omega|/\omega$, and some fan multiplier attains
equality in the first inequality. For $n$ equal sectors the minimum of $\max_l|w_l|$ over all fan multipliers is therefore $|\sin\omega|/(n\sin(\omega/n))$.
\item[(b)] Let $0<\gamma$, $2\gamma<\omega$ and $(\omega-2\gamma)/n<\pi$, and consider the partition into two edge sectors of opening $\gamma$ and $n$ equal
sectors in between. Over all fan multipliers for it with $w=0$ on both edge sectors, the minimum of $\max_l|w_l|$ is
\begin{equation}\label{eq:kng}
\kappa_{n,\gamma}(\omega):=\frac{|\sin\omega|}{n\sin\frac{\omega-2\gamma}{n}},
\end{equation}
which tends to $|\sin\omega|/\omega$ as $n\to\infty$, $\gamma\to0$. On each edge sector such a fan multiplier is the reflection multiplier of the adjacent side (by (c)).
\item[(c)] If a fan multiplier has $w_1=0$, then it coincides on $S_1$ with the reflection multiplier \eqref{eq:mell} of the side $\theta=\theta_0$; if it
has $w_N=0$, then it coincides on $S_N$ with the reflection multiplier of the side $\theta=\theta_0+\omega$.
\end{enumerate}
\end{lemma}
\begin{proof}
We use the notation $N_l$ and \eqref{eq:R} of the proof of Lemma~\ref{lem:fanid}.

\emph{(a).} By \eqref{fan} and the triangle inequality,
\[2|\sin\omega|=\Big|\sum_lw_l\Delta E_l\Big|\le\sum_l|w_l|\,|\Delta E_l|\le\max_l|w_l|\cdot2\sum_l\sin(\vartheta_l-\vartheta_{l-1}),\]
which is the first inequality. The second follows from $\sin x<x$ for $x>0$, which gives $\sum_l\sin(\vartheta_l-\vartheta_{l-1})<\sum_l(\vartheta_l-\vartheta_{l-1})=\omega$,
together with $|\sin\omega|>0$ ($\omega\ne\pi$). For equality, put $u:=\Delta E/|\Delta E|$ and
$\kappa_*:=|\Delta E|/\sum_l|\Delta E_l|=|\sin\omega|/\sum_l\sin(\vartheta_l-\vartheta_{l-1})$, and choose
\[w_l:=\kappa_*\,u\,\frac{\overline{\Delta E_l}}{|\Delta E_l|}.\]
Then $w_l\Delta E_l=\kappa_*u|\Delta E_l|$, so $\sum_lw_l\Delta E_l=\kappa_*u\sum_l|\Delta E_l|=u|\Delta E|=\Delta E$: \eqref{fan} holds and $|w_l|=\kappa_*$ for all $l$. By the
converse in Lemma~\ref{lem:fanid} there is exactly one fan multiplier with these weights. (Both inequalities in the triangle estimate become equalities because all terms $w_l\Delta E_l$ point
in the same direction $u$ and all $|w_l|$ are equal; see Figure~\ref{fig:fan}(b), (c).) For $n$ equal sectors, $\vartheta_l-\vartheta_{l-1}=\omega/n$ and
$\sum_l\sin(\vartheta_l-\vartheta_{l-1})=n\sin(\omega/n)$.

\emph{(b).} By hypothesis the edge sectors have opening $\gamma<\omega/2<\pi$ and the $n$ inner sectors opening $(\omega-2\gamma)/n\in(0,\pi)$, so
this partition satisfies the hypothesis of Definition~\ref{def:fan}. The edge sectors have $w=0$ and contribute nothing to $\sum_lw_l\Delta E_l$. The triangle inequality over the
$n$ inner sectors, each of opening $(\omega-2\gamma)/n$, gives $\max|w_l|\ge|\sin\omega|/(n\sin\frac{\omega-2\gamma}{n})=\kappa_{n,\gamma}(\omega)$, and the choice of
(a) on the inner sectors, with $w=0$ on the edge sectors, attains it. Since $x-x^3/(6n^2)\le n\sin(x/n)\le x$ for $x\ge0$, with $x=\omega-2\gamma$ we get
$\kappa_{n,\gamma}(\omega)\to|\sin\omega|/\omega$ as $n\to\infty$ and $\gamma\to0$. The last sentence of (b) is (c).

\emph{(c).} If $w_1=0$, then $b_1=1$, $N_1=N_0+\Delta E_1=E_1-E_0$ and
$a_1=N_1-E_1=-E_0$, so $m(z)=\bar z-e^{-2i\theta_0}z$ on $S_1$. This is the reflection multiplier $m_\ell$ of the side $\theta=\theta_0$: take $z_0=0$ and the
inward normal $\nu=ie^{i\theta_0}$ of that side in Lemma~\ref{lem:refl}, so that $\bar\nu^2=-e^{-2i\theta_0}$. Likewise $w_N=0$ gives $b_N=1$, $a_N=N_N-E_N=-E_N$, and
$m(z)=\bar z-e^{-2i(\theta_0+\omega)}z$ on $S_N$, the reflection multiplier of the side $\theta=\theta_0+\omega$ ($\nu=-ie^{i(\theta_0+\omega)}$).
\end{proof}

\begin{definition}[passive/active split]\label{def:split}
Let $\mathcal T$ be a triangulation of $\overline\Om$ and $m$ a continuous, piecewise linear multiplier on $\mathcal T$. Its weight \eqref{eq:W} is constant
on each $T\in\mathcal T$. A \emph{passive/active split} is a partition $\mathcal T=\mathcal P\cup\mathcal A$ into \emph{passive} triangles $\mathcal P$ and
\emph{active} triangles $\mathcal A$ such that every triangle meeting $\partial\Om$ is passive. The \emph{triangle weights} and the \emph{passive constant} are
\begin{equation}\label{eq:wT}
w_T:=w|_T=1-\dzb m|_T\quad(T\in\mathcal T),\qquad \kappa_{\mathcal P}:=\max_{T\in\mathcal P}|w_T| .
\end{equation}
\end{definition}
In what follows such a split is given and a number $\kappa_0>\kappa_{\mathcal P}$ is fixed.

\begin{remark}\label{rem:fan}
The limit $|\sin\omega|/\omega$ in Lemma~\ref{lem:fan} coincides with the corner contribution to $\kess$ in \eqref{eq:kess}: the pointwise estimate
available from Lemma~\ref{lem:id} is sharp at corners. All uses of Lemma~\ref{lem:fan} satisfy its hypotheses: in Lemma~\ref{lem:goodmult} the inner
sectors have opening $(\omega_j-2\gamma)/n<\pi$ and the edge sectors opening $\gamma<\omega_j/2<\pi$, and the fans of Section~\ref{sec:num} consist of
equal sectors of opening $\pi/(2n)\le\pi/16$. The certified results do not depend on Lemmas~\ref{lem:fanid} and~\ref{lem:fan}, since the verifier
computes every $w_T$ exactly from the nodal values of $m$ (Section~\ref{sec:impl}); the lemmas only predict the value of $\kappa_{\mathcal P}$.
\end{remark}

For $\zeta\in\C$ let
\begin{equation}\label{q}
\begin{gathered}
q_T(\zeta):=\kappa_0|\zeta|^2-\Re(w_T\zeta^2)\quad(T\in\mathcal P),\qquad q_T(\zeta):=\kappa_0|\zeta|^2\quad(T\in\mathcal A),\\
\|f\|_q^2:=\sum_{T\in\mathcal T}\int_Tq_T(f)\,dA ,
\end{gathered}
\end{equation}
and let $q(f)$ be the function with $q(f):=q_T(f)$ on each $T$. The polarisations of $q_T$ and the corresponding form on $L^2$ are
\begin{equation}\label{eq:qip}
\begin{gathered}
b_T(\zeta,\xi):=\kappa_0\Re(\zeta\bar\xi)-\Re(w_T\zeta\xi)\quad(T\in\mathcal P),\qquad b_T(\zeta,\xi):=\kappa_0\Re(\zeta\bar\xi)\quad(T\in\mathcal A),\\
\langle f,g\rangle_q:=\sum_{T\in\mathcal T}\int_Tb_T(f,g)\,dA,
\end{gathered}
\end{equation}
and the pointwise lower bound of $q$ is the piecewise constant function $c$ with
\begin{equation}\label{eq:cfun}
c:=\kappa_0-|w_T|\ \text{ on }T\in\mathcal P,\qquad c:=\kappa_0\ \text{ on }T\in\mathcal A,\qquad c_{\min}:=\kappa_0-\kappa_{\mathcal P}>0 .
\end{equation}
\begin{lemma}[the form $q$]\label{lem:qform}
Identify $\zeta=\xi_1+i\xi_2\in\C$ with $\xi=(\xi_1,\xi_2)^{\mathsf T}\in\R^2$, and for $T\in\mathcal P$ write $w_T=:\sigma_T+i\tau_T$ with
$\sigma_T,\tau_T\in\R$. Put
\[W_T:=\begin{pmatrix}\sigma_T&-\tau_T\\-\tau_T&-\sigma_T\end{pmatrix},\quad A_T:=\kappa_0I-W_T\quad(T\in\mathcal P),\qquad A_T:=\kappa_0I\quad(T\in\mathcal A).\]
\begin{enumerate}
\item[(a)] $q_T(\zeta)=\xi^{\mathsf T}A_T\xi$. For $T\in\mathcal P$, $W_T$ is symmetric with trace $0$ and determinant $-|w_T|^2$, $W_T^2=|w_T|^2I$, its
eigenvalues are $\pm|w_T|$, and the eigenvalues $\kappa_0\mp|w_T|$ of $A_T$ are positive.
\item[(b)] On every $T\in\mathcal T$, $c\,|\zeta|^2\le q_T(\zeta)\le(\kappa_0+\kappa_{\mathcal P})|\zeta|^2$; on $T\in\mathcal P$,
$(\kappa_0-|w_T|)|\zeta|^2\le q_T(\zeta)\le(\kappa_0+|w_T|)|\zeta|^2$.
\item[(c)] $\langle\cdot,\cdot\rangle_q$ is a real inner product on $L^2$ with $\langle f,f\rangle_q=\|f\|_q^2$, and
$c_{\min}\|f\|^2\le\|f\|_q^2\le(\kappa_0+\kappa_{\mathcal P})\|f\|^2$ for $f\in L^2$.
\end{enumerate}
\end{lemma}
\begin{proof}
(a) For $T\in\mathcal P$, $\Re(w_T\zeta^2)=\sigma_T(\xi_1^2-\xi_2^2)-2\tau_T\xi_1\xi_2=\xi^{\mathsf T}W_T\xi$, and $\kappa_0|\zeta|^2=\kappa_0\xi^{\mathsf T}\xi$;
on $\mathcal A$ the claim is immediate. The trace and determinant of $W_T$ are read off; a symmetric $2\times2$ matrix with trace $0$ and determinant
$-|w_T|^2$ has eigenvalues $\pm|w_T|$ and satisfies $W_T^2=|w_T|^2I$. The eigenvalues $\kappa_0\mp|w_T|$ of $A_T$ are positive because
$|w_T|\le\kappa_{\mathcal P}<\kappa_0$.
(b) By (a), $\xi^{\mathsf T}A_T\xi$ lies between the smallest and the largest eigenvalue of $A_T$ times $|\xi|^2=|\zeta|^2$; these are
$\kappa_0-|w_T|=c$ and $\kappa_0+|w_T|\le\kappa_0+\kappa_{\mathcal P}$ on $T\in\mathcal P$, and $\kappa_0=c$ on $T\in\mathcal A$.
(c) $b_T$ is real bilinear and symmetric with $b_T(\zeta,\zeta)=q_T(\zeta)$, so $\langle\cdot,\cdot\rangle_q$ is symmetric and bilinear with
$\langle f,f\rangle_q=\|f\|_q^2$. Integrating (b) and using $c\ge c_{\min}$ gives the norm equivalence, which also shows that
$\langle\cdot,\cdot\rangle_q$ is positive definite.
\end{proof}

The next proposition is the key reformulation of the bound $\kappa(\Om)\le\kappa_0$.
\begin{proposition}[reduction to the active triangles]\label{prop:goal}
Let $m$ be a continuous, piecewise linear multiplier, $\mathcal T=\mathcal P\cup\mathcal A$ a passive/active split, and $\kappa_0>\kappa_{\mathcal P}$.
\begin{enumerate}
\item[(a)] For every $f\in\A$,
\begin{equation}\label{goal}
\Re\!\int_\Om f^2-\kappa_0\|f\|^2=\sum_{T\in\mathcal A}\Re\Big(w_T\!\int_Tf^2\Big)-\|f\|_q^2 .
\end{equation}
\item[(b)] $\kappa(\Om)\le\kappa_0$ if and only if $\sum_{T\in\mathcal A}\Re(w_T\int_Tf^2)\le\|f\|_q^2$ for all $f\in\Az$.
\item[(c)] The quadratic form $f\mapsto\sum_{T\in\mathcal A}\Re(w_T\int_Tf^2)$ is compact on $\A$.
\end{enumerate}
\end{proposition}
\begin{proof}
(a) Since $f^2\in L^1$ and $w$ is bounded and piecewise constant, Lemma~\ref{lem:id} reads
$\int_\Om f^2=\sum_{T\in\mathcal T}w_T\int_Tf^2$. Take real parts and subtract $\kappa_0\|f\|^2=\sum_T\kappa_0\int_T|f|^2$. Because
$\Re(w_T\int_Tf^2)=\int_T\Re(w_Tf^2)$, the contribution of $T$ is
\[\Re\Big(w_T\!\int_Tf^2\Big)-\kappa_0\!\int_T|f|^2=\begin{cases}-\int_Tq_T(f)\,dA,&T\in\mathcal P,\\[2pt] \Re\big(w_T\int_Tf^2\big)-\int_Tq_T(f)\,dA,&T\in\mathcal A;\end{cases}\]
summing over $T$ gives \eqref{goal}.
(b) By \eqref{eq:kappa-re}, $\kappa(\Om)\le\kappa_0$ if and only if $\Re\int f^2\le\kappa_0\|f\|^2$ for all $f\in\Az$; by \eqref{goal} this is the stated condition.
(c) Every active triangle is closed and does not meet $\partial\Om$, so the union of the active triangles is a compact subset of $\Om$; the claim is
Lemma~\ref{lem:bergman}(c).
\end{proof}
Nothing is lost on passive triangles; they are needed only to keep the active triangles, on which $|w_T|$ may be large, away from the boundary.

\section{A computable certificate}\label{sec:cert}
\emph{Standing assumptions.} In Sections~\ref{sec:cert} and~\ref{sec:conv}, $\Om$ is a polygonal domain, $\mathcal T$ a triangulation of $\overline\Om$,
$m$ a continuous piecewise linear multiplier on $\mathcal T$ with triangle weights $w_T$, and $\mathcal T=\mathcal P\cup\mathcal A$ a split into passive and active
triangles in which every triangle meeting $\partial\Om$ is passive; $\kappa_0>\kappa_{\mathcal P}$ is fixed. The forms $q_T$, the norm $\|\cdot\|_q$ \eqref{q}, the real
inner product $\langle\cdot,\cdot\rangle_q$, the piecewise constant function $c$ and $c_{\min}=\kappa_0-\kappa_{\mathcal P}$ are those of Section~\ref{sec:mult}. We use the
matrices $W_T$, $A_T$ of Lemma~\ref{lem:qform}: identifying $\C$ with $\R^2$, so that $\Re(\zeta\bar\zeta')$ is the Euclidean product $\zeta\cdot\zeta'$,
one has $q_T(\zeta)=\zeta^{\mathsf T}A_T\zeta$ and, by \eqref{eq:qip}, $\langle f,g\rangle_q=\sum_T\int_Tf\cdot A_Tg\,dA$; the bounds $c\,|f|^2\le q(f)$ and
$c_{\min}\|f\|^2\le\|f\|_q^2\le(\kappa_0+\kappa_{\mathcal P})\|f\|^2$ are Lemma~\ref{lem:qform}(b), (c).

By Proposition~\ref{prop:goal}(b), (c), $\kappa(\Om)\le\kappa_0$ holds if and only if the compact quadratic form
$f\mapsto\sum_{T\in\mathcal A}\Re(w_T\int_Tf^2)$ is bounded by $\|f\|_q^2$ on $\Az$; this section reduces this inequality to a finite matrix
inequality whose data are integrals of explicit functions.

\subsection{Cells and Taylor functionals}
\begin{definition}[cells and Taylor functionals]\label{def:cells}
A \emph{cell decomposition} consists of finitely many closed triangles $C$, the \emph{cells}, each contained in an active triangle, which together partition
$\bigcup_{T\in\mathcal A}T$ up to null sets; each cell $C$ carries a \emph{centre} $a_C$ and radii $0<r_C<R_C$ with $C\subset\overline{D(a_C,r_C)}$ and
$D(a_C,R_C)\subset\Om$. Several cells may share a centre. We put
\begin{equation}\label{eq:cell}
\theta_C:=\frac{r_C}{R_C}\in(0,1),\qquad w_C:=w_T\quad\text{for the active triangle }T\supset C .
\end{equation}
Fix a degree $p\ge0$ and a scale $s_a>0$ for each centre $a$. Enumerate the pairs $(a,n)$, $a$ a centre and $0\le n\le p$, as $(a_\mu,n_\mu)$, $\mu=1,\dots,N$.
The \emph{Taylor functionals} and the \emph{scaled monomials} are
\begin{equation}\label{eq:ell}
\ell_\mu(f):=s_{a_\mu}^{n_\mu}\,\frac{f^{(n_\mu)}(a_\mu)}{n_\mu!}\quad(f\in\A),\qquad\phi_\mu(z):=s_{a_\mu}^{-n_\mu}(z-a_\mu)^{n_\mu}\qquad(\mu=1,\dots,N),
\end{equation}
and $\ell(f):=(\ell_1(f),\dots,\ell_N(f))\in\C^N$. For a cell $C$ let $J_C:=\{\mu:a_\mu=a_C\}$ and let
\begin{equation}\label{eq:tC}
t_C:=\sum_{\mu\in J_C}\ell_\mu(f)\,\phi_\mu
\end{equation}
be the Taylor polynomial of $f$ of degree $p$ at $a_C$.
\end{definition}
The scales $s_a$ affect only the conditioning of the matrices below; their choice is described in Section~\ref{sec:impl}.

\begin{lemma}[Taylor remainder]\label{lem:taylor}
Let $C$ be a cell. Then $\|f-t_C\|_{L^2(C)}\le\theta_C^{\,p+2}\|f\|_{L^2(D(a_C,R_C))}$ for every $f\in\A$.
\end{lemma}
\begin{proof}
Write $a:=a_C$, $r:=r_C$, $R:=R_C$, $\theta:=\theta_C$ and $\hat f_n:=f^{(n)}(a)/n!$. Since $D(a,R)\subset\Om$, $f(z)=\sum_{n\ge0}\hat f_n(z-a)^n$, the series
converging absolutely and uniformly on $\overline{D(a,\varrho)}$ for every $\varrho<R$ (Lemma~\ref{lem:holo}(b)). In polar coordinates $z=a+se^{i\phi}$,
\[\int_{D(a,\varrho)}(z-a)^n\,\overline{(z-a)^k}\,dA=\int_0^\varrho s^{n+k+1}ds\int_0^{2\pi}e^{i(n-k)\phi}d\phi=\begin{cases}\pi\varrho^{2n+2}/(n+1),&k=n,\\0,&k\ne n.\end{cases}\]
For $\varrho<R$, uniform convergence allows term-by-term integration of $|f|^2=\sum_{n,k}\hat f_n\overline{\hat f_k}(z-a)^n\overline{(z-a)^k}$, so
$\|f\|^2_{D(a,\varrho)}=\sum_n|\hat f_n|^2\pi\varrho^{2n+2}/(n+1)$; letting $\varrho\uparrow R$ (monotone convergence on both sides), the same holds for $\varrho=R$.
Applying this to $f-t_C=\sum_{n\ge p+1}\hat f_n(z-a)^n$ on $D(a,r)$ and to $f$ on $D(a,R)$,
\[\|f-t_C\|^2_{D(a,r)}=\sum_{n\ge p+1}|\hat f_n|^2\frac{\pi R^{2n+2}}{n+1}\,\theta^{2n+2}\le\theta^{2p+4}\sum_{n\ge p+1}|\hat f_n|^2\frac{\pi R^{2n+2}}{n+1}\le\theta^{2p+4}\|f\|^2_{D(a,R)},\]
because $\theta<1$ and $2n+2\ge2p+4$ for $n\ge p+1$. Finally $C\subset\overline{D(a,r)}$, and the circle $\partial D(a,r)$ has measure zero, so
$\|f-t_C\|_{L^2(C)}\le\|f-t_C\|_{D(a,r)}$.
\end{proof}

\subsection{The certificate quantities}
We identify $\C^N$ with $\R^{2N}$ by
\begin{equation}\label{eq:ident}
x=(x_1,\dots,x_N)\longmapsto(\Re x_1,\Im x_1,\dots,\Re x_N,\Im x_N);
\end{equation}
then $\Re(y^*x)=y\cdot x$ for $x,y\in\C^N$, with $y\cdot x$ the Euclidean inner product of $\R^{2N}$. A real quadratic form on $\C^N$ is thereby a real quadratic form on $\R^{2N}$, which we identify with its
symmetric $2N\times2N$ matrix; $F(x,y)$ denotes the associated symmetric bilinear form, $\|F\|$ the spectral norm, and $F\preceq F'$ the Loewner order.
For $\eta>0$ let $Q$ be the real quadratic form
\begin{equation}\label{eq:Q}
\begin{split}Q(x):=\sum_{\text{cells }C}\Big[&\Re\Big(w_C\sum_{\mu,\nu\in J_C}x_\mu x_\nu\int_C\phi_\mu\phi_\nu\,dA\Big)\\
&+\eta|w_C|\sum_{\mu,\nu\in J_C}\bar x_\mu x_\nu\int_C\bar\phi_\mu\phi_\nu\,dA\Big]\qquad(x\in\C^N).\end{split}
\end{equation}
Inserting $x=\ell(f)$ into \eqref{eq:Q} and using \eqref{eq:tC}, i.e.
$\int_Ct_C^2\,dA=\sum_{\mu,\nu\in J_C}x_\mu x_\nu\int_C\phi_\mu\phi_\nu\,dA$ and $\|t_C\|_{L^2(C)}^2=\sum_{\mu,\nu\in J_C}\bar x_\mu x_\nu\int_C\bar\phi_\mu\phi_\nu\,dA$, gives
\begin{equation}\label{eq:Qell}
Q(\ell(f))=\sum_{\text{cells }C}\Big[\Re\Big(w_C\!\int_Ct_C^2\,dA\Big)+\eta|w_C|\,\|t_C\|_{L^2(C)}^2\Big]\qquad(f\in\A).
\end{equation}
For $v\in\C$ let $W(v)$ be the real symmetric matrix of the quadratic form $\xi\mapsto\Re(v\xi^2)$ on $\R^2\cong\C$,
\begin{equation}\label{eq:Wv}
W(v):=\begin{pmatrix}\Re v&-\Im v\\-\Im v&-\Re v\end{pmatrix};
\end{equation}
thus $W_T=W(w_T)$ on $\mathcal P$.
The \emph{dual form} of $q_T$ and its coefficients are
\begin{equation}\label{qstar}
\begin{gathered}
q_T^*(\xi):=\hat\kappa_T|\xi|^2+\Re(\hat w_T\xi^2)\qquad(\xi\in\C),\\
\hat\kappa_T:=\frac{\kappa_0}{\kappa_0^2-|w_T|^2},\quad\hat w_T:=\frac{w_T}{\kappa_0^2-|w_T|^2}\quad(T\in\mathcal P),\qquad\hat\kappa_T:=\frac1{\kappa_0},\quad\hat w_T:=0\quad(T\in\mathcal A),
\end{gathered}
\end{equation}
and $q^*(k)$ denotes the function equal to $q_T^*(k)$ on $T$. Thus $q_T^*(\xi)=(\kappa_0|\xi|^2+\Re(w_T\xi^2))/(\kappa_0^2-|w_T|^2)$ on $\mathcal P$ and
$q_T^*(\xi)=|\xi|^2/\kappa_0$ on $\mathcal A$.
\begin{lemma}[dual form]\label{lem:dualform}
Let $T\in\mathcal T$.
\begin{enumerate}
\item[(a)] For $v,\zeta\in\C$ the vector $W(v)\zeta$ is the complex number $\bar v\bar\zeta$, and $W(v)^2=|v|^2I$.
\item[(b)] $A_T^{-1}=\hat\kappa_TI+W(\hat w_T)$; hence $q_T^*(\xi)=\xi^{\mathsf T}A_T^{-1}\xi$, and $A_T^{-1}\xi$ is the number $\hat\kappa_T\xi+\overline{\hat w_T}\,\bar\xi$.
\item[(c)] $\Re(\zeta\bar\xi)^2\le q_T(\zeta)\,q_T^*(\xi)$ for all $\zeta,\xi\in\C$, with equality for $\xi=A_T\zeta$, i.e.\ for $\bar\xi=\kappa_0\bar\zeta-w_T\zeta$ if
$T\in\mathcal P$ and $\xi=\kappa_0\zeta$ if $T\in\mathcal A$. Thus $q_T^*(\xi)=\max_{\zeta\ne0}\Re(\zeta\bar\xi)^2/q_T(\zeta)$.
\end{enumerate}
\end{lemma}
\begin{proof}
(a) Write $v=v_1+iv_2$, $\zeta=\xi_1+i\xi_2$. Then $W(v)\zeta=(v_1\xi_1-v_2\xi_2,\,-v_2\xi_1-v_1\xi_2)$, which is $\bar v\bar\zeta=(v_1-iv_2)(\xi_1-i\xi_2)$.
Applying this twice, $W(v)^2\zeta=\overline{\bar v\bar\zeta}\,\bar v=|v|^2\zeta$.

(b) On $\mathcal A$, $A_T=\kappa_0I$ and $\hat\kappa_T=1/\kappa_0$, $\hat w_T=0$. On $\mathcal P$, by (a),
$(\kappa_0I-W(w_T))(\kappa_0I+W(w_T))=(\kappa_0^2-|w_T|^2)I$, and $\kappa_0^2-|w_T|^2>0$ because $|w_T|\le\kappa_{\mathcal P}<\kappa_0$. Hence
$A_T^{-1}=(\kappa_0I+W(w_T))/(\kappa_0^2-|w_T|^2)=\hat\kappa_TI+W(\hat w_T)$, $W$ being real-linear in $v$. By \eqref{eq:Wv},
$\xi^{\mathsf T}A_T^{-1}\xi=\hat\kappa_T|\xi|^2+\Re(\hat w_T\xi^2)=q_T^*(\xi)$, and by (a) $A_T^{-1}\xi=\hat\kappa_T\xi+\overline{\hat w_T}\,\bar\xi$.

(c) $A_T$ is positive definite (its eigenvalues are $\kappa_0\mp|w_T|>0$ on $\mathcal P$). The Cauchy--Schwarz inequality for the inner product
$(\zeta,\zeta')\mapsto\zeta^{\mathsf T}A_T\zeta'$ gives
\[\Re(\zeta\bar\xi)^2=\big(\zeta^{\mathsf T}A_T(A_T^{-1}\xi)\big)^2\le\big(\zeta^{\mathsf T}A_T\zeta\big)\big(\xi^{\mathsf T}A_T^{-1}A_TA_T^{-1}\xi\big)=q_T(\zeta)\,q_T^*(\xi),\]
with equality if and only if $\zeta$ and $A_T^{-1}\xi$ are parallel, e.g.\ for $\xi=A_T\zeta$. On $\mathcal P$, (a) gives $A_T\zeta=\kappa_0\zeta-\bar w_T\bar\zeta$, i.e.\
$\bar\xi=\kappa_0\bar\zeta-w_T\zeta$.
\end{proof}
By Lemma~\ref{lem:dualform}(c) pointwise and the Cauchy--Schwarz inequality in $L^2(\Om)$,
$(\Re\int_\Om f\bar k)^2\le(\int_\Om q(f)^{1/2}q^*(k)^{1/2})^2\le\int_\Om q(f)\int_\Om q^*(k)$, that is,
\begin{equation}\label{dual}\Big(\Re\!\int_\Om f\,\bar k\,dA\Big)^2\le\|f\|_q^2\int_\Om q^*(k)\,dA\qquad(f,k\in L^2(\Om)).\end{equation}

\subsection{Representers}
\begin{definition}[bump representer]\label{def:bump}
For $a\in\Om$, $R'>0$ with $D(a,R')\subset\Om$, and an integer $n\ge0$ let
\begin{equation}\label{eq:bump}
b_{a,n}(z):=\begin{cases}\mathcal N^{-1}\psi(|z-a|)\,(z-a)^n,&|z-a|<R',\\0,&\text{otherwise},\end{cases}\qquad\psi(\varrho):=1-\frac{\varrho^2}{R'^2},
\end{equation}
with the normalisation $\mathcal N:=\pi R'^{2n+2}/((n+1)(n+2))$.
For each centre $a$ fix a radius $R'_a\le R_C$ for all cells $C$ with centre $a$ (so that $D(a,R'_a)\subset\Om$), and let $b_\mu:=s_{a_\mu}^{n_\mu}\,b_{a_\mu,n_\mu}$
with $R'=R'_{a_\mu}$ ($\mu=1,\dots,N$).
\end{definition}
\begin{lemma}\label{lem:bump}
In the setting of Definition~\ref{def:bump}, $b_{a,n}$ is Lipschitz on $\C$ and $\langle f,b_{a,n}\rangle=f^{(n)}(a)/n!$ for all $f\in\A$. In particular
$\langle f,b_\mu\rangle=\ell_\mu(f)$ for $f\in\A$, $\mu=1,\dots,N$.
\end{lemma}
\begin{proof}
Because $\psi(R')=0$, $b_{a,n}$ is continuous and piecewise smooth with bounded gradient, hence Lipschitz. Let $\hat f_k:=f^{(k)}(a)/k!$. The Taylor series of
$f$ at $a$ converges to $f$ in $L^2(D(a,R'))$ (proof of Lemma~\ref{lem:taylor}) and $b_{a,n}$ is bounded, so in polar coordinates
\begin{align*}\langle f,b_{a,n}\rangle&=\frac1{\mathcal N}\sum_k\hat f_k\int_{D(a,R')}(z-a)^k\,\overline{(z-a)^n}\,\psi(|z-a|)\,dA\\
&=\frac{\hat f_n}{\mathcal N}\,2\pi\!\int_0^{R'}\varrho^{2n+1}\Big(1-\frac{\varrho^2}{R'^2}\Big)d\varrho=\hat f_n,\end{align*}
since the angular integral vanishes for $k\ne n$ and $2\pi R'^{2n+2}\big(\frac1{2n+2}-\frac1{2n+4}\big)=\mathcal N$. As $s_{a_\mu}^{n_\mu}$ is real,
$\langle f,b_\mu\rangle=s_{a_\mu}^{n_\mu}\langle f,b_{a_\mu,n_\mu}\rangle=\ell_\mu(f)$.
\end{proof}
A function $k\in L^2(\Om)$ is a \emph{representer} of $\ell_\mu$ if $\langle f,k\rangle=\ell_\mu(f)$ for all $f\in\Az$. For Lipschitz functions $\chi_\mu$ on
$\overline\Om$ with $\chi_\mu|_{\partial\Om}=0$, called \emph{correctors}, and constants $\gamma_\mu\in\C$ put
\begin{equation}\label{rep}
k_\mu:=b_\mu-\dz\chi_\mu-\gamma_\mu\qquad(\mu=1,\dots,N).
\end{equation}
The correctors are arbitrary; in practice they are finite element approximations of the solution of an elliptic problem (Lemma~\ref{lem:optcorr}).
\begin{lemma}[representers]\label{lem:rep}
\begin{enumerate}
\item[(a)] For every Lipschitz $\chi$ on $\overline\Om$ with $\chi|_{\partial\Om}=0$ and every $\gamma\in\C$, $\langle f,\dz\chi+\gamma\rangle=0$ for all $f\in\Az$.
Consequently every $k_\mu$ of the form \eqref{rep} is a representer of $\ell_\mu$.
\item[(b)] Conversely, every representer of $\ell_\mu$ in $L^2(\Om)$ is of the form $b_\mu-\dz\chi-\gamma$ with $\chi\in H^1_0(\Om;\C)$ and $\gamma\in\C$.
\end{enumerate}
\end{lemma}
\begin{proof}
(a) $\langle f,\dz\chi\rangle=\int f\,\overline{\dz\chi}=\int f\,\dzb\bar\chi=0$ by \eqref{eq:ibp} with $F=f\in L^1(\Om)$ and $u=\bar\chi$, and
$\langle f,\gamma\rangle=\bar\gamma\int f=0$. With Lemma~\ref{lem:bump}, $\langle f,k_\mu\rangle=\langle f,b_\mu\rangle=\ell_\mu(f)$.
(b) If $k$ is a representer, then $b_\mu-k\perp\Az$ by Lemma~\ref{lem:bump}, so $b_\mu-k\in(\Az)^\perp=\dz H^1_0(\Om)\oplus\C$ (Lemma~\ref{lem:decomp}).
\end{proof}
The \emph{Gram form} of representers $k_1,\dots,k_N$ is the real quadratic form
\begin{equation}\label{gram}
\mathcal R(y):=\int_\Om q^*\Big(\sum_{\mu=1}^Ny_\mu k_\mu\Big)dA\qquad(y\in\C^N),
\end{equation}
regarded via \eqref{eq:ident} as a real quadratic form on $\R^{2N}$.
\begin{lemma}[matrix of the Gram form]\label{lem:grammatrix}
Let $R_1,R_2\in\C^{N\times N}$ be given by
\[(R_1)_{\mu\nu}:=\sum_{T\in\mathcal T}\hat\kappa_T\int_T\bar k_\mu k_\nu\,dA,\qquad(R_2)_{\mu\nu}:=\sum_{T\in\mathcal T}\hat w_T\int_Tk_\mu k_\nu\,dA .\]
Then $R_1$ is Hermitian, $R_2$ is complex
symmetric, $\mathcal R(y)=y^*R_1y+\Re(y^{\mathsf T}R_2y)$, and $\mathcal R\succeq0$. In the coordinates \eqref{eq:ident}, the $2\times2$ block of the matrix of
$\mathcal R$ for the coordinates $(\Re y_\mu,\Im y_\mu)$ and $(\Re y_\nu,\Im y_\nu)$ is
\[\begin{pmatrix}\Re(R_1+R_2)_{\mu\nu}&-\Im(R_1+R_2)_{\mu\nu}\\ \Im(R_1-R_2)_{\mu\nu}&\Re(R_1-R_2)_{\mu\nu}\end{pmatrix}.\]
\end{lemma}
\begin{proof}
For $K:=\sum_\mu y_\mu k_\mu$ one has $|K|^2=\sum_{\mu,\nu}\bar y_\mu y_\nu\bar k_\mu k_\nu$ and $K^2=\sum_{\mu,\nu}y_\mu y_\nu k_\mu k_\nu$; inserting these into
\eqref{qstar} and integrating over each $T$ gives $\mathcal R(y)=y^*R_1y+\Re(y^{\mathsf T}R_2y)$. The symmetries are evident, and $\mathcal R\ge0$ because $q_T^*$
is positive semidefinite (Lemma~\ref{lem:dualform}(b)). For the block, write $y=u+iv$ with $u,v\in\R^N$ and expand: $y^*R_1y=u^{\mathsf T}(\Re R_1)u+v^{\mathsf T}(\Re
R_1)v-2u^{\mathsf T}(\Im R_1)v$ (as $\Im R_1$ is antisymmetric) and $\Re(y^{\mathsf T}R_2y)=u^{\mathsf T}(\Re R_2)u-v^{\mathsf T}(\Re R_2)v-2u^{\mathsf T}(\Im R_2)v$;
the stated symmetric matrix has this quadratic form.
\end{proof}

\emph{Real-linear representers.} The form $q^*$ is not invariant under $k\mapsto ik$ on passive triangles, and the problem that defines an optimal corrector is
only real-linear (Lemma~\ref{lem:optcorr}); for complex-linear combinations as in \eqref{gram} the corrector error cannot vanish in general
(Proposition~\ref{prop:complexlin}). We therefore also admit separate representers for the real and the imaginary parts of the Taylor functionals.
\begin{definition}[real-linear representers and their Gram form]\label{def:reallin}
Let
\begin{equation}\label{eq:lambda}
\lambda_{2\mu-1}:=\Re\ell_\mu,\qquad\lambda_{2\mu}:=\Im\ell_\mu\qquad(\mu=1,\dots,N),
\end{equation}
so that $\ell(f)=(\lambda_1(f),\dots,\lambda_{2N}(f))$ in the coordinates \eqref{eq:ident}. A \emph{real representer} of $\lambda_j$ is a function $\tilde k\in L^2(\Om)$ with
$\Re\langle f,\tilde k\rangle=\lambda_j(f)$ for all $f\in\Az$. The admissible real representers are
\begin{equation}\label{eq:tildek}
\tilde k_j:=\hat b_j-\dz\chi_j-\gamma_j\quad(j=1,\dots,2N),\qquad\hat b_{2\mu-1}:=b_\mu,\quad\hat b_{2\mu}:=i\,b_\mu ,
\end{equation}
with Lipschitz correctors $\chi_j$ vanishing on $\partial\Om$ and constants $\gamma_j\in\C$, chosen independently for each $j$. The Gram form of real representers
$\tilde k_1,\dots,\tilde k_{2N}$ is
\begin{equation}\label{eq:G}
\mathsf G(\alpha):=\int_\Om q^*\Big(\sum_{j=1}^{2N}\alpha_j\tilde k_j\Big)dA\qquad(\alpha\in\R^{2N}).
\end{equation}
\end{definition}
\begin{lemma}\label{lem:reallin}
\begin{enumerate}
\item[(a)] Each $\tilde k_j$ of the form \eqref{eq:tildek} is a real representer of $\lambda_j$. For any real representers $\tilde k_1,\dots,\tilde k_{2N}$ of
$\lambda_1,\dots,\lambda_{2N}$, one has $y\cdot\ell(f)=\Re\langle f,\sum_jy_j\tilde k_j\rangle$ for all $y\in\R^{2N}$ and $f\in\Az$.
\item[(b)] $\mathsf G(\alpha)=\alpha^{\mathsf T}\mathsf G\alpha$ with the positive semidefinite matrix
$\mathsf G_{ij}:=\sum_T\big[\hat\kappa_T\Re\int_T\overline{\tilde k_i}\,\tilde k_j\,dA+\Re\big(\hat w_T\int_T\tilde k_i\tilde k_j\,dA\big)\big]$.
\item[(c)] If $\tilde k_{2\mu-1}=k_\mu$ and $\tilde k_{2\mu}=ik_\mu$ for all $\mu$, then $\mathsf G$ is the matrix of the form \eqref{gram} in the coordinates \eqref{eq:ident}.
\end{enumerate}
\end{lemma}
\begin{proof}
(a) By Lemma~\ref{lem:bump}, $\Re\langle f,b_\mu\rangle=\Re\ell_\mu(f)$ and $\Re\langle f,ib_\mu\rangle=\Re(-i\ell_\mu(f))=\Im\ell_\mu(f)$; the corrector and the constant
do not contribute by Lemma~\ref{lem:rep}(a). Since the $y_j$ are real, $\Re\langle f,\sum_jy_j\tilde k_j\rangle=\sum_jy_j\lambda_j(f)=y\cdot\ell(f)$.
(b) For $K:=\sum_j\alpha_j\tilde k_j$ with real $\alpha_j$, $|K|^2=\sum_{i,j}\alpha_i\alpha_j\Re(\overline{\tilde k_i}\tilde k_j)$ and $K^2=\sum_{i,j}\alpha_i\alpha_j\tilde k_i\tilde k_j$;
insert into \eqref{qstar} and integrate. $\mathsf G\succeq0$ because $q_T^*\ge0$.
(c) Then $\sum_j\alpha_j\tilde k_j=\sum_\mu(\alpha_{2\mu-1}+i\alpha_{2\mu})k_\mu=\sum_\mu y_\mu k_\mu$, where $\alpha$ corresponds to $y$ under \eqref{eq:ident}.
\end{proof}

\subsection{The reliability theorem}
The certificate controls the Taylor remainders by the \emph{remainder density}
\begin{equation}\label{eq:rho}
\rho(z):=\sum_{\text{cells }C:\ z\in D(a_C,R_C)}|w_C|\,(1+\eta^{-1})\,\theta_C^{\,2p+4}\qquad(z\in\Om).
\end{equation}
In the next two lemmas and in Theorem~\ref{thm:rel}, $\mathcal R$ denotes either the Gram form \eqref{gram} of representers \eqref{rep}, regarded as a form on $\R^{2N}$,
or the Gram form $\mathsf G$ \eqref{eq:G} of real representers \eqref{eq:tildek}; and $\mathcal R_{\rm up}$ is a real symmetric $2N\times2N$ matrix with
$\mathcal R_{\rm up}\succeq\mathcal R$.
\begin{lemma}[Taylor data are controlled by $\|f\|_q$]\label{lem:step1}
For every $f\in\Az$ there is $\alpha\in\R^{2N}$ with
\begin{equation}\label{eq:rel-b}
\ell(f)=\mathcal R_{\rm up}\alpha,\qquad\alpha^{\mathsf T}\mathcal R_{\rm up}\alpha\le\|f\|_q^2 .
\end{equation}
\end{lemma}
\begin{proof}
Let $y\in\R^{2N}$. For complex representers let $y$ also denote the corresponding vector of $\C^N$ \eqref{eq:ident} and put $K:=\sum_\mu y_\mu k_\mu$; since
$\langle f,y_\mu k_\mu\rangle=\bar y_\mu\ell_\mu(f)$, \eqref{eq:ident} gives $y\cdot\ell(f)=\Re(y^*\ell(f))=\Re\langle f,K\rangle$. For real representers put
$K:=\sum_jy_j\tilde k_j$; then $y\cdot\ell(f)=\Re\langle f,K\rangle$ by Lemma~\ref{lem:reallin}(a). In both cases $\int_\Om q^*(K)=\mathcal R(y)$, so
\eqref{dual} and $\mathcal R\preceq\mathcal R_{\rm up}$ give
\[(y\cdot\ell(f))^2\le\|f\|_q^2\,\mathcal R(y)\le\|f\|_q^2\,y^{\mathsf T}\mathcal R_{\rm up}y\qquad(y\in\R^{2N}).\]
$\mathcal R_{\rm up}$ is positive semidefinite because $\mathcal R_{\rm up}\succeq\mathcal R\succeq0$ (Lemmas~\ref{lem:grammatrix} and~\ref{lem:reallin}). The
finite-dimensional duality Lemma~\ref{lem:findim}, with $S:=\mathcal R_{\rm up}$, $x:=\ell(f)$ and $t:=\|f\|_q^2$, gives $\alpha$ (e.g.\ $\alpha=S^+x$, $S^+$ the
Moore--Penrose pseudoinverse) with \eqref{eq:rel-b}.
\end{proof}
\begin{remark}
Neither invertibility of $\mathcal R_{\rm up}$ nor a bound on its condition number is needed in Lemma~\ref{lem:step1}: directions in the kernel of
$\mathcal R_{\rm up}$ are exactly those in which $y\cdot\ell(f)=0$ for all $f\in\Az$.
\end{remark}
\begin{lemma}\label{lem:epsle1}
Suppose that there is at least one cell, let $Q$ be given by \eqref{eq:Q} with $\eta>0$, and let $\varepsilon\in\R$. If
$Q(\mathcal R_{\rm up}\alpha)\le(1-\varepsilon)\,\alpha^{\mathsf T}\mathcal R_{\rm up}\alpha$ for all $\alpha\in\R^{2N}$, then $\varepsilon\le1$.
\end{lemma}
\begin{proof}
Suppose $\varepsilon>1$, and let $f\in\Az$. Lemma~\ref{lem:step1} applies to $f$ and to $if\in\Az$: $x:=\ell(f)=\mathcal R_{\rm up}\alpha$ and $\ell(if)=ix=\mathcal R_{\rm up}\alpha'$
with $\alpha^{\mathsf T}\mathcal R_{\rm up}\alpha,\ \alpha'^{\mathsf T}\mathcal R_{\rm up}\alpha'\le\|f\|_q^2$. In \eqref{eq:Q} the indefinite parts of $Q(x)$ and $Q(ix)$ cancel:
$Q(x)+Q(ix)=2\eta\sum_{\text{cells}}|w_C|\int_C|\sum_{\mu\in J_C}x_\mu\phi_\mu|^2\ge0$. By hypothesis,
$0\le Q(x)+Q(ix)\le(1-\varepsilon)(\alpha^{\mathsf T}\mathcal R_{\rm up}\alpha+\alpha'^{\mathsf T}\mathcal R_{\rm up}\alpha')\le0$, so $\alpha^{\mathsf T}\mathcal R_{\rm up}\alpha=0$;
as $\mathcal R_{\rm up}\succeq0$, this gives $\mathcal R_{\rm up}\alpha=0$, i.e.\ $\ell(f)=0$. Thus $\ell$ vanishes on $\Az$; in particular, for $\mu$ with $n_\mu=0$,
$f(a_\mu)=\ell_\mu(f)=0$ for all $f\in\Az$, which contradicts Lemma~\ref{lem:nonvanish}.
\end{proof}

\begin{theorem}[reliability]\label{thm:rel}
Let $\mathcal R$ be the Gram form \eqref{gram} of representers \eqref{rep}, regarded as a form on $\R^{2N}$, or the Gram form $\mathsf G$ \eqref{eq:G} of real
representers \eqref{eq:tildek}, and let $\mathcal R_{\rm up}\succeq\mathcal R$ be a real symmetric $2N\times2N$ matrix. Let $Q$ be given by \eqref{eq:Q} with $\eta>0$, and let
$\varepsilon\ge\sup_{z\in\Om}\rho(z)/c(z)$, with $\rho$ from \eqref{eq:rho} and $c$ from \eqref{eq:cfun}. If
\begin{equation}\label{crit}Q(\mathcal R_{\rm up}\alpha)\le(1-\varepsilon)\,\alpha^{\mathsf T}\mathcal R_{\rm up}\alpha\qquad\text{for all }\alpha\in\R^{2N},\end{equation}
then $\kappa(\Om)\le\kappa_0$. If moreover $\kappa_0<1$, then $\beta(\Om)\ge\sqrt{(1-\kappa_0)/2}>0$.
\end{theorem}
\begin{proof}
It suffices to show
\begin{equation}\label{eq:rel-a}
\sum_{T\in\mathcal A}\Re\Big(w_T\!\int_Tf^2\Big)\le\|f\|_q^2\qquad(f\in\Az).
\end{equation}
By Proposition~\ref{prop:goal}(b), \eqref{eq:rel-a} gives $\kappa(\Om)\le\kappa_0$, and Corollary~\ref{cor:bounds}(a) gives the bound for $\beta$. (Since $\kappa\le1$ always, by the Cauchy--Schwarz inequality,
only $\kappa_0<1$ is of interest.) If there are no active triangles, \eqref{eq:rel-a} is trivial. Otherwise fix $f\in\Az$.

\emph{Step 1: the Taylor part.} By Lemma~\ref{lem:step1} there is $\alpha$ with \eqref{eq:rel-b}. By \eqref{crit}, and since $\varepsilon\le1$ by Lemma~\ref{lem:epsle1},
\begin{equation}\label{eq:rel-c}
Q(\ell(f))=Q(\mathcal R_{\rm up}\alpha)\le(1-\varepsilon)\,\alpha^{\mathsf T}\mathcal R_{\rm up}\alpha\le(1-\varepsilon)\|f\|_q^2 .
\end{equation}

\emph{Step 2: one cell.} Let $C$ be a cell, $t:=t_C$, $r:=f-t$, and $\|\cdot\|:=\|\cdot\|_{L^2(C)}$. Since $f^2=t^2+2tr+r^2$, the Cauchy--Schwarz inequality gives
\begin{align*}\Re\Big(w_C\!\int_Cf^2\Big)&=\Re\Big(w_C\!\int_Ct^2\Big)+\Re\Big(w_C\!\int_C(2tr+r^2)\Big)\\&\le\Re\Big(w_C\!\int_Ct^2\Big)+|w_C|\big(2\|t\|\|r\|+\|r\|^2\big).\end{align*}
By Young's inequality $2\|t\|\|r\|\le\eta\|t\|^2+\eta^{-1}\|r\|^2$ (expand $(\sqrt\eta\|t\|-\|r\|/\sqrt\eta)^2\ge0$), so
\[\Re\Big(w_C\!\int_Cf^2\Big)\le\Re\Big(w_C\!\int_Ct_C^2\Big)+\eta|w_C|\,\|t_C\|_{L^2(C)}^2+|w_C|(1+\eta^{-1})\|f-t_C\|_{L^2(C)}^2 .\]

\emph{Step 3: summation.} The cells partition the active triangles up to null sets, and $w_C=w_T$ for $C\subset T$; so
$\sum_{T\in\mathcal A}\Re(w_T\int_Tf^2)$ is the sum of $\Re(w_C\int_Cf^2)$ over the cells. Summing Step~2 over the cells, the terms
$\Re(w_C\int_Ct_C^2)+\eta|w_C|\|t_C\|_{L^2(C)}^2$ add up to $Q(\ell(f))$ by \eqref{eq:Qell}. By Lemma~\ref{lem:taylor},
$\|f-t_C\|_{L^2(C)}^2\le\theta_C^{2p+4}\int_\Om\mathbf 1_{D(a_C,R_C)}|f|^2$; interchanging the finite sum with the integral, the remainder terms add up to at most
$\int_\Om\rho|f|^2$, with $\rho$ from \eqref{eq:rho}. Pointwise $\rho\le\varepsilon c$ and $c|f|^2\le q(f)$, so $\int\rho|f|^2\le\varepsilon\int c|f|^2\le\varepsilon\|f\|_q^2$.
With \eqref{eq:rel-c},
\[\sum_{T\in\mathcal A}\Re\Big(w_T\!\int_Tf^2\Big)\le Q(\ell(f))+\int_\Om\rho|f|^2\,dA\le(1-\varepsilon)\|f\|_q^2+\varepsilon\|f\|_q^2=\|f\|_q^2,\]
which is \eqref{eq:rel-a}.
\end{proof}
Condition \eqref{crit} is a finite-dimensional matrix inequality whose data are integrals of explicit functions. No assumption on the spectrum or on the
quality of the discretisation enters: poor correctors or a low Taylor degree only make \eqref{crit} fail.
\begin{remark}[scalar weight]\label{rem:scalar}
Since $q_T(\zeta)\ge c\,|\zeta|^2$, Theorem~\ref{thm:rel} remains true if $\|f\|_q^2$ is replaced by $\int c|f|^2$ and $\mathcal R$ by the Hermitian form with matrix
$\int c^{-1}\bar k_\mu k_\nu$. This variant is simpler, but it discards the part $|w_T||f|^2-\Re(w_Tf^2)\ge0$ on passive triangles and is no longer equivalent to
$\kappa\le\kappa_0$; see Remark~\ref{rem:rings}.
\end{remark}

\section{Convergence}\label{sec:conv}
Because \eqref{goal} is an identity, the only losses in Theorem~\ref{thm:rel} come from the discretisation of the compact part. They can be made arbitrarily small.
Throughout this section the standing assumptions of Section~\ref{sec:cert} hold, the representers are real-linear (Definition~\ref{def:reallin}), $|\cdot|$ is the
Euclidean norm on $\R^{2N}$, and $\|\cdot\|$ is the spectral norm of a quadratic form on $\R^{2N}$. We regard $L^2(\Om)=L^2(\Om;\C)$ as a real Hilbert space with the
inner product $\langle\cdot,\cdot\rangle_q$; its norm $\|\cdot\|_q$ is equivalent to the $L^2$ norm, and $\Az$ is a closed subspace (Lemma~\ref{lem:bergman}(b)).

\subsection{Structure of the Gram form}
\begin{lemma}[the map $\sharp$]\label{lem:sharp}
For $k\in L^2(\Om)$ let
\begin{equation}\label{eq:sharp}
k^\sharp:=A_T^{-1}k=\hat\kappa_Tk+\overline{\hat w_T}\,\bar k\quad\text{on }T\in\mathcal T,
\end{equation}
i.e.\ $k^\sharp=(\kappa_0k+\bar w_T\bar k)/(\kappa_0^2-|w_T|^2)$ on $\mathcal P$ and $k^\sharp=k/\kappa_0$ on $\mathcal A$. Then
\[\langle f,k^\sharp\rangle_q=\Re\!\int_\Om f\bar k\,dA,\qquad\|k^\sharp\|_q^2=\int_\Om q^*(k)\,dA\qquad(f,k\in L^2(\Om)),\]
and $\sharp$ is a real-linear bijection of $L^2(\Om)$ with inverse $\mathsf Ae:=\kappa_0e-\bar w_T\bar e$ on $\mathcal P$ and $\mathsf Ae:=\kappa_0e$ on $\mathcal A$.
\end{lemma}
\begin{proof}
The second expression in \eqref{eq:sharp} is Lemma~\ref{lem:dualform}(b). Since $\langle f,g\rangle_q=\sum_T\int_Tf\cdot A_Tg\,dA$,
$\langle f,k^\sharp\rangle_q=\sum_T\int_Tf\cdot k\,dA=\Re\int_\Om f\bar k$, and $\|k^\sharp\|_q^2=\sum_T\int_T(A_T^{-1}k)\cdot k\,dA=\int_\Om q^*(k)$ by
Lemma~\ref{lem:dualform}(b). The pointwise inverse is $e\mapsto A_Te$, which is $\kappa_0e-W(w_T)e=\kappa_0e-\bar w_T\bar e$ on $\mathcal P$ by
Lemma~\ref{lem:dualform}(a). Both maps are bounded on $L^2$ because $|w_T|\le\kappa_{\mathcal P}<\kappa_0$.
\end{proof}
\begin{lemma}[splitting of the Gram form]\label{lem:split}
Let $\tilde k_1,\dots,\tilde k_{2N}\in L^2(\Om)$ be real representers of $\lambda_1,\dots,\lambda_{2N}$, let $P_0$ be the $\langle\cdot,\cdot\rangle_q$-orthogonal projection of
$L^2(\Om)$ onto $\Az$, and for $\alpha\in\R^{2N}$ let
\begin{equation}\label{eq:fe}
\tilde K_\alpha:=\sum_{j=1}^{2N}\alpha_j\tilde k_j,\qquad f_\alpha:=P_0\tilde K_\alpha^\sharp,\qquad e_\alpha:=\tilde K_\alpha^\sharp-f_\alpha .
\end{equation}
\begin{enumerate}
\item[(a)] $f_\alpha$ is the unique element of $\Az$ with $\langle f,f_\alpha\rangle_q=\alpha\cdot\ell(f)$ for all $f\in\Az$; in particular, it does not depend on the choice of
the real representers.
\item[(b)] $\mathsf G=\mathcal R^0+\mathcal E$ with the positive semidefinite quadratic forms
\begin{equation}\label{eq:R0E}
\mathcal R^0(\alpha):=\|f_\alpha\|_q^2,\qquad\mathcal E(\alpha):=\|e_\alpha\|_q^2\qquad(\alpha\in\R^{2N}),
\end{equation}
and $\mathcal R^0$ does not depend on the choice of the real representers.
\item[(c)] $\ell(f_\alpha)=\mathcal R^0\alpha$ for all $\alpha\in\R^{2N}$.
\end{enumerate}
\end{lemma}
\begin{proof}
$f_\alpha$ and $e_\alpha$ depend linearly on $\alpha$, so \eqref{eq:R0E} defines quadratic forms.
(a) For $f\in\Az$, $\langle f,e_\alpha\rangle_q=0$, so by Lemmas~\ref{lem:sharp} and~\ref{lem:reallin}(a)
$\langle f,f_\alpha\rangle_q=\langle f,\tilde K_\alpha^\sharp\rangle_q=\Re\langle f,\tilde K_\alpha\rangle=\alpha\cdot\ell(f)$. Two elements of $\Az$ with this property differ by an
element of $\Az$ orthogonal to $\Az$.
(b) By Lemma~\ref{lem:sharp} and Pythagoras, $\mathsf G(\alpha)=\int q^*(\tilde K_\alpha)=\|\tilde K_\alpha^\sharp\|_q^2=\|f_\alpha\|_q^2+\|e_\alpha\|_q^2$. By (a), $\mathcal R^0$
does not depend on the representers.
(c) For $\alpha'\in\R^{2N}$, (a) with $f=f_\alpha$ and $\alpha'$ in place of $\alpha$ gives $\alpha'\cdot\ell(f_\alpha)=\langle f_\alpha,f_{\alpha'}\rangle_q=\alpha'^{\mathsf T}\mathcal R^0\alpha$,
the last step by polarisation of $\mathcal R^0$.
\end{proof}
\begin{lemma}[optimal correctors]\label{lem:optcorr}
For $j=1,\dots,2N$ there are $\chi_j^*\in H^1_0(\Om;\C)$ and $\gamma_j^*\in\C$ such that $\tilde k_j^*:=\hat b_j-\dz\chi_j^*-\gamma_j^*$ is a real representer of $\lambda_j$ and
$(\tilde k_j^*)^\sharp\in\Az$. For the real representers $\tilde k_1^*,\dots,\tilde k_{2N}^*$ one has $\mathcal E=0$ in Lemma~\ref{lem:split}.
\end{lemma}
\begin{proof}
Let $e_j^b:=\hat b_j^\sharp-P_0\hat b_j^\sharp$ be the component of $\hat b_j^\sharp$ that is $\langle\cdot,\cdot\rangle_q$-orthogonal to $\Az$, and $h_j:=\mathsf Ae_j^b$. For
$f\in\Az$, Lemma~\ref{lem:sharp} gives $\Re\langle f,h_j\rangle=\langle f,e_j^b\rangle_q=0$. Applying this to $if\in\Az$ gives $0=\Re\langle if,h_j\rangle=-\Im\langle f,h_j\rangle$.
Hence $\langle f,h_j\rangle=0$ for all $f$ in the complex space $\Az$, and $h_j\in(\Az)^\perp=\dz H^1_0(\Om)\oplus\C$ (Lemma~\ref{lem:decomp}). Write
$h_j=\dz\chi_j^*+\gamma_j^*$ with $\chi_j^*\in H^1_0(\Om;\C)$ and $\gamma_j^*\in\C$. Then $\tilde k_j^*=\hat b_j-h_j$, so $\Re\langle f,\tilde k_j^*\rangle=\Re\langle f,\hat b_j\rangle=\lambda_j(f)$
for $f\in\Az$ (Lemma~\ref{lem:bump}), and $(\tilde k_j^*)^\sharp=\hat b_j^\sharp-e_j^b=P_0\hat b_j^\sharp\in\Az$. As $\sharp$ is real-linear, $\tilde K_\alpha^\sharp\in\Az$ for every
$\alpha\in\R^{2N}$, i.e.\ $e_\alpha=0$.
\end{proof}
The $\chi_j^*$ need not be Lipschitz, so the $\tilde k_j^*$ need not be admissible in \eqref{eq:tildek}. They are approximated by finite elements.
\begin{lemma}[finite element correctors]\label{lem:fecorr}
Let $\chi_j^*,\gamma_j^*$ be as in Lemma~\ref{lem:optcorr}. Let $V_h\subset H^1_0(\Om)$, $h\to0$, be spaces of continuous piecewise polynomials vanishing on $\partial\Om$
with $\inf_{v\in V_h}|\chi-v|_1\to0$ for every $\chi\in H^1_0(\Om)$ (for instance Lagrange elements on shape-regular meshes). Choose $\chi_{j,h}\in V_h+iV_h$ with
$|\chi_j^*-\chi_{j,h}|_1\to0$ and put $\tilde k_{j,h}:=\hat b_j-\dz\chi_{j,h}-\gamma_j^*$. Then the $\tilde k_{j,h}$ are admissible in \eqref{eq:tildek}, the form $\mathcal R^0$ of
Lemma~\ref{lem:split} is the same for all $h$, and the corrector error $\mathcal E_h$ of the $\tilde k_{j,h}$ satisfies
\[\|\mathcal E_h\|\le\frac1{4c_{\min}}\sum_{j=1}^{2N}|\chi_j^*-\chi_{j,h}|_1^2\longrightarrow0\qquad(h\to0).\]
\end{lemma}
\begin{proof}
Elements of $V_h$ are Lipschitz and vanish on $\partial\Om$, so the $\tilde k_{j,h}$ are admissible; they are real representers by Lemma~\ref{lem:reallin}(a), and
$\mathcal R^0$ does not depend on them by Lemma~\ref{lem:split}(b). Write $\tilde k_{j,h}=\tilde k_j^*+\dz(\chi_j^*-\chi_{j,h})$. As $\dz(\chi_j^*-\chi_{j,h})\in(\Az)^\perp$
(Lemma~\ref{lem:decomp}), its image under $\sharp$ is $\langle\cdot,\cdot\rangle_q$-orthogonal to $\Az$ by Lemma~\ref{lem:sharp}, while $(\tilde k_j^*)^\sharp\in\Az$. Hence
$e_{\alpha,h}=\big(\sum_j\alpha_j\dz(\chi_j^*-\chi_{j,h})\big)^\sharp$ and, by Lemma~\ref{lem:sharp},
\[\mathcal E_h(\alpha)=\int_\Om q^*\Big(\sum_j\alpha_j\dz(\chi_j^*-\chi_{j,h})\Big)dA\le\frac{|\alpha|^2}{c_{\min}}\sum_j\|\dz(\chi_j^*-\chi_{j,h})\|^2=\frac{|\alpha|^2}{4c_{\min}}\sum_j|\chi_j^*-\chi_{j,h}|_1^2 .\]
Here we used $q_T^*(\xi)\le|\xi|^2/(\kappa_0-|w_T|)\le|\xi|^2/c_{\min}$ (the largest eigenvalue of $A_T^{-1}$), the Cauchy--Schwarz inequality in $\R^{2N}$, and
$\|\dz\varphi\|=\frac12|\varphi|_1$ (Lemma~\ref{lem:decomp}). As $\mathcal E_h\succeq0$, this bounds $\|\mathcal E_h\|$.
\end{proof}
\begin{lemma}[Taylor data versus the active form]\label{lem:Qapprox}
Let the cells satisfy $\theta_C=\frac12$ for every cell $C$, let $L$ be the maximal number of disks $D(a_C,R_C)$ that contain a common point, and
$w_{\max}:=\max_{T\in\mathcal T}|w_T|$. Then for all $\eta>0$, $p\ge0$ and $f\in\Az$,
\begin{equation}\label{eq:B-S1}
\begin{gathered}Q(\ell(f))\le\sum_{T\in\mathcal A}\Re\Big(w_T\!\int_Tf^2\Big)+\big(\Lambda_12^{-p}+\Lambda_2\eta\big)\|f\|_q^2,\\ \Lambda_1:=\frac{5w_{\max}(1+L)}{16\,c_{\min}},\quad\Lambda_2:=\frac{w_{\max}(16+L)}{8\,c_{\min}}.\end{gathered}
\end{equation}
\end{lemma}
\begin{proof}
Let $C$ be a cell, $t:=t_C$, $r:=f-t$, $\|\cdot\|:=\|\cdot\|_{L^2(C)}$ and $D_C:=D(a_C,R_C)$. From $t^2=f^2-2fr+r^2$, the Cauchy--Schwarz inequality, and
$\|t\|^2\le2\|f\|^2+2\|r\|^2$,
\[\Re\Big(w_C\!\int_Ct^2\Big)+\eta|w_C|\|t\|^2\le\Re\Big(w_C\!\int_Cf^2\Big)+|w_C|\big(2\|f\|\|r\|+\|r\|^2\big)+2\eta|w_C|\big(\|f\|^2+\|r\|^2\big).\]
By Lemma~\ref{lem:taylor} and $\theta_C=\frac12$, $\|r\|\le2^{-p-2}\|f\|_{D_C}$; and $2\|f\|\|f\|_{D_C}\le\|f\|^2+\|f\|_{D_C}^2$. Now sum over the cells, using
\eqref{eq:Qell} on the left. The cells are disjoint up to null sets and partition the active triangles, so $\sum_C\|f\|_{L^2(C)}^2\le\|f\|^2$ and
$\sum_C\Re(w_C\int_Cf^2)=\sum_{T\in\mathcal A}\Re(w_T\int_Tf^2)$; every point lies in at most $L$ disks, so $\sum_C\|f\|_{D_C}^2\le L\|f\|^2$. Hence
\[Q(\ell(f))\le\sum_{T\in\mathcal A}\Re\Big(w_T\!\int_Tf^2\Big)+w_{\max}\Big[2^{-p-2}(1+L)+2^{-2p-4}L+2\eta\big(1+2^{-2p-4}L\big)\Big]\|f\|^2 .\]
Using $2^{-2p-4}\le2^{-p}/16$ and $\|f\|^2\le\|f\|_q^2/c_{\min}$ gives \eqref{eq:B-S1}.
\end{proof}

\subsection{Existence of certificates}
\begin{theorem}\label{thm:B}
Let the standing assumptions of Section~\ref{sec:cert} hold (in particular every triangle meeting $\partial\Om$ is passive), and let $\kappa_0>\max\{\kappa(\Om),\kappa_{\mathcal P}\}$.
Then there exist cells with $\theta_C\le\frac12$, $\eta>0$, a degree $p$, a number $\varepsilon\ge\sup_\Om\rho/c$, and real representers \eqref{eq:tildek}, with one
corrector $\chi_j$ from a finite element space on a sufficiently fine mesh and one constant $\gamma_j$ for each real Taylor functional $\lambda_j$, $j=1,\dots,2N$,
such that \eqref{crit} holds for $\mathcal R_{\rm up}=\mathcal R=\mathsf G$, the Gram form \eqref{eq:G} of these representers.
\end{theorem}
\begin{proof}
\emph{Step 0: cells.} An active triangle $T$ does not meet $\partial\Om$, so it is a compact subset of $\Om$ and $d_T:=\dist(T,\partial\Om)>0$. Subdivide $T$ into cells of
diameter at most $d_T/2$; for each cell $C$ let $a_C$ be one of its vertices, $r_C:=\operatorname{diam}C$ and $R_C:=2r_C$. Then $C\subset\overline{D(a_C,r_C)}$,
$D(a_C,R_C)\subset D(a_C,d_T)\subset\Om$, and $\theta_C=\frac12$. Fix these cells, and let $L$ and $w_{\max}$ be as in Lemma~\ref{lem:Qapprox}.

\emph{Step 1: choice of $\eta$ and $p$.} Let $\sigma:=(\kappa_0-\kappa(\Om))/(\kappa_0+w_{\max})\in(0,1)$ be the relative gap. For $f\in\Az$, \eqref{goal}, the
inequality $\Re\int f^2\le\kappa\|f\|^2$ (by \eqref{eq:kappa-re}) and $\|f\|_q^2\le(\kappa_0+w_{\max})\|f\|^2$ give
\begin{multline*}\sum_{T\in\mathcal A}\Re\Big(w_T\!\int_Tf^2\Big)=\|f\|_q^2+\Re\!\int_\Om f^2-\kappa_0\|f\|^2\le\|f\|_q^2-(\kappa_0-\kappa)\|f\|^2\\ \le\Big(1-\frac{\kappa_0-\kappa}{\kappa_0+w_{\max}}\Big)\|f\|_q^2=(1-\sigma)\|f\|_q^2 .\end{multline*}
Choose $\eta>0$ with $\Lambda_2\eta\le\sigma/4$, and then $p$ so large that $\Lambda_12^{-p}\le\sigma/4$ and
\[\varepsilon:=(1+\eta^{-1})\,L\,w_{\max}\,2^{-2p-4}/c_{\min}\le\sigma/8 .\]
Then Lemma~\ref{lem:Qapprox} gives
\begin{equation}\label{eq:B-S2}
Q(\ell(f))\le\big(1-\tfrac\sigma2\big)\|f\|_q^2\qquad(f\in\Az).
\end{equation}
Since $\rho(z)\le Lw_{\max}(1+\eta^{-1})2^{-2p-4}$ and $c\ge c_{\min}$, this $\varepsilon$ satisfies $\varepsilon\ge\sup\rho/c$, as required in Theorem~\ref{thm:rel}. From now on the
cells, $\eta$ and $p$, and hence $N$, $\ell$, $Q$ and $\varepsilon$, are fixed. Let $\mathsf q:=\|Q\|$.

\emph{Step 2: correctors.} For every family of real representers, $\mathsf G=\mathcal R^0+\mathcal E$ with $\mathcal R^0$ independent of the family
(Lemma~\ref{lem:split}). Take the finite element representers $\tilde k_{j,h}$ of Lemma~\ref{lem:fecorr}; their corrector error $\mathcal E=\mathcal E_h$ satisfies
$\|\mathcal E\|\to0$ as $h\to0$, while $N$, $Q$ and $\mathcal R^0$ stay fixed.

\emph{Step 3: verification of \eqref{crit}.} Let $\alpha\in\R^{2N}$, $x:=\mathcal R^0\alpha=\ell(f_\alpha)$ (Lemma~\ref{lem:split}(c)) and $y:=\mathcal E\alpha$, so that
$\mathsf G\alpha=x+y$. For a positive semidefinite matrix $S$, $|S\alpha|^2=\alpha^{\mathsf T}S^2\alpha\le\|S\|\,\alpha^{\mathsf T}S\alpha$; hence
$|x|^2\le\|\mathcal R^0\|\,\|f_\alpha\|_q^2$ and $|y|^2\le\|\mathcal E\|\,\|e_\alpha\|_q^2$. Since $f_\alpha\in\Az$, \eqref{eq:B-S2} gives $Q(x)\le(1-\frac\sigma2)\|f_\alpha\|_q^2$.
With $Q(x+y)=Q(x)+2Q(x,y)+Q(y)$, $|Q(x,y)|\le\mathsf q|x||y|$, $|Q(y)|\le\mathsf q|y|^2$ and $2ab\le a^2+b^2$,
\[Q(\mathsf G\alpha)\le\big(1-\tfrac\sigma2\big)\|f_\alpha\|_q^2+\tau\big(\|f_\alpha\|_q^2+\|e_\alpha\|_q^2\big)+\mathsf q\|\mathcal E\|\,\|e_\alpha\|_q^2,\qquad\tau:=\mathsf q\big(\|\mathcal R^0\|\,\|\mathcal E\|\big)^{1/2}.\]
Since $\alpha^{\mathsf T}\mathsf G\alpha=\|f_\alpha\|_q^2+\|e_\alpha\|_q^2$, \eqref{crit} holds as soon as $1-\frac\sigma2+\tau\le1-\varepsilon$ and $\tau+\mathsf q\|\mathcal E\|\le1-\varepsilon$.
Because $\varepsilon\le\sigma/8$, both hold if $\tau\le\sigma/4$ and $\tau+\mathsf q\|\mathcal E\|\le1-\sigma/8$; by Step~2 this is the case on a sufficiently fine mesh.
\end{proof}

\subsection{Convergence of the bounds}
\begin{lemma}[existence of good multipliers]\label{lem:goodmult}
Let $\Om$ be a polygonal domain with corners $c_j$ of opening $\omega_j$, and let $\kappa_1>\kess$. There are a conforming triangulation $\mathcal T$ of
$\overline\Om$, a continuous piecewise linear multiplier $m$ on $\mathcal T$, and a set $\mathcal P\subset\mathcal T$ containing all triangles that meet $\partial\Om$,
such that $\kappa_{\mathcal P}=\max_{T\in\mathcal P}|w_T|<\kappa_1$.
\end{lemma}
\begin{proof}
\emph{Choice of the fan parameters.} By Lemma~\ref{lem:fan}(b), $\kappa_{n,\gamma}(\omega_j)\to|\sin\omega_j|/\omega_j\le\kess<\kappa_1$ as $n\to\infty$, $\gamma\to0$.
Choose $n\ge2$ and $\gamma$ with $0<\gamma<\min_j\omega_j/2$ such that $\max_j\kappa_{n,\gamma}(\omega_j)<\kappa_1$. Since $\omega_j<2\pi$ and $n\ge2$,
$(\omega_j-2\gamma)/n<\pi$, so the hypotheses of Lemma~\ref{lem:fan}(b) hold at every corner.

\emph{The multiplier in a boundary layer.} Choose $r_1>0$ so small that for every corner $c_j$, $\Om\cap D(c_j,2r_1)$ is a sector as in Lemma~\ref{lem:fan} (translated
to $c_j$), and the disks $D(c_j,2r_1)$ are pairwise disjoint. At $c_j$ let $F_j$ be the union of the triangles with vertices $c_j$, $c_j+r_1e^{i\vartheta_{l-1}}$,
$c_j+r_1e^{i\vartheta_l}$, one for each sector of Lemma~\ref{lem:fan}(b); each lies in a single sector, and on it let $m(z):=\tilde m(z-c_j)$, where $\tilde m$ is the fan
multiplier of Lemma~\ref{lem:fan}(b). (Since $\dzb$ commutes with translations, the weights are unchanged.) For a side $e$ joining $c_j$ to $c_{j'}$ let $Z_e$ be
the trapezoid bounded by $e$, by the chords of the two edge-sector triangles of $F_j$ and $F_{j'}$ adjacent to $e$, and by the segment parallel to $e$ at distance
$r_1\sin\gamma$ that joins the endpoints of these chords. On $Z_e$ let $m$ be the reflection multiplier of $e$. By Lemma~\ref{lem:fan}(b), on the edge-sector triangles
of $F_j$ the fan multiplier is the reflection multiplier of the adjacent side; so on the common chord both definitions give the same affine function. Hence $m$ is
well defined and continuous on the boundary layer $\mathcal B:=\bigcup_jF_j\cup\bigcup_eZ_e$, affine on each triangle of each $F_j$ and on each $Z_e$, and zero on
$\partial\Om$. For $r_1$ small, with $\gamma$ fixed, the trapezoids lie in $\overline\Om$ and meet the fans and each other only along chords and vertices. Then
$\mathcal B$ is a closed polygonal neighbourhood of $\partial\Om$ in $\overline\Om$; its inner boundary consists of the chords
$[c_j+r_1e^{i\vartheta_{l-1}},c_j+r_1e^{i\vartheta_l}]$ of the inner sectors and the inner edges of the trapezoids, and it has positive distance from $\partial\Om$.

\emph{Triangulation and passive set.} Split each $Z_e$ into two triangles by a diagonal. Triangulate the polygonal region $\overline{\Om\setminus\mathcal B}$
(multiply connected if $\Om$ is) without additional vertices on its boundary; its triangles do not meet $\partial\Om$. Together with the triangles of $\mathcal B$ this is a
conforming triangulation $\mathcal T$ of $\overline\Om$. Let $m$ be the continuous piecewise linear function on $\mathcal T$ with the nodal values already fixed on
$\mathcal B$ (the nodes of $\overline{\Om\setminus\mathcal B}$ are vertices of $\mathcal B$, or interior vertices where any value, e.g.\ $0$, may be taken). On the
triangles of $\mathcal B$ it coincides with the affine pieces above. Hence $m$ is Lipschitz and vanishes on $\partial\Om$. Let $\mathcal P$ be the set of
triangles of $\mathcal B$; it contains all triangles meeting $\partial\Om$. On the edge-sector triangles and on the trapezoids $w=0$, and on the inner fan triangles
$|w|=\kappa_{n,\gamma}(\omega_j)$. So $\kappa_{\mathcal P}=\max_j\kappa_{n,\gamma}(\omega_j)<\kappa_1$.
\end{proof}
\begin{corollary}[convergence]\label{cor:conv}
For every $\kappa_0$ with $\kappa(\Om)<\kappa_0<1$ there are a multiplier, a passive set and discrete data for which \eqref{crit} holds. Hence the certified lower bounds
$\sqrt{(1-\kappa_0)/2}$ converge to $\beta(\Om)$ along a suitable sequence of parameters.
\end{corollary}
\begin{proof}
The interval $(\kappa(\Om),1)$ is not empty, because $\kappa(\Om)=1-2\beta(\Om)^2<1$ for a bounded Lipschitz domain. By Lemma~\ref{lem:kess}, $\kess\le\kappa(\Om)<\kappa_0$, so
Lemma~\ref{lem:goodmult} with $\kappa_1:=\kappa_0$ gives $\mathcal T$, $m$ and $\mathcal P$ with $\kappa_{\mathcal P}<\kappa_0$, i.e.\ $\kappa_0>\max\{\kappa(\Om),\kappa_{\mathcal P}\}$.
Theorem~\ref{thm:B} provides discrete data for which \eqref{crit} holds, and Theorem~\ref{thm:rel} certifies $\beta(\Om)\ge\sqrt{(1-\kappa_0)/2}$. As
$\kappa_0\downarrow\kappa(\Om)$, these bounds increase to $\sqrt{(1-\kappa(\Om))/2}=\beta(\Om)$ by Theorem~\ref{thm:beta}.
\end{proof}
The upper bounds from trial functions converge as well; this is the usual Rayleigh--Ritz argument, which we include for completeness.
\begin{proposition}[convergence of the upper bounds]\label{prop:upperconv}
Let $V_1\subset V_2\subset\dots\subset\Az$ be finite-dimensional subspaces with dense union, and $\kappa_k:=\max_{0\ne f\in V_k}|\int f^2|/\|f\|^2$. Then
$\kappa_k\uparrow\kappa(\Om)$, so the upper bounds $\sqrt{(1-\kappa_k)/2}$ for $\beta$ decrease to $\beta(\Om)$.
\end{proposition}
\begin{proof}
The maximum exists because the unit sphere of $V_k$ is compact. The sequence $(\kappa_k)$ is nondecreasing since the $V_k$ are nested, and $\kappa_k\le\kappa(\Om)$. The
Rayleigh quotient is continuous on $\Az\setminus\{0\}$, because $|\int f^2-\int g^2|=|\int(f-g)(f+g)|\le\|f-g\|\,\|f+g\|$. Given $0\ne f\in\Az$, let $f_k$ be the
$L^2$-orthogonal projection of $f$ onto $V_k$; by density of $\bigcup_kV_k$ and nestedness, $f_k\to f$. Then $\lim_k\kappa_k\ge\lim_k|\int f_k^2|/\|f_k\|^2=|\int f^2|/\|f\|^2$;
taking the supremum over $f$, $\kappa_k\uparrow\kappa(\Om)$. The last statement follows from Theorem~\ref{thm:beta}.
\end{proof}
\begin{remark}[corner designs and the scalar weight]\label{rem:rings}
\emph{Designs.} Theorem~\ref{thm:B} puts no condition on the multiplier away from the passive set. In computations we use the fan profile of
Lemma~\ref{lem:fan}(a) on a few innermost rings of the corner mesh, declared passive, and outside them a multiplier that vanishes along the sides and is
concentrated, with large modulus, in a thin wedge along the bisector (Figure~\ref{fig:mesh}); then few cells are active. The proof of Theorem~\ref{thm:B} gives no
rate of convergence.

\emph{Scalar weight.} With the scalar weight of Remark~\ref{rem:scalar} the situation is different. That criterion yields
$\sum_{\mathcal A}\Re(w_T\int_Tf^2)\le\int c|f|^2$ on $\Az$, and since $\int c|f|^2=\kappa_0\|f\|^2-\sum_{\mathcal P}|w_T|\|f\|_T^2$ this is the inequality
$J(f)\le\kappa_0\|f\|^2$, with $J(f):=\sum_{\mathcal P}|w_T|\|f\|_T^2+\sum_{\mathcal A}\Re(w_T\int_Tf^2)$. Accordingly, every $\kappa_0$ for which that criterion holds satisfies
$\kappa_0\ge\sup_{\Az}J/\|f\|^2$, a number which may exceed $\max(\kappa,\kappa_{\mathcal P})$. A corner is invariant under dilations, so such an excess does not disappear when
the passive zone shrinks; it is observed for three passive rings (about $0.011$ on the square) and not for two.
\end{remark}

\section{Verified implementation}\label{sec:impl}
A run of the verifier with input $\kappa_0\in\mathbb Q$ either proves \eqref{crit}, and thereby $\kappa(\Om)\le\kappa_0$, or stops without a statement.
Floating-point arithmetic is used only as an oracle that proposes data; everything on the proof path is recomputed exactly or enclosed. The subsections below
follow the hypotheses of Theorem~\ref{thm:rel}: the domain and the mesh, the weights, cells and correctors, the enclosure of the Gram matrix, and the matrix
inequality \eqref{crit}; then the upper bounds for $\beta$, the formal verification, and the trust base.

\subsection{Exact data and geometric checks}
Polygon and mesh vertices are dyadic rationals; nodal values of $m$ are Gaussian dyadic rationals. All tests in this subsection and in the next one are
exact integer computations.

\emph{Boundary.} The boundary is given as closed polygons $\Gamma_0,\dots,\Gamma_J$. The verifier checks that each $\Gamma_k$ is simple
(non-adjacent sides do not meet, adjacent sides meet only in their common vertex), that $\Gamma_0$ has positive and $\Gamma_1,\dots,\Gamma_J$
have negative signed area, that no two of the $\Gamma_k$ meet, and, by a crossing-number test at one vertex (which suffices because the curves
are disjoint), that each $\Gamma_k$, $k\ge1$, lies inside $\Gamma_0$ and outside every $\Gamma_l$, $l\ne0,k$. By the Jordan curve theorem,
$\Om$, the interior of $\Gamma_0$ minus the closed regions $\overline{H_k}$ bounded by $\Gamma_k$ ($k\ge1$), is then a connected polygonal
domain; the components of $\C\setminus\partial\Om$ are $\Om$, the $H_k$ and the exterior of $\Gamma_0$; $\Om$ lies to the left of every
directed side; and $2|\Om|$ is the sum of the doubled signed areas of the $\Gamma_k$.

\emph{Tiling.} For the triangles of the mesh the verifier checks:
(T1) every triangle is positively oriented; (T2) no directed edge, i.e.\ ordered pair of vertex indices, occurs in two triangles; (T3) every directed edge
$(i,j)$ whose reverse $(j,i)$ does not occur lies on a side of the boundary and has the direction of that side; (T4) the doubled areas of the triangles add up
to $2|\Om|$.
\begin{lemma}[tiling certificate]\label{lem:tiling}
Let $\Om$ be bounded by $\Gamma_0,\dots,\Gamma_J$ as above, and let a finite set of triangles satisfy (T1)--(T4). Then the closed triangles cover $\overline\Om$, their
interiors are pairwise disjoint and contained in $\Om$, the mesh is conforming, and the unpaired directed edges cover $\partial\Om$ exactly once.
\end{lemma}
\begin{proof}
For $z$ on no edge let $\operatorname{ind}(z)$ be the number of triangles containing $z$. By (T1), $\operatorname{ind}(z)$ is the winding number about $z$ of the
$1$-chain $c_{\mathcal T}:=\sum_T\partial T$. In $c_{\mathcal T}$ the edges $(i,j)$ and $(j,i)$ cancel, so by (T2) and (T3) $c_{\mathcal T}$ is a sum of segments on
$\partial\Om$, each directed like the side containing it. Hence $\operatorname{ind}$ extends to a locally constant function on $\C\setminus\partial\Om$: it
equals a constant $\operatorname{ind}_\Om$ on $\Om$, a constant $\operatorname{ind}_k$ on each $H_k$, and $0$ on the exterior of $\Gamma_0$. Where $\partial\Om$ is
crossed at a point that is not a vertex, from $\Om$ to the other side, the winding number drops by the number of segments of $c_{\mathcal T}$ through that point,
because all of them have $\Om$ on their left; hence $\operatorname{ind}_\Om\ge\operatorname{ind}_k\ge0$. By (T4),
$|\Om|=\int\operatorname{ind}\,dA=\operatorname{ind}_\Om|\Om|+\sum_k\operatorname{ind}_k|H_k|$. If $\operatorname{ind}_\Om=0$, all $\operatorname{ind}_k$ vanish and the
right-hand side is $0$; so $\operatorname{ind}_\Om\ge1$, and then $\operatorname{ind}_\Om=1$ and all $\operatorname{ind}_k=0$. Thus $\operatorname{ind}=1$ almost
everywhere in $\Om$ and $\operatorname{ind}=0$ almost everywhere outside $\overline\Om$; the interiors of the triangles are pairwise disjoint and contained in $\Om$,
and their union, being closed and containing almost every point of $\Om$, is $\overline\Om$. The mesh is conforming: an edge not on $\partial\Om$ is shared, with
both end points, by exactly one other triangle, and $\operatorname{ind}=1$ near its interior points leaves no room for a vertex of a third triangle. The unpaired
edges cover $\partial\Om$ exactly once, because $\operatorname{ind}_\Om-\operatorname{ind}_k=1$ and $\operatorname{ind}_\Om-0=1$.
\end{proof}
The argument uses neither convexity nor pairwise intersection tests between triangles, and it covers multiply connected $\Om$. By Lemma~\ref{lem:tiling},
a function given by nodal values and affine on each triangle is continuous, and every mesh vertex on $\partial\Om$ is an end point of an unpaired edge. The
verifier checks that $m$ vanishes at all these end points, so that $m=0$ on $\partial\Om$, and that every triangle with such a vertex is passive or has $w_T=0$; by
conformity these include all triangles meeting $\partial\Om$. It computes $w_T$ exactly.

\subsection{Weights, cells and correctors}
\emph{Weights.} Positivity of $q_T$ is checked through dyadic bounds $|w_T|\le w_T^+<\kappa_0$; the coefficients $\hat\kappa_T$ and $\hat w_T$ of $q_T^*$
\eqref{qstar} are exact rationals, and $c$ is replaced by the smaller function $\kappa_0-w_T^+$.

\emph{Centres and cells.} Centres and radii are dyadic. For each centre $a$ the verifier checks that $a$ lies in a closed triangle of the mesh, hence in
$\overline\Om$, and that its distance to every side of every $\Gamma_k$ is at least a radius $R_a$ with $R_a\ge R_C$ for every cell $C$ with centre $a$ (so $R'_a\le R_C\le R_a$ for the bump radius $R'_a$ of Definition~\ref{def:bump}); since
$D(a,R_a)$ is connected, contains $a$ and does not meet $\partial\Om$, it lies in $\Om$. Each cell is a sub-triangle of the red refinement of an active triangle,
given by its refinement path; the paths of the cells of one triangle form a complete prefix-free code (no path is a prefix of another, and
$\sum4^{-\text{length}}=1$), so the cells partition the active triangles. $C\subset\overline{D(a_C,r_C)}$ is checked at the vertices of $C$,
$\theta_C^{2p+4}=(r_C^2/R_C^2)^{p+2}$ is rational, and $\rho$ \eqref{eq:rho} is bounded on each mesh triangle by the sum over all cells whose disks are not certified
disjoint from it; so $\varepsilon$ is an exact rational.

\emph{Scaling.} The scale $s_a$ of Definition~\ref{def:cells} is a dyadic number of the order of the radius of the cells with centre $a$, and each functional is in
addition multiplied by a power of two. Such diagonal changes of the coordinates of $\ell(f)$ leave \eqref{crit} invariant: if $\ell$ is replaced by $\Sigma\ell$ with a
positive diagonal matrix $\Sigma$, then $Q$ becomes $\Sigma^{-1}Q\Sigma^{-1}$ and $\mathcal R$ becomes $\Sigma\mathcal R\Sigma$, and \eqref{crit} for $\alpha$ becomes \eqref{crit} for $\Sigma\alpha$.

\emph{Correctors.} The correctors $\chi_j$, $j=1,\dots,2N$, are continuous piecewise quadratic functions on the $n_{\rm ref}$-fold red refinement
of the mesh, with integer nodal values (after scaling by a power of two) produced by a floating-point Galerkin solve of the weighted problem
that defines the optimal corrector (Lemma~\ref{lem:optcorr}); the verifier checks that they vanish at all nodes on $\partial\Om$. Since that problem is only
real-linear, the real and the imaginary part of each Taylor functional get their own corrector, as in Definition~\ref{def:reallin}. Any such $\chi_j$ is admissible in
\eqref{eq:tildek}.

\subsection{Enclosure of the Gram matrix}
The verifier assembles the Gram matrix $\mathsf G$ of the real-linear representers (Lemma~\ref{lem:reallin}(b)) directly as a real symmetric
$2N\times2N$ ball matrix; no complex matrix is converted. For a triangle $\tau$ of the refined mesh, $\hat\kappa_\tau$, $\hat w_\tau$ and $A_\tau$ denote the values of
the mesh triangle containing $\tau$; thus $q_\tau^*(\xi)=\xi^{\mathsf T}A_\tau^{-1}\xi$ with $A_\tau^{-1}=\hat\kappa_\tau I+W(\hat w_\tau)$ (Lemma~\ref{lem:dualform}). All
integrals below are evaluated exactly or enclosed by the following formulas; $\lambda_0,\lambda_1,\lambda_2$ are the barycentric coordinates of a triangle $\tau$ with
area $|\tau|$.
\begin{lemma}[exact integration formulas]\label{lem:intformulas}
\begin{enumerate}
\item[(a)] For integers $e_0,e_1,e_2\ge0$,
\[\int_\tau\lambda_0^{e_0}\lambda_1^{e_1}\lambda_2^{e_2}\,dA=\frac{2|\tau|\,e_0!\,e_1!\,e_2!}{(e_0+e_1+e_2+2)!}.\]
\item[(b)] If $g,g'$ are affine on $\tau$ with vertex values $g_k,g'_k$, and $s:=\sum_kg_k$, $s':=\sum_kg'_k$, then $\int_\tau gg'\,dA=\frac{|\tau|}{12}\big(\sum_kg_kg'_k+ss'\big)$.
\item[(c)] Let $a\in\C$, $U:=z-a$, and $U_k$ the vertices of $\tau$ minus $a$, and let $h_m$ be the complete homogeneous symmetric polynomial of degree $m$. Then
\begin{gather*}\int_\tau\lambda_kU^m\,dA=\frac{2|\tau|\,h_m(U_0,U_1,U_2,U_k)}{(m+1)(m+2)(m+3)},\\
\int_\tau\lambda_k\lambda_lU^m\,dA=\frac{2|\tau|(1+\mathbf 1_{k=l})\,h_m(U_0,U_1,U_2,U_k,U_l)}{(m+1)(m+2)(m+3)(m+4)}.\end{gather*}
\item[(d)] (Green's formula.) If $G\subset\C$ is a bounded region with piecewise smooth, positively oriented boundary and $H$ is $C^1$ on $\overline G$, then
$\int_G\dzb H\,dA=\frac1{2i}\oint_{\partial G}H\,dz$.
\end{enumerate}
\end{lemma}
\begin{proof}
(a) The affine map from the reference triangle $\{x,y\ge0,\ x+y\le1\}$ onto $\tau$ has Jacobian $2|\tau|$, and on the reference triangle the Dirichlet integral
$\int x^{e_1}y^{e_2}(1-x-y)^{e_0}\,dx\,dy=e_0!e_1!e_2!/(e_0+e_1+e_2+2)!$ follows from two Beta integrals.
(b) $g=\sum_kg_k\lambda_k$, and by (a) $\int_\tau\lambda_k\lambda_l=|\tau|(1+\mathbf 1_{k=l})/12$.
(c) Insert $U=\sum_l\lambda_lU_l$ and expand $U^m=\sum_{|e|=m}\frac{m!}{e!}\prod_l(\lambda_lU_l)^{e_l}$. By (a),
$\int_\tau\lambda_kU^m=\frac{2|\tau|\,m!}{(m+3)!}\sum_{|e|=m}(e_k+1)U^e$, and $\sum_{|e|=m}(e_k+1)U^e=h_m(U_0,U_1,U_2,U_k)$, since $e_k+1$ counts the ways to split
$e_k$ between the two arguments $U_k$. Likewise $\int_\tau\lambda_k\lambda_lU^m=\frac{2|\tau|\,m!}{(m+4)!}\sum_{|e|=m}\frac{(e+\mathbf 1_k+\mathbf 1_l)!}{e!}U^e$, where the
factor is $(e_k+1)(e_l+1)$ if $k\ne l$ and $(e_k+2)(e_k+1)$ if $k=l$; these are the numbers of splittings counted by $h_m(U_0,U_1,U_2,U_k,U_l)$, multiplied by
$1$ resp.\ $2$.
(d) With $dz=dx+i\,dy$ and $\dzb=\frac12(\partial_x+i\partial_y)$ this is Green's theorem $\oint(g_1\,dx+g_2\,dy)=\int(\partial_xg_2-\partial_yg_1)\,dA$ for $g_1=H$, $g_2=iH$.
\end{proof}
Each entry of $\mathsf G$ splits by bilinearity into representer--representer, corrector--corrector and representer--corrector parts, and terms with the constants
$\gamma_j$.

\emph{Representer--representer.} Let $R_1^b$, $R_2^b$ be the matrices $R_1$, $R_2$ of Lemma~\ref{lem:grammatrix} with $k_\mu$ replaced by $b_\mu$.
Since $\hat b_{2\mu-1}=b_\mu$ and $\hat b_{2\mu}=ib_\mu$, the entries $\sum_T\Re\int_T(\hat\kappa_T\overline{\hat b_j}\hat b_{j'}+\hat w_T\hat b_j\hat b_{j'})$ for
$(j,j')=(2\mu-1,2\nu-1)$, $(2\mu-1,2\nu)$, $(2\mu,2\nu-1)$, $(2\mu,2\nu)$ are $\Re(R_1^b+R_2^b)_{\mu\nu}$, $-\Im(R_1^b+R_2^b)_{\mu\nu}$, $\Im(R_1^b-R_2^b)_{\mu\nu}$
and $\Re(R_1^b-R_2^b)_{\mu\nu}$, the block of Lemma~\ref{lem:grammatrix}.

\emph{Corrector--corrector:} $\sum_\tau\int_\tau(\dz\chi_i)^{\mathsf T}A_\tau^{-1}\dz\chi_j$ over the fine triangles $\tau$. Since $\dz\chi_i$ is affine on $\tau$,
Lemma~\ref{lem:intformulas}(b) makes this part an exact quadratic form in the integer nodal data. The weights are replaced by integers $\check\kappa_\tau$,
$\check w_\tau$ (after a fixed dyadic scaling) with $\check\kappa_\tau-\hat\kappa_\tau\ge|\check w_\tau-\hat w_\tau|$; since $W$ is real-linear with $W(v)^2=|v|^2I$,
this gives $\check\kappa_\tau I+W(\check w_\tau)\succeq A_\tau^{-1}$. As every element contribution is positive semidefinite and monotone in the weight, the
resulting integer matrix is a Loewner upper bound of this part, and hence the whole form is enlarged.

\emph{Representer--corrector:} integrals $\int_\tau\lambda_k\psi_a(z-a)^n\,dA$, with $\psi_a(z):=1-|z-a|^2/R_a'^2$. If $\tau\subset\overline{D(a,R'_a)}$
(checked at the vertices), they are exact rationals by Lemma~\ref{lem:intformulas}(c), with $\bar U=\sum_l\lambda_l\bar U_l$ in $|U|^2=U\bar U$. Triangles that are
certified disjoint from the disk are skipped, and the triangles cut by the circle are treated by the integrator described next.

\emph{Moments over curved regions.} Let $D'_a:=D(a,R'_a)$. The representer--representer blocks need the moments $\int_G(z-z_0)^m\overline{(z-z_0)}^{\,n}\,dA$ over $G=D'_a\cap D'_{a'}$ and over
$G=T\cap D'_a\cap D'_{a'}$ for passive $T$ with $w_T\ne0$, and the representer--corrector part needs them for the cut triangles. By
Lemma~\ref{lem:intformulas}(d) with $H:=(z-z_0)^m\overline{(z-z_0)}^{\,n+1}/(n+1)$ they become integrals over $\partial G$, which consists of edge segments clipped to
the disks (end points from a quadratic equation) and of circular arcs clipped to the half-planes of the edges and to the other disks (end points from $\arccos$). On
a segment the integrand is a polynomial in the parameter; on an arc $z=z_*+r e^{i\phi}$ it is a trigonometric polynomial $\sum_kc_ke^{ik\phi}$, integrated exactly by
$\int_{\phi_0}^{\phi_1}e^{ik\phi}d\phi=(e^{ik\phi_1}-e^{ik\phi_0})/(ik)$ for $k\ne0$ and $=\phi_1-\phi_0$ for $k=0$. All of this is done in ball arithmetic \cite{Arb}; there is no quadrature and no
truncation error, and a configuration that cannot be decided in ball arithmetic (a tangency, coinciding arc end points) aborts the run. The cell moments needed
for $Q$ are computed in the same way.

\emph{Constants.} The constants $\gamma_j$ are eliminated by the Schur complement with respect to the $2\times2$ matrix $\int_\Om A^{-1}\,dA$, where $A^{-1}$ is the
function equal to $A_T^{-1}$ on $T$; this is exact for the corresponding (optimal, admissible) choice of $\gamma_j$.

Hence the ball matrix contains a real symmetric matrix $\mathsf G^+\succeq\mathsf G$. As built-in tests, the enclosed integrals of each
representer over all fine triangles are summed and compared with the closed forms of its zeroth and first moments over the disk, and
$\sum_\tau\int_\tau\dz\chi_j=0$ is checked exactly. The bump $\psi$ makes the representers continuous. (The certificates with the scalar weight
of Remark~\ref{rem:scalar} use complex-linear representers and a Hermitian $N\times N$ Gram matrix $R$; there the real form of Lemma~\ref{lem:grammatrix} with
$R_1=R$, $R_2=0$ is used.)

\subsection{The matrix inequality}
The matrix of $Q$ is assembled as a real symmetric $2N\times2N$ ball matrix from its complex symmetric part $\sum_{\text{cells}}w_C\int_C\phi_\mu\phi_\nu$ and its
Hermitian part $\eta\sum_{\text{cells}}w_C^+\int_C\bar\phi_\mu\phi_\nu$; using the dyadic bound $w_C^+\ge|w_C|$ only enlarges $Q$ and keeps Theorem~\ref{thm:rel} valid.
A dyadic factor $F$ with $2N$ rows is proposed in floating point, and a rigorous bound $r_F\ge0$ for the spectral norm of the ball matrix $\mathsf G^+-FF^{\mathsf T}$
is computed. Then
\begin{equation}\label{eq:Rup}
\mathcal R_{\rm up}:=FF^{\mathsf T}+r_FI\succeq\mathsf G^+\succeq\mathsf G ,
\end{equation}
and \eqref{crit} is checked through the following two lemmas, without inverting the nearly singular matrix $\mathsf G^+$.
\begin{lemma}[low-rank-plus-identity criterion]\label{lem:lowrank}
Let $F$ be a real $2N\times k$ matrix, $r_F\ge0$, $\mathcal R_{\rm up}$ as in \eqref{eq:Rup}, $Q$ a real quadratic form on $\R^{2N}$ with $\|Q\|\le\mathsf q$, and $\Lambda$ an
upper bound for the largest eigenvalue of the real symmetric matrix $M_F$ of the form $v\mapsto Q(Fv)$ on $\R^k$. Then
\begin{equation}\label{eq:BF}
Q(\mathcal R_{\rm up}\alpha)\le B_F\,\alpha^{\mathsf T}\mathcal R_{\rm up}\alpha\quad(\alpha\in\R^{2N}),\qquad B_F:=\max(\Lambda,r_F\mathsf q)+\sqrt{r_F}\,\mathsf q\,\|F\| .
\end{equation}
Consequently \eqref{crit} holds whenever $B_F\le1-\varepsilon$.
\end{lemma}
\begin{proof}
Given $\alpha$, put $v_1:=F^{\mathsf T}\alpha$ and $v_2:=\sqrt{r_F}\,\alpha$. Then $\mathcal R_{\rm up}\alpha=Fv_1+\sqrt{r_F}\,v_2$ and
$|v_1|^2+|v_2|^2=\alpha^{\mathsf T}\mathcal R_{\rm up}\alpha$. (By Lemma~\ref{lem:step1}, every $\ell(f)$ with $\|f\|_q\le1$ is such a vector with
$|v_1|^2+|v_2|^2\le1$.) Now $Q(Fv_1)\le\Lambda|v_1|^2$, $|Q(Fv_1,v_2)|\le\mathsf q\|F\||v_1||v_2|$ and $Q(v_2)\le\mathsf q|v_2|^2$, so, with $2|v_1||v_2|\le|v_1|^2+|v_2|^2$,
\[Q(Fv_1+\sqrt{r_F}\,v_2)\le\Lambda|v_1|^2+r_F\mathsf q|v_2|^2+\sqrt{r_F}\,\mathsf q\|F\|\big(|v_1|^2+|v_2|^2\big)\le B_F\big(|v_1|^2+|v_2|^2\big).\qedhere\]
\end{proof}
\begin{lemma}[eigenvalue bound by congruence]\label{lem:eigbound}
Let $M$ be a real symmetric $k\times k$ matrix and $V$ a real $k\times k$ matrix such that $Z:=V^{\mathsf T}V$ is strictly diagonally dominant, with margins
$d_i:=Z_{ii}-\sum_{j\ne i}|Z_{ij}|>0$. Let $Y:=V^{\mathsf T}MV$. If $\Lambda\ge0$ and $\Lambda\ge\max_i(Y_{ii}+\sum_{j\ne i}|Y_{ij}|)/d_i$, then $\Lambda I-M$ is positive
semidefinite.
\end{lemma}
\begin{proof}
The symmetric matrix $\Lambda Z-Y$ has diagonal entries $\Lambda Z_{ii}-Y_{ii}$ and off-diagonal row sums $\sum_{j\ne i}|\Lambda Z_{ij}-Y_{ij}|\le\Lambda\sum_{j\ne i}|Z_{ij}|+\sum_{j\ne i}|Y_{ij}|$
(here $\Lambda\ge0$ is used), so it is diagonally dominant with nonnegative diagonal: $\Lambda Z_{ii}-Y_{ii}-\Lambda\sum_{j\ne i}|Z_{ij}|-\sum_{j\ne i}|Y_{ij}|=\Lambda d_i-(Y_{ii}+\sum_{j\ne i}|Y_{ij}|)\ge0$.
By Gershgorin's theorem $\Lambda Z-Y=V^{\mathsf T}(\Lambda I-M)V$ is positive semidefinite. $Z$ is positive definite by the same theorem, so $V$ is invertible, and
$\Lambda I-M$ is positive semidefinite.
\end{proof}
The verifier applies Lemma~\ref{lem:eigbound} to $M=M_F$ with an approximate eigenvector matrix $V$; $Y$ and $Z$ are evaluated in ball arithmetic. Restricting to
$\Lambda\ge0$ is no loss, since only $\max(\Lambda,r_F\mathsf q)$ enters \eqref{eq:BF}.

\subsection{Upper bounds for $\beta$}
\begin{lemma}[enclosed Rayleigh quotients]\label{lem:ub}
Let $g\in\A$, and let $I_1:=\int_\Om g\,dA$, $I_2:=\int_\Om g^2\,dA$ and $I_3:=\int_\Om|g|^2\,dA$ satisfy $I_3>|I_1|^2/|\Om|$. Then $f:=g-I_1/|\Om|\in\Az\setminus\{0\}$ and
\[\kappa(\Om)\ge\frac{|I_2-I_1^2/|\Om||}{I_3-|I_1|^2/|\Om|},\qquad\beta(\Om)^2\le\frac12\Big(1-\frac{|I_2-I_1^2/|\Om||}{I_3-|I_1|^2/|\Om|}\Big).\]
\end{lemma}
\begin{proof}
By Lemma~\ref{lem:mean} with $I=I_1$, $f\in\Az$, $\int f^2=I_2-I_1^2/|\Om|$ and $\|f\|^2=I_3-|I_1|^2/|\Om|>0$. The first expression in \eqref{kappa}
gives the bound for $\kappa$, and Corollary~\ref{cor:bounds}(b) the bound for $\beta$.
\end{proof}
The trial functions $g$ are polynomials plus corner functions $g_\delta$ as in \eqref{eq:gdelta}, i.e.\ $((z-c)/(z-c'))^{-1+\delta}$, with coefficients taken from a
floating-point Rayleigh--Ritz computation and treated as exact rationals. The integrals $I_1$, $I_2$, $I_3$ are reduced to boundary integrals and enclosed with the
rigorous integrator of \cite{Arb}, and Lemma~\ref{lem:ub} is applied.

\subsection{Formal verification in Lean}
Parts of the argument and of the computation are checked in Lean~4 with Mathlib \cite{Mathlib}.

\emph{The theory.} Theorem~\ref{thm:rel} with the form $q$ is formalised for the mean-zero Bergman space of an arbitrary bounded open set $\Om\subset\C$: the
multiplier identity (Lemma~\ref{lem:id}), for multipliers that are Lipschitz, vanish outside $\Om$ and are affine on each triangle, is derived from an
integration-by-parts lemma $\int_\Om F\,\partial_{\bar z}u=0$ for holomorphic $F\in L^1(\Om)$ and Lipschitz $u$ vanishing outside $\Om$, proved by truncation. The
Taylor lemma (Lemma~\ref{lem:taylor}), the representer identity \eqref{rep} with disk representers and Lipschitz correctors, the dual form \eqref{dual} and the remaining
steps of the proof of Theorem~\ref{thm:rel} are proved as well; all of these depend only on Lean's three standard axioms. Theorem~\ref{thm:beta} is not formalised.

\emph{The certificates.} The following table records, for each certificate, which of three parts is checked in Lean: the facts about the mesh, the last stage
(\eqref{crit} for the exported enclosures), and the step from the rational $\kappa_0$ to the stated decimal bound for $\beta$.
\begin{center}\small
\begin{tabular}{llll}
\toprule
certificate ($\kappa_0$) & mesh facts & last stage & $\kappa_0\mapsto\beta$\\
\midrule
square, $64$ and $256$ sectors ($0.63676$, $0.63664$) & Lean & Lean & Lean\\
rectangle $0.7002$, L-shape $0.8152$, $\Om_H$ $0.9522$ & verifier & Lean & Lean\\
other rows of Tables~\ref{tab:square},~\ref{tab:other}; scalar weight & verifier & verifier & Lean\\
square, $512$ sectors ($0.63663$); rectangle $0.70012^\ast$ & verifier & verifier & script\\
rectangle family (Table~\ref{tab:rect}) & verifier & verifier & script\\
\bottomrule
\end{tabular}
\end{center}
Here $\Om_H$ is the square with a hole of Section~\ref{sec:num}, and ``script'' means exact rational arithmetic by the scripts in the repository \cite{Repo}.

\emph{Last stage.} For the five certificates in the first two rows the exported data of the last stage are checked by a verified checker: for these,
$\mathcal R_{\rm up}=FF^{\mathsf T}+r_FI$ is positive definite, dominates every Gram form in the exported enclosure, and satisfies \eqref{crit} for every matrix of $Q$ in
the exported enclosure. The soundness of the checker is kernel-checked; its evaluation on the data uses compiled code. The best square certificate ($512$ sectors,
$\kappa_0=0.63663$, which gives the headline bound $\beta\ge0.426245$), the best rectangle certificate ($\kappa_0=0.70012$) and the certificates of
Table~\ref{tab:rect} were not exported and are \emph{not} among the Lean-checked ones; they rest on the verifier alone.

\emph{Mesh facts.} For the $64$- and $256$-sector certificates all facts about the mesh are proved from the exact data: the triangles cover the square and meet only
in null sets, the multiplier is continuous, affine on each triangle and hence Lipschitz, and the cells and disks have the stated properties ($4.7\cdot10^6$ and
$7.6\cdot10^7$ triangle pairs checked). The resulting Lean theorems state $\kappa((-1,1)^2)\le15919/25000$ and $\kappa((-1,1)^2)\le3979/6250$ under a single remaining
hypothesis: the existence of the correctors (Lipschitz, vanishing outside the square) for which the Gram form and the form $Q$ lie in the enclosures computed with
Arb. The corrector coefficients ($2.8\cdot10^7$ and $4.5\cdot10^8$ integers) are not exported. For the other three Lean-checked certificates only the last stage is
checked; the mesh facts rest on the verifier.

\emph{From $\kappa_0$ to $\beta$.} The step from each rational $\kappa_0$ to the stated decimal bound for $\beta$ is kernel-checked for $37$ certified statements,
namely all rows of Tables~\ref{tab:square} and~\ref{tab:other} except the $512$-sector square and the starred rectangle row, together with the scalar-weight
certificates; for these two rows and for the rectangle family of Table~\ref{tab:rect} this step is checked in exact rational arithmetic by the scripts in the
repository \cite{Repo}, not in Lean.

The formal verification is therefore partial: in every case the hypothesis connecting the unexported correctors to the Arb enclosures remains,
Theorem~\ref{thm:beta} and the computation of the enclosures are not formalised, and the headline certificate for the square is not formally checked.

\subsection{Trust base}
All matrix-weight certificates of Section~\ref{sec:num} were produced by revision~4 of the verifier script
\texttt{lbb\_polygon\_certify.py}.\footnote{SHA-256 of the script: \texttt{f8846e0ed8f6f09525aa41eb40c83600b4b86d452fab79dd7cebcd0fa459d96b}.}
Not verified by machine are the statements proved in Sections~\ref{sec:fried}--\ref{sec:conv}, the closed formulas implemented, and the software stack
(integer and rational arithmetic of the Python interpreter, FLINT and Arb). The correctness of the Gram enclosure is additionally tested against an independent adaptive
quadrature on samples of entries, and the verifier is checked to fail for $\kappa_0$ slightly below the certified values and for deliberately corrupted correctors.

\section{Numerical results}\label{sec:num}
\emph{Test domains.} Besides the square $(-1,1)^2$ and the rectangle $(0,2)\times(0,1)$ we use the rectangles $\Om_r$, $0<r\le1$, the L-shaped domain $\Om_L$ and the
square with a hole $\Om_H$,
\begin{equation}\label{eq:domains}
\Om_r:=(0,1)\times(0,r),\qquad \Om_L:=(-1,1)^2\setminus[0,1]^2,\qquad \Om_H:=(-1,1)^2\setminus[-\tfrac12,\tfrac12]^2 .
\end{equation}
The constants $\kappa$ and $\beta$ are invariant under translations, rotations and dilations, so the square and the rectangle $(0,2)\times(0,1)$ have the constants of
$\Om_1$ and $\Om_{1/2}$. For all these domains $\kess=2/\pi=0.6366198\ldots$, the maximum in the definition of $\kess$ being attained at the right angles.

\begin{figure}[t]
\centering\includegraphics[width=.95\textwidth]{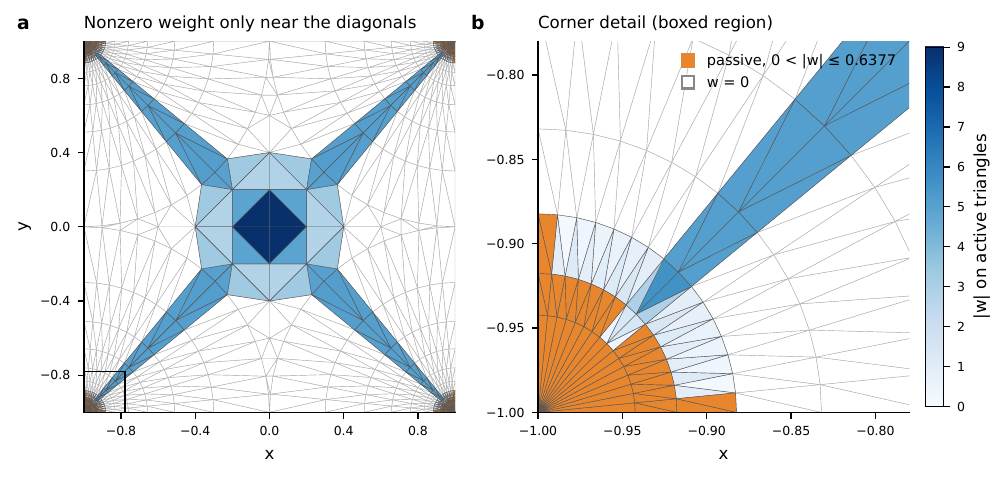}
\caption{An illustrative multiplier for the square ($16$ sectors, two passive rings, $1184$ triangles), not one of the designs used for the
certified results (which have three or four passive rings, Section~\ref{sec:num}): modulus of the weight $w=1-\dzb m$. Of the triangles, $736$ have $w=0$, $192$ are passive with
$0<|w|\le0.637644$, and $256$ are active.}\label{fig:mesh}
\end{figure}

\emph{Certificate data.} The multipliers have the following structure (Remark~\ref{rem:rings}; Figure~\ref{fig:mesh} shows a simpler illustrative design).
\begin{itemize}
\item A fan of $n$ equal sectors per right angle, and of $3n$ sectors at a re-entrant corner, so that every sector has opening $\pi/(2n)\le\pi/16$ and Lemma~\ref{lem:fan} applies;
the fans occupy the three or four innermost rings of each corner mesh, which are passive.
\item The reflection multiplier along the sides, and a wedge of active cells along each bisector.
\end{itemize}
The constant $\kappa_{\mathcal P}$ of the mesh is a floor for the $\kappa_0$ it can certify; it equals $0.637644$, $0.636684$, $0.636636$, $0.636624$, $0.636621$ for $n=16,64,128,256,512$.
The Taylor degree is $p=12$. The correctors are continuous piecewise quadratics on one or two uniform refinements of the mesh; ``corrector nodes'' is the number of interior nodes of
this space. $N$ is the number of complex Taylor functionals ($2N$ real ones). Lower bounds for $\beta$ are rounded down, upper bounds up. The \emph{margin} of a certificate is $1$
minus the verified upper bound of the quotient of the two sides of \eqref{crit}.

All enclosures reported below rest on the same two arguments, which we state once.
\begin{proposition}[verified enclosures]\label{prop:enclosures}
For the square, the rectangle $(0,2)\times(0,1)$, $\Om_L$ and $\Om_H$, the constants $\kappa$ and $\beta$ lie in the intervals of Table~\ref{tab:summary}; for $\Om_r$,
$r\in\{\frac14,\frac38,\frac12,\frac9{16},\frac58,\frac34,\frac78,1\}$, they lie in the intervals of Table~\ref{tab:rect}.
\end{proposition}
\begin{proof}
\emph{Upper bounds for $\kappa$, lower bounds for $\beta$.} Each entry $\kappa\le\kappa_0$ is a certificate of Theorem~\ref{thm:rel}, checked by the verifier of
Section~\ref{sec:impl}; the certificates are listed in Tables~\ref{tab:square}, \ref{tab:other} and~\ref{tab:rect-runs}. By Theorem~\ref{thm:beta},
$\beta=\big(\frac{1-\kappa}2\big)^{1/2}\ge\big(\frac{1-\kappa_0}2\big)^{1/2}$, and this value is rounded down.

\emph{Lower bounds for $\kappa$, upper bounds for $\beta$.} For the square and for $\Om_r$ with $r\ge\frac34$ the lower bound is $\kappa\ge\kess=2/\pi$ (Lemma~\ref{lem:kess}).
In all other cases it is the enclosed Rayleigh quotient of one explicit trial function (Lemma~\ref{lem:ub}), which is a lower bound
for $\kappa$ by \eqref{kappa}. By Theorem~\ref{thm:beta}, $\beta\le\big(\frac{1-\kappa_{\rm low}}2\big)^{1/2}$ for every lower bound $\kappa_{\rm low}$ of $\kappa$, and this value is rounded up.
\end{proof}

\subsection{Square}
\begin{figure}[t]
\centering\includegraphics[width=.6\textwidth]{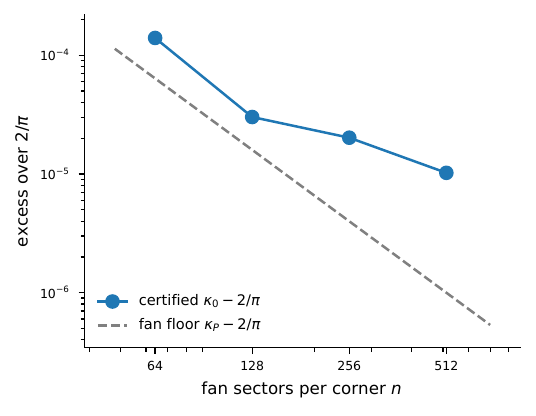}
\caption{Square: distance of the certified upper bound $\kappa_0$ to $\kess=2/\pi$ (dots) and of the fan floor $\kappa_{\mathcal P}$ (dashed) against the number $n$ of fan sectors per corner.
The point $n=64$ uses one refinement of the corrector mesh, the others two.}\label{fig:conv}
\end{figure}

Table~\ref{tab:square} lists certified values for $\Om=(-1,1)^2$. With $512$ sectors the verifier certifies $\kappa\le0.63663$.
By Proposition~\ref{prop:enclosures},
\[0.426245\le\beta\big((-1,1)^2\big)\le\sqrt{\tfrac12-\tfrac1\pi}=0.426251\ldots,\qquad \tfrac{2}{\pi}\le\kappa\le0.63663 ,\]
compared with $\beta\ge\sin(\pi/8)=0.38268\ldots$ from \cite{HorganPayne,CostabelDauge}. The certified values approach $\kess$ as the number of sectors grows (Figure~\ref{fig:conv}):
the distance to $2/\pi$ is $1.4\cdot10^{-4}$, $3.0\cdot10^{-5}$, $2.0\cdot10^{-5}$, $1.0\cdot10^{-5}$ for $n=64,128,256,512$. For $n\ge256$ the result is limited not by the fan floor
$\kappa_{\mathcal P}$ but by the corrector space; a third refinement of the corrector mesh for $n=512$ ($3.1\cdot10^6$ nodes) exceeds the index range of the sparse direct solver
used to propose the correctors. The data are consistent with $\kappa=\kess$, i.e.\ $\beta^2=\frac12-\frac1\pi$, for the square. The enclosure leaves an interval of length less than
$6.24\cdot10^{-6}$ for $\beta$ and does not decide the question.
\begin{table}[t]
\caption{Square $(-1,1)^2$: certified upper bounds $\kappa_0$ for $\kappa$ and resulting lower bounds for $\beta$.}\label{tab:square}
\begin{tabular}{rrrrlll}
\toprule
sectors $n$ & $\kappa_{\mathcal P}$ & $N$ & corrector nodes & $\kappa_0$ & $\beta\ge$ & margin\\
\midrule
16 & $0.637644$ & 637 & 12\,353 & $0.638$ & $0.425440$ & $1.9\cdot10^{-2}$\\
16 & $0.637644$ & 637 & 49\,793 & $0.63775$ & $0.425587$ & $2.8\cdot10^{-2}$\\
64 & $0.636684$ & 286 & 24\,545 & $0.63676$ & $0.426168$ & $4.3\cdot10^{-2}$\\
128 & $0.636636$ & 286 & 49\,121 & $0.63671$ & $0.426198$ & $2.5\cdot10^{-2}$\\
128 & $0.636636$ & 286 & 196\,673 & $0.63665$ & $0.426233$ & $3.2\cdot10^{-2}$\\
256 & $0.636624$ & 286 & 393\,281 & $0.63664$ & $0.426239$ & $7.3\cdot10^{-2}$\\
512 & $0.636621$ & 286 & 786\,497 & $0.63663$ & $0.426245$ & $2.2\cdot10^{-2}$\\
\bottomrule
\end{tabular}
\end{table}

\subsection{Rectangle, L-shaped domain, and a doubly connected domain}
Table~\ref{tab:other} lists certified values for $(0,2)\times(0,1)$, for $\Om_L$ and for $\Om_H$ (see \eqref{eq:domains}).
The domain $\Om_H$ is not star-shaped, so no estimate of Horgan--Payne type applies, whereas Theorem~\ref{thm:beta} and the method apply without change. For $\Om_H$ the verifier checks
the orientation and nesting of the boundary components, the tiling, and the distances to all sides exactly.
The upper bounds for $\beta$ come from enclosed Rayleigh quotients of explicit trial functions (polynomials and corner functions; for $\Om_H$ also negative powers of $z$, and corner functions
at the eight corners whose branch cuts lie outside $\Om_H$): $\kappa((0,2)\times(0,1))\ge0.700056156290$, $\kappa(\Om_L)\ge0.815074344196$ and $\kappa(\Om_H)\ge0.952086244272$, the latter with $176$ coefficients
($0.951895821127$ with $36$). In all three cases $\kappa>\kess$, i.e.\ $\beta^2$ lies strictly below the part of the spectrum generated by the corners (Lemma~\ref{lem:cosserat} below).
For the rectangle the enclosure contains the value $0.387262$ computed in \cite{BCDG}. The margins of the best runs are small (between $1.5\cdot10^{-5}$ and $3.6\cdot10^{-5}$); $0.7001$ is not
certified for the rectangle with two refinements of the corrector mesh, and $0.8152$ and $0.9522$ are not certified on the coarser ones.
\begin{table}[t]
\caption{Certified upper bounds $\kappa_0$ for $\kappa$ on two corrector meshes, and resulting lower bounds for $\beta$ ($n=16$; $n=8$ for $\Om_H$). $^\ast$: the same three-ring design ($r_{\min}=0.002$) on the congruent copy $(-1,1)\times(-\frac12,\frac12)$, with three refinements of the
corrector mesh (the run listed for $r=\frac12$ in Table~\ref{tab:rect-runs}).}\label{tab:other}
\begin{tabular}{lrrlll}
\toprule
domain & $N$ & corrector nodes & $\kappa_0$ & $\beta\ge$ & margin\\
\midrule
rectangle $(0,2)\times(0,1)$ & 1027 & 20\,937 & $0.7004$ & $0.387040$ & $1.6\cdot10^{-4}$\\
 & 1027 & 84\,433 & $0.7002$ & $0.387169$ & $8.9\cdot10^{-5}$\\
 & 1027 & 339\,105 & $0.70012^\ast$ & $0.387220$ & $1.8\cdot10^{-5}$\\
L-shape $(-1,1)^2\setminus[0,1]^2$ & 1508 & 36\,669 & $0.8153$ & $0.303891$ & $5.1\cdot10^{-5}$\\
 & 1508 & 147\,609 & $0.8152$ & $0.303973$ & $3.6\cdot10^{-5}$\\
square with hole $\Om_H$ & 1300 & 23\,400 & $0.9524$ & $0.154272$ & $1.0\cdot10^{-4}$\\
 & 1300 & 94\,608 & $0.9522$ & $0.154596$ & $1.5\cdot10^{-5}$\\
\bottomrule
\end{tabular}
\end{table}

\subsection{Rectangles of varying aspect ratio}
\begin{figure}[t]
\centering\includegraphics[width=.95\textwidth]{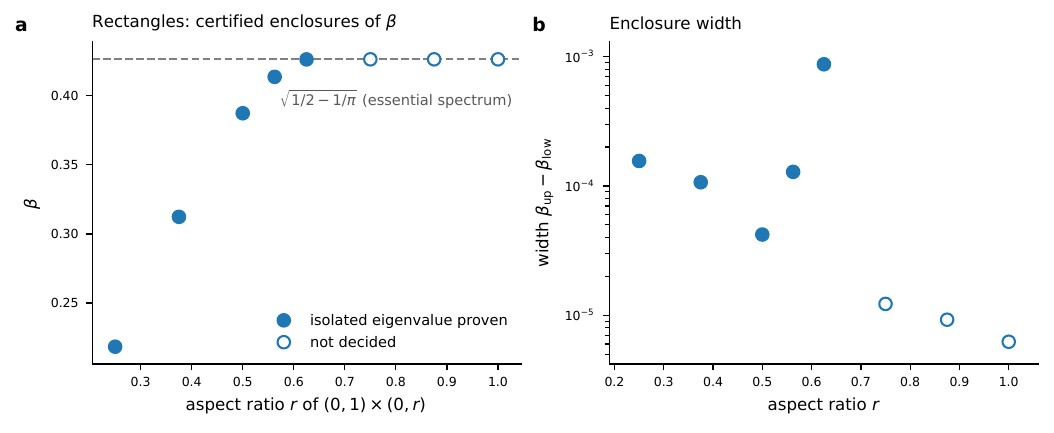}
\caption{Rectangles $\Om_r=(0,1)\times(0,r)$: (a) certified enclosures of $\beta$ (the two bounds coincide at this scale); filled: $\beta^2$ proven to be an isolated eigenvalue; dashed:
$\sqrt{1/2-1/\pi}$. (b) Width $\beta_{\rm up}-\beta_{\rm low}$ of the enclosure, computed from the upper endpoint
$\beta_{\rm up}=\sqrt{(1-\kappa_{\rm low})/2}$ in full precision (rounded up to $12$ digits) and the lower endpoint $\beta_{\rm low}$ rounded
down to six digits as in Table~\ref{tab:rect}.}\label{fig:rect}
\end{figure}

For the rectangles $\Om_r$ of \eqref{eq:domains}, $\kess=2/\pi$ for all $r$, while $\beta(\Om_r)\to0$ as $r\to0$. Hence $\beta(\Om_r)^2$ is an isolated eigenvalue of the Cosserat operator
(Definition~\ref{def:cosserat}) for sufficiently elongated rectangles, whereas for the square the question is open. Table~\ref{tab:rect} and Figure~\ref{fig:rect} give certified enclosures for
eight dyadic aspect ratios. The upper bounds for $\kappa$ are certified with three-ring designs ($16$ sectors; $64$ for $r=\frac58$) for $r\in\{\frac14,\frac38,\frac12,\frac9{16},\frac58\}$
and with fine fans ($512$ sectors) for $r\in\{\frac34,\frac78,1\}$ (Table~\ref{tab:rect-runs}). The fine fans perform poorly when a genuine eigenfunction exists, since their passive
rings cover regions where it has mass. The lower bounds for $\kappa$ are enclosed Rayleigh quotients of trial functions (polynomials and corner functions).

\begin{definition}[Cosserat operator]\label{def:cosserat}
For $q\in L^2_0(\Om)$ let $u_q\in H^1_0(\Om)^2$ be the Riesz representer of Lemma~\ref{lem:riesz}, i.e.\ $(\nabla u_q,\nabla v)=(q,\Div v)$ for all $v\in H^1_0(\Om)^2$; thus $\Delta u_q=\nabla q$
with the Dirichlet Laplacian. The \emph{Cosserat operator} \cite{CCDL} is
\begin{equation}\label{eq:cosserat}
\mathcal C q:=\Div u_q\qquad(q\in L^2_0(\Om)),
\end{equation}
in operator notation $\mathcal C=\Div\Delta^{-1}\nabla$.
\end{definition}
Part~(a) of the next lemma is \cite[(1.4)--(1.5)]{CCDL} 
(see also \cite[Lemma~A.1]{HSV2012}). 
Part~(b) rests on the description of the essential spectrum on polygons in \cite[Thm.~3.3, (3.10)]{CCDL}, 
whose bottom is $\inf\sigma_{\rm ess}(\mathcal C)=\frac12-\frac{\kess}2$ (in \cite[(3.10)]{CCDL} the interval is to be read with $|\sin\omega_c|$ at re-entrant corners).
\begin{lemma}\label{lem:cosserat}
Let $\Om$ be a polygonal domain.
\begin{enumerate}
\item[(a)] $\mathcal C$ is a bounded self-adjoint operator on $L^2_0(\Om)$ with $(\mathcal Cq,q)=|u_q|_1^2\le\|q\|^2$, and $\beta(\Om)^2=\min\sigma(\mathcal C)$.
\item[(b)] $\kappa(\Om)>\kess$ if and only if $\beta(\Om)^2<\inf\sigma_{\rm ess}(\mathcal C)$. In this case $\beta(\Om)^2$ is an isolated eigenvalue of finite multiplicity of $\mathcal C$.
\end{enumerate}
\end{lemma}
\begin{proof}
(a) For $v\in H^1_0(\Om)^2$, $\int_\Om\Div v\,dA=0$ (by density of $C_c^\infty$), so $\mathcal C$ maps into $L^2_0(\Om)$. For $q,q'\in L^2_0(\Om)$,
$(\mathcal Cq,q')=(\Div u_q,q')=(\nabla u_{q'},\nabla u_q)$, which is symmetric in $q,q'$. By \eqref{eq:beta-b}, $(\mathcal Cq,q)=(q,\Div u_q)=|u_q|_1^2$, and
$|u_q|_1^2\le\|q\|\,\|\Div u_q\|\le\|q\|\,|u_q|_1$ because $\|\Div v\|\le|v|_1$ by \eqref{eq:beta-a}; hence $|u_q|_1\le\|q\|$ and $\|\mathcal Cq\|\le|u_q|_1\le\|q\|$.
For a bounded self-adjoint operator the minimum of the spectrum is the infimum of $(\mathcal Cq,q)/\|q\|^2$ over $q\ne0$, which equals $\inf_q|u_q|_1^2/\|q\|^2=\beta(\Om)^2$ by
\eqref{eq:beta-b}.

(b) By Theorem~\ref{thm:beta}, $\beta^2=\frac{1-\kappa}2$, so $\kappa>\kess$ if and only if $\beta^2<\frac12-\frac{\kess}2=\inf\sigma_{\rm ess}(\mathcal C)$, the last equality by \cite[Thm.~3.3, (3.10)]{CCDL}. In that case
$\beta^2=\min\sigma(\mathcal C)$ lies in the spectrum but not in the essential spectrum. For a self-adjoint operator every such point is an isolated eigenvalue of finite multiplicity
\cite[Sec.~VII.3]{ReedSimon}.
\end{proof}
For $r\le0.53127$, in particular for $r\in\{\frac14,\frac38,\frac12\}$, the conclusion of the next proposition already follows from the explicit upper bound
$\beta(\Om_r)^2\le1-\tanh(t_r)/t_r$, $t_r:=\pi r/2$, of \cite[Sect.~5.1, eq.~(5.3)]{CCDL}, 
which is below $\frac12-\frac1\pi$ exactly for these $r$. The cases $r=\frac9{16}$ and $r=\frac58$ are new.
\begin{proposition}\label{prop:rect}
For $r\in\{\frac14,\frac38,\frac12,\frac9{16},\frac58\}$, $\kappa(\Om_r)>2/\pi$; hence $\beta(\Om_r)^2<\frac12-\frac1\pi$ is an isolated eigenvalue of finite multiplicity of the
Cosserat operator \eqref{eq:cosserat}.
\end{proposition}
\begin{proof}
For these $r$ the certified lower bounds $\kappa\ge\kappa_{\rm low}$ of Table~\ref{tab:rect} (enclosed Rayleigh quotients, Proposition~\ref{prop:enclosures}) exceed $2/\pi=\kess$. The
tightest case is $r=\frac58$, where $\kappa_{\rm low}=0.636712977439$ exceeds $2/\pi$ by $9.3\cdot10^{-5}$. Lemma~\ref{lem:cosserat}(b) gives the claim, with
$\beta^2=\frac{1-\kappa}2<\frac12-\frac1\pi$ by Theorem~\ref{thm:beta}.
\end{proof}
For $r\in\{\frac34,\frac78,1\}$ no trial function found exceeds $0.6366187<2/\pi$; the
certified upper bounds leave windows of length $1.0$--$2.0\cdot10^{-5}$ above $2/\pi$ for $\kappa$, so whether $\beta^2$ equals the bottom of the essential spectrum is not decided there.
\begin{table}[t]
\caption{Rectangles $\Om_r=(0,1)\times(0,r)$: verified enclosures. For $r\ge\frac34$ the lower bound for $\kappa$ is $\kess$ (Lemma~\ref{lem:kess}). The last column states whether $\beta^2$ is proven to be an
isolated eigenvalue below the essential spectrum (Proposition~\ref{prop:rect}).}\label{tab:rect}
\renewcommand{\arraystretch}{1.35}\begin{tabular}{llllll}
\toprule
$r$ & $\kappa\ge$ & $\kappa\le$ & $\beta\ge$ & $\beta\le$ & isolated\\
\midrule
$\frac{1}{4}$ & $0.904564387222$ & $0.9047$ & $0.218288$ & $0.218445$ & yes\\
$\frac{3}{8}$ & $0.804867674457$ & $0.805$ & $0.312249$ & $0.312356$ & yes\\
$\frac{1}{2}$ & $0.700056156290$ & $0.70012$ & $0.387220$ & $0.387263$ & yes\\
$\frac{9}{16}$ & $0.657887461604$ & $0.6581$ & $0.413461$ & $0.413590$ & yes\\
$\frac{5}{8}$ & $0.636712977439$ & $0.6382$ & $0.425323$ & $0.426197$ & yes\\
$\frac{3}{4}$ & $2/\pi$ & $0.63664$ & $0.426239$ & $0.426252$ & open\\
$\frac{7}{8}$ & $2/\pi$ & $0.636635$ & $0.426242$ & $0.426252$ & open\\
$1$ & $2/\pi$ & $0.63663$ & $0.426245$ & $0.426252$ & open\\
\bottomrule
\end{tabular}
\end{table}
\begin{table}[t]
\caption{Rectangles: the certificates behind the upper bounds $\kappa\le\kappa_0$ of Table~\ref{tab:rect} (matrix weight, $p=12$). Sectors per
right angle, passive rings, refinements $n_{\rm ref}$ of the corrector mesh, number $N$ of complex Taylor functionals, interior corrector
nodes, margin, and the run identifier under which the certificate is stored in \cite{Repo}.}\label{tab:rect-runs}
\centering\footnotesize\setlength{\tabcolsep}{4pt}
\renewcommand{\arraystretch}{1.4}\begin{tabular}{lrrrrrll}
\toprule
$r$ & sectors & rings & $n_{\rm ref}$ & $N$ & nodes & margin & run\\
\midrule
$\frac14$ & 16 & 3 & 2 & 1131 & 97\,073 & $1.0\cdot10^{-4}$ & \texttt{R5:rect\_r1\_4\_t16:nref2:9047/10000}\\
$\frac38$ & 16 & 3 & 2 & 1079 & 88\,537 & $1.0\cdot10^{-4}$ & \texttt{R5:rect\_r3\_8\_t16:nref2:161/200}\\
$\frac12$ & 16 & 3 & 3 & 1027 & 339\,105 & $1.8\cdot10^{-5}$ & \texttt{R5:rectc\_r1\_2\_t16:nref3:17503/25000}\\
$\frac9{16}$ & 16 & 3 & 3 & 1027 & 332\,993 & $5.8\cdot10^{-5}$ & \texttt{R5:rectc\_r9\_16\_t16:nref3:6581/10000}\\
$\frac58$ & 64 & 3 & 2 & 1053 & 315\,377 & $1.1\cdot10^{-3}$ & \texttt{R5:rect\_r5\_8\_t64:nref2:3191/5000}\\
$\frac34$ & 512 & 4 & 2 & 312 & 787\,377 & $3.5\cdot10^{-2}$ & \texttt{R5:rect\_r3\_4\_n512:nref2:3979/6250}\\
$\frac78$ & 512 & 4 & 2 & 312 & 787\,065 & $8.9\cdot10^{-2}$ & \texttt{R5:rect\_r7\_8\_n512:nref2:127327/200000}\\
$1$ & 512 & 4 & 2 & 286 & 786\,497 & $2.2\cdot10^{-2}$ & \texttt{R5:sqf\_n512:nref2:63663/100000}\\
\bottomrule
\end{tabular}
\end{table}

\subsection{Comparison with the scheme of Gallistl}
Gallistl's Rayleigh--Ritz scheme \cite{Gallistl19} discretises the pressure by piecewise polynomials of a fixed degree and yields values $\mu_h=\beta_h^2$ that decrease monotonically under
mesh refinement and are upper bounds for $\mu:=\beta^2$ when the discrete eigenvalue problem is solved exactly. Its a posteriori analysis \cite{Gallistl21} adds the quantity
\[\mu_h^{\rm lb}:=\mu_h-2\sqrt{2\mu_h}\,\eta_{\rm G},\]
where $\eta_{\rm G}$ is a computable error estimator. This quantity is a lower bound for $\mu$ under two hypotheses: $\mu$ is
an eigenvalue separated from the essential spectrum, and the mesh is fine enough that a discrete quantity is at least $\frac12$ \cite[Sect.~4]{Gallistl21}.
Table~\ref{tab:gallistl} compares the values on the finest adaptive meshes reported in \cite[Tables~1 and~3]{Gallistl21} (polynomial degree $4$) with our enclosures.
\begin{table}[t]
\caption{Squared inf-sup constant $\mu=\beta^2$. Left: values of \cite[Tables~1, 3]{Gallistl21} on the finest reported adaptive mesh. $\mu_h$ is a Rayleigh--Ritz
upper bound, computed in floating point. $\mu_h^{\rm lb}$ is a lower bound only under the two hypotheses of \cite{Gallistl21}; for the square the separation hypothesis is not known
to hold, and \cite{Gallistl21} labels those values experimental. Right: our verified enclosures (Tables~\ref{tab:summary} and~\ref{tab:rect}), lower ends rounded down, upper ends rounded up; the width is that of the unrounded enclosure.}\label{tab:gallistl}
\small
\begin{tabular}{lllll}
\toprule
domain & $\mu_h^{\rm lb}$ \cite{Gallistl21} & $\mu_h$ \cite{Gallistl21} & verified enclosure of $\mu$ (this paper) & width\\
\midrule
$(0,2)\times(0,1)$ & $0.148262$ & $0.149977$ & $[0.149940,\ 0.149972]$ & $3.2\cdot10^{-5}$\\
square & $0.174966$ & $0.182943$ & $[0.181685,\ 0.181691]$ & $5.1\cdot10^{-6}$\\
\bottomrule
\end{tabular}
\end{table}
For the $2:1$ rectangle, the reference value $\mu=0.1499718$ of \cite{BCDG} used in \cite{Gallistl21} lies in our enclosure. Our enclosure is about $50$ times narrower than the interval
$[\mu_h^{\rm lb},\mu_h]$ and holds without hypotheses. For the square the two intervals differ in kind. Under the conjecture $\mu=\frac12-\frac1\pi$ the values of \cite{Gallistl21}
converge slowly towards that number, whereas the certificate excludes every $\mu<0.181685$ unconditionally.
For aspect ratio $1.61$, \cite{Gallistl19} computed $\mu_h\le0.18159009<\frac12-\frac1\pi$ in floating point. Since $\mu_h$ is an upper bound, this indicates that $\mu$ is an isolated eigenvalue for a
rectangle slightly less elongated than the golden ratio, which \cite[p.~454]{CCDL} had mentioned as a possible transition point. Proposition~\ref{prop:rect} proves this rigorously for the nearby
ratio $1.6$ ($r=\frac58$), where $\mu\le0.181644<\frac12-\frac1\pi$. The two approaches are complementary. The Rayleigh--Ritz values give upper bounds for $\beta$ without corner
analysis, and the present certificate supplies the unconditional lower bounds.

\subsection{Summary and further observations}
\begin{table}[t]
\caption{Verified enclosures of $\kappa$ and $\beta$ (Proposition~\ref{prop:enclosures}). The lower bound for $\kappa$ on the square is Lemma~\ref{lem:kess}; the others are enclosed Rayleigh quotients.}\label{tab:summary}
\begin{tabular}{lllll}
\toprule
domain & $\kappa\ge$ & $\kappa\le$ & $\beta\ge$ & $\beta\le$\\
\midrule
square $(-1,1)^2$ & $2/\pi$ & $0.63663$ & $0.426245$ & $0.426252$\\
rectangle $(0,2)\times(0,1)$ & $0.700056156290$ & $0.70012$ & $0.387220$ & $0.387263$\\
L-shape $(-1,1)^2\setminus[0,1]^2$ & $0.815074344196$ & $0.8152$ & $0.303973$ & $0.304078$\\
square with hole $\Om_H$ & $0.952086244272$ & $0.9522$ & $0.154596$ & $0.154781$\\
\bottomrule
\end{tabular}
\end{table}
Table~\ref{tab:summary} collects the enclosures; $\beta$ is determined to four significant digits for the square and to three for the other domains.

\emph{Scalar weight.} A first series of certificates used the scalar weight of Remark~\ref{rem:scalar}. They are valid but weaker: $\kappa\le0.63772$ for the square ($n=16$, two passive rings,
$N=1649$, $150\,977$ corrector nodes), $0.7008$ for the rectangle, $0.8154$ for the L-shaped domain and $0.955$ for $\Om_H$. In floating point, on the square with three passive rings the
smallest $\kappa_0$ passing the scalar criterion stays between $8.5\cdot10^{-3}$ and $1.9\cdot10^{-2}$ above $\kappa_{\mathcal P}$, whereas with the form $q$ it is within $6.2\cdot10^{-5}$, in
agreement with Remark~\ref{rem:rings}.

\emph{Other polygons.} In floating point (not certified) the criterion with the form $q$ passes at $\kappa_0=0.505404$ for a regular pentagon ($\kappa_{\mathcal P}=0.505372$), at $0.8282$ for an
equilateral triangle ($\kappa_{\mathcal P}=0.828044$; an enclosed trial function gives $\kappa\ge0.826974445484$), and at $0.965953$ for a right triangle with legs $2$ and $1$. The present
implementation requires dyadic vertices (verifier) and axis-parallel sides (data builder).

\emph{Cost.} The certified runs for the square take between $200$ seconds ($n=64$) and about two hours ($n=512$, $40$\,GB) on shared-memory servers; the runs of Tables~\ref{tab:other} and~\ref{tab:rect} between ten minutes
and $100$ minutes.

\section{Application: guaranteed error bounds for the Stokes equations}\label{sec:app}
\emph{Notation.} In this section $\Om\subset\R^2$ is a polygonal domain as in Section~\ref{sec:fried}; everything is two-dimensional. To keep the letters $f$, $u$, $w$ and $p$ of
the previous sections (holomorphic functions, Lipschitz functions, the weight \eqref{eq:W}, the Taylor degree) free, vector fields of the Stokes problem are set in bold and the
pressure is denoted by $P$. Let $\boldsymbol f\in L^2(\Om)^2$ and let $(\boldsymbol u,P)\in H^1_0(\Om)^2\times L^2_0(\Om)$ solve
\begin{equation}\label{stokes}
(\nabla\boldsymbol u,\nabla\boldsymbol v)-(P,\Div\boldsymbol v)=(\boldsymbol f,\boldsymbol v)\quad\forall\boldsymbol v\in H^1_0(\Om)^2,\qquad (\Div\boldsymbol u,q)=0\quad\forall q\in L^2_0(\Om).
\end{equation}
Here $(\cdot,\cdot)$ is the $L^2$ inner product of scalar, vector or matrix fields.

For any approximation $(\boldsymbol u_h,P_h)\in H^1_0(\Om)^2\times L^2(\Om)$, an upper bound for the velocity error requires control of $\Div\boldsymbol u_h$, and any bound for the pressure
error requires a lower bound for $\beta$; these are the two places where the inf-sup constant enters a posteriori error control. Bounds of this kind belong to the
class of functional, or equilibrated-flux, a posteriori estimates: an auxiliary flux $\sigma\in H(\Div)$ enters as in the hypercircle method of Prager and Synge \cite{PragerSynge},
which Liu, Nakao, You and Oishi extend to explicit error estimates for the Stokes problem with exactly divergence-free discrete velocities \cite{LiuNakaoYouOishi}.
Guaranteed bounds that combine an equilibrated flux with the inf-sup constant are known: \cite[Thm.~5.1]{HSV2012} 
bounds $(|\boldsymbol u-\boldsymbol u_h|_1^2+\beta^2\|P-P_h\|^2)^{1/2}$ with the factor $2/(\sqrt5-1)$, and in the functional a posteriori estimates of Repin
\cite{Repin2002,Repin2008} the estimates for the pressure, and for velocity approximations that are not divergence free, contain the inf-sup constant
\cite{Repin2002,RepinStenberg}. Theorem~\ref{thm:apost} records the separated velocity and pressure bounds in the form used below; we claim no novelty for it and include its short
proof for completeness. What the present paper contributes is a certified value of $\beta_0$; with an explicit $\beta_0\le\beta$ the bounds are fully computable.

\emph{The divergence-free subspace.} Let
\begin{equation}\label{eq:Z}
Z:=\{\boldsymbol v\in H^1_0(\Om)^2:\Div\boldsymbol v=0\},\qquad Z^\perp:=\{\boldsymbol v\in H^1_0(\Om)^2:(\nabla\boldsymbol v,\nabla\boldsymbol z)=0\ \ \forall\boldsymbol z\in Z\},
\end{equation}
and let $\Pi_Z$ be the $(\nabla\cdot,\nabla\cdot)$-orthogonal projection onto $Z$. The following is the standard characterisation of $1/\beta$ as the norm of the right inverse of the
divergence on $Z^\perp$; see \cite[Lemma~A.1]{HSV2012} and \cite[(1.5)]{CCDL}. 
We include the short proof.
\begin{lemma}\label{lem:zperp}
Let $\Om$ be a bounded domain with $\beta(\Om)>0$ (for instance a polygonal domain), and for $q\in L^2_0(\Om)$ let $u_q\in H^1_0(\Om)^2$ be the Riesz representer of
Lemma~\ref{lem:riesz} (which keeps its non-bold notation), $(\nabla u_q,\nabla\boldsymbol v)=(q,\Div\boldsymbol v)$ for all $\boldsymbol v\in H^1_0(\Om)^2$. Then $Z$ is closed, $Z^\perp=\{u_q:q\in L^2_0(\Om)\}$, and
\[|\boldsymbol v|_1\le\beta(\Om)^{-1}\|\Div\boldsymbol v\|\qquad\text{for all }\boldsymbol v\in Z^\perp .\]
\end{lemma}
\begin{proof}
$Z$ is closed because $\|\Div\boldsymbol v\|\le|\boldsymbol v|_1$ by \eqref{eq:beta-a}. Let $U:=\{u_q:q\in L^2_0(\Om)\}$. The map $q\mapsto u_q$ is linear, and
$\beta\|q\|\le|u_q|_1\le\|q\|$ by \eqref{eq:beta-b} and the Cauchy--Schwarz inequality; since $\beta>0$, $U$ is closed. Since $\int_\Om\Div\boldsymbol v\,dA=0$ for
$\boldsymbol v\in H^1_0(\Om)^2$, the function $\Div\boldsymbol v$ lies in $L^2_0(\Om)$; hence $\boldsymbol v\perp U$ if and only if $(q,\Div\boldsymbol v)=0$ for all $q\in L^2_0(\Om)$, if and only
if $\Div\boldsymbol v=0$. Thus $U^\perp=Z$ and $U=U^{\perp\perp}=Z^\perp$. Given $\boldsymbol v\in Z^\perp$, write $\boldsymbol v=u_q$. Then
$|\boldsymbol v|_1^2=(\nabla u_q,\nabla\boldsymbol v)=(q,\Div\boldsymbol v)\le\|q\|\,\|\Div\boldsymbol v\|\le\beta^{-1}|\boldsymbol v|_1\|\Div\boldsymbol v\|$.
\end{proof}

\begin{theorem}\label{thm:apost}
Let $\Om$ be a polygonal domain, $\boldsymbol f\in L^2(\Om)^2$, $(\boldsymbol u,P)$ the solution of \eqref{stokes}, and $(\boldsymbol u_h,P_h)$ any element of $H^1_0(\Om)^2\times L^2(\Om)$.
Let $\beta_0\in(0,\beta(\Om)]$, let $C_{\mathrm P}\ge\lambda_D^{-1/2}$, where $\lambda_D:=\inf_{0\ne v\in H^1_0(\Om)}\|\nabla v\|^2/\|v\|^2$ is the first Dirichlet eigenvalue of $-\Delta$,
and let $\sigma\in H(\Div,\Om)^2$ (row-wise). Define
\begin{equation}\label{eq:etas}
G_h:=\nabla\boldsymbol u_h-P_hI,\qquad \eta_R:=\|\sigma-G_h\|+C_{\mathrm P}\|\boldsymbol f+\Div\sigma\|,\qquad \eta_D:=\|\Div\boldsymbol u_h\|,
\end{equation}
and $\|q\|_{L^2/\R}:=\min_{c\in\R}\|q-c\|$. Then
\begin{equation}\label{apost}
|\boldsymbol u-\boldsymbol u_h|_1\le\Big(\eta_R^2+\frac{\eta_D^2}{\beta_0^2}\Big)^{1/2},\qquad
\|P-P_h\|_{L^2/\R}\le\frac1{\beta_0}\Big(\eta_R+\frac{\eta_D}{\beta_0}\Big).
\end{equation}
\end{theorem}
\begin{proof}
Let $\boldsymbol e:=\boldsymbol u-\boldsymbol u_h$ and $e_P:=P-P_h$.

\emph{Step 1: error equation.} Since $(P_hI,\nabla\boldsymbol v)=(P_h,\operatorname{tr}\nabla\boldsymbol v)=(P_h,\Div\boldsymbol v)$, we have
$(G_h,\nabla\boldsymbol v)=(\nabla\boldsymbol u_h,\nabla\boldsymbol v)-(P_h,\Div\boldsymbol v)$. Define the residual
$\operatorname{Res}(\boldsymbol v):=(\boldsymbol f,\boldsymbol v)-(G_h,\nabla\boldsymbol v)$. With the first equation in \eqref{stokes}, for all $\boldsymbol v\in H^1_0(\Om)^2$,
\[(\nabla\boldsymbol e,\nabla\boldsymbol v)-(e_P,\Div\boldsymbol v)=\big[(\nabla\boldsymbol u,\nabla\boldsymbol v)-(P,\Div\boldsymbol v)\big]-(G_h,\nabla\boldsymbol v)=\operatorname{Res}(\boldsymbol v).\]

\emph{Step 2: residual bound.} Each row $\sigma_i$ of $\sigma$ is in $H(\Div,\Om)$, so for $\boldsymbol v\in H^1_0(\Om)^2$
\[(\sigma,\nabla\boldsymbol v)=\sum_i(\sigma_i,\nabla v_i)=-\sum_i(\Div\sigma_i,v_i)=-(\Div\sigma,\boldsymbol v).\]
Adding $0=(\Div\sigma,\boldsymbol v)+(\sigma,\nabla\boldsymbol v)$ gives $\operatorname{Res}(\boldsymbol v)=(\boldsymbol f+\Div\sigma,\boldsymbol v)+(\sigma-G_h,\nabla\boldsymbol v)$. By the
Cauchy--Schwarz inequality and $\|\boldsymbol v\|\le\lambda_D^{-1/2}|\boldsymbol v|_1\le C_{\mathrm P}|\boldsymbol v|_1$ (definition of $\lambda_D$, componentwise),
$|\operatorname{Res}(\boldsymbol v)|\le(C_{\mathrm P}\|\boldsymbol f+\Div\sigma\|+\|\sigma-G_h\|)|\boldsymbol v|_1=\eta_R|\boldsymbol v|_1$. Let $\boldsymbol r\in H^1_0(\Om)^2$ be the Riesz
representative, $(\nabla\boldsymbol r,\nabla\boldsymbol v)=\operatorname{Res}(\boldsymbol v)$ for all $\boldsymbol v$; then $|\boldsymbol r|_1=\|\operatorname{Res}\|_{-1}\le\eta_R$.

\emph{Step 3: velocity.} The second equation in \eqref{stokes} with $q=\Div\boldsymbol u\in L^2_0(\Om)$ gives $\Div\boldsymbol u=0$. Split $\boldsymbol e=\boldsymbol e_0+\boldsymbol e_\perp$
with $\boldsymbol e_0:=\Pi_Z\boldsymbol e\in Z$ and $\boldsymbol e_\perp\in Z^\perp$ (Lemma~\ref{lem:zperp}). Then $\Div\boldsymbol e_\perp=\Div\boldsymbol e=-\Div\boldsymbol u_h$, and
Lemma~\ref{lem:zperp} gives $|\boldsymbol e_\perp|_1\le\eta_D/\beta\le\eta_D/\beta_0$. For $\boldsymbol v\in Z$, $(e_P,\Div\boldsymbol v)=0$ and
$(\nabla\boldsymbol e_\perp,\nabla\boldsymbol v)=0$, so the error equation gives $(\nabla\boldsymbol e_0,\nabla\boldsymbol v)=\operatorname{Res}(\boldsymbol v)=(\nabla\boldsymbol r,\nabla\boldsymbol v)$.
Thus $\boldsymbol e_0\in Z$ and $\boldsymbol e_0-\boldsymbol r\perp Z$, i.e.\ $\boldsymbol e_0=\Pi_Z\boldsymbol r$, and $|\boldsymbol e_0|_1\le|\boldsymbol r|_1\le\eta_R$. By orthogonality,
$|\boldsymbol e|_1^2=|\boldsymbol e_0|_1^2+|\boldsymbol e_\perp|_1^2\le\eta_R^2+\eta_D^2/\beta_0^2$.

\emph{Step 4: pressure.} By the error equation, $(e_P,\Div\boldsymbol v)=(\nabla\boldsymbol e,\nabla\boldsymbol v)-\operatorname{Res}(\boldsymbol v)=(\nabla(\boldsymbol e-\boldsymbol r),\nabla\boldsymbol v)$
for all $\boldsymbol v$. Let $e_P^0:=e_P-|\Om|^{-1}\int_\Om e_P\,dA\in L^2_0(\Om)$. Then $\|e_P^0\|=\|e_P\|_{L^2/\R}$ (orthogonal projection onto constants), and
$(e_P^0,\Div\boldsymbol v)=(e_P,\Div\boldsymbol v)$ because $\int_\Om\Div\boldsymbol v\,dA=0$. The definition \eqref{beta} of $\beta$ gives
\[\beta\|e_P^0\|\le\sup_{\boldsymbol v\ne0}\frac{(e_P^0,\Div\boldsymbol v)}{|\boldsymbol v|_1}=\sup_{\boldsymbol v\ne0}\frac{(\nabla(\boldsymbol e-\boldsymbol r),\nabla\boldsymbol v)}{|\boldsymbol v|_1}=|\boldsymbol e-\boldsymbol r|_1 .\]
Since $\boldsymbol e-\boldsymbol r=\boldsymbol e_\perp-(I-\Pi_Z)\boldsymbol r$, we get $|\boldsymbol e-\boldsymbol r|_1\le|\boldsymbol e_\perp|_1+|(I-\Pi_Z)\boldsymbol r|_1\le\eta_D/\beta_0+\eta_R$. With
$\beta\ge\beta_0$ this is the second bound in \eqref{apost}.
\end{proof}
The velocity bound depends on $\beta_0$ through $\eta_D/\beta_0$; the pressure bound depends on it through $\eta_R/\beta_0$ and $\eta_D/\beta_0^2$. Without an explicit $\beta_0$ neither bound
is available for a method with $\Div\boldsymbol u_h\neq0$. For exactly divergence-free velocities the velocity bound is independent of $\beta$, but the pressure bound still is not.

\subsection*{Computation}
\emph{Discretisation and flux.} We take Taylor--Hood $\mathbb P_2/\mathbb P_1$ approximations on uniform triangulations with mesh size $h$ (each square of side $h$ is split into two
triangles, the diagonals alternating so that no triangle has all vertices on $\partial\Om$). The flux $\sigma$ is a continuous piecewise cubic tensor field, hence in $H(\Div,\Om)^2$,
and minimises $\|\sigma-G_h\|^2+t\|\boldsymbol f+\Div\sigma\|^2$ for $t\in\{1,4\}$; the smaller of the two resulting bounds is reported. The discrete solution and the flux are computed
in floating point; their coefficients are exact data for what follows.

\emph{Verified evaluation of $\eta_R$ and $\eta_D$.} All terms in $\eta_R$ and $\eta_D$ are integrals of squares of piecewise polynomials. The data $\boldsymbol f$ are polynomials with
integer or dyadic coefficients and all vertices are dyadic. Each integral is therefore an exact quadratic form in the nodal values on a sufficiently fine Lagrange lattice of the element,
with reference mass matrices computed in exact rational arithmetic. Everything is evaluated in interval arithmetic (\texttt{IntervalArithmetic.jl}~1.0.8, as used by \texttt{VFEM.jl}),
and every operation is flagged as guaranteed. An independent floating-point evaluation by collapsed Gauss quadrature lies inside the enclosures (relative agreement $10^{-12}$).

\emph{The Poincar\'e constant.} For the square and the rectangle $C_{\mathrm P}$ is exact, from $\lambda_D=\pi^2(a^{-2}+b^{-2})$ on $(0,a)\times(0,b)$. For $\Om_L$ and $\Om_H$ it comes from the
verified Crouzeix--Raviart lower eigenvalue bound \cite{Liu2015}, as implemented in \texttt{VFEM.jl} ($\lambda_D\ge\lambda_h^{\rm CR}/(1+(0.1893\,h)^2\lambda_h^{\rm CR})$ with a
verified algebraic eigensolver), on the triangulation with $h=1/16$:
\[\lambda_D(\Om_L)\ge9.549224,\qquad\lambda_D(\Om_H)\ge34.743177,\]
so that $\lambda_D(\Om_L)^{-1/2}\le0.3236056$ and $\lambda_D(\Om_H)^{-1/2}\le0.1696545$,
and we may take $C_{\mathrm P}=0.32361$ and $C_{\mathrm P}=0.16966$, respectively (on $h=1/32$ the same routine gives $9.609018$ and $35.424448$). Domain monotonicity would only give
$\lambda_D\ge\pi^2/2$. On the finest meshes the term $C_{\mathrm P}\|\boldsymbol f+\Div\sigma\|$ contributes less than $0.3\%$ of $\eta_R$ in all four examples; for $\Om_L$ and $\Om_H$, where
$C_{\mathrm P}$ is not known exactly, it contributes less than $0.4\%$ on every mesh (and up to $9\%$ on the coarsest meshes of the square and the rectangle, where $C_{\mathrm P}$ is exact). So the value
of $C_{\mathrm P}$ is immaterial here.

\emph{Comparison values of $\beta$.} As comparison we use the explicit lower bounds of Costabel and Dauge:
\begin{equation}\label{eq:cd-bounds}
\beta(\Om)\ge\sin\frac{\omega_{\rm HP}}2,\qquad\qquad
\beta(\Om)\ge\frac{R_{\rm in}}{\sqrt2\,R_{\rm out}}\Big(1+\sqrt{1-R_{\rm in}^2/R_{\rm out}^2}\Big)^{-1/2}.
\end{equation}
The first is the Horgan--Payne inequality with respect to the barycenter, with Horgan--Payne angle $\omega_{\rm HP}$; it holds for rectangles \cite[Thm.~2.2(i)]{CostabelDauge}. 
The second holds for domains contained in a ball of radius $R_{\rm out}$ and star-shaped with respect to a concentric ball of radius $R_{\rm in}$ \cite[Thm.~2.3]{CostabelDauge}. 
For the square $\omega_{\rm HP}=\pi/4$, which gives the Horgan--Payne value $\sin(\pi/8)$. For the rectangle $(0,2)\times(0,1)$, $\omega_{\rm HP}=\frac\pi2-\arctan2=\arctan\frac12$, which gives
$\beta\ge\sin(\frac12\arctan\frac12)\ge0.22975$. The L-shape $\Om_L$ is star-shaped with respect to the disk of radius $R_{\rm in}=\frac12$ centred at $(-\frac12,-\frac12)$ and contained in
the concentric disk of radius $R_{\rm out}=\sqrt{5/2}$, which gives $\beta\ge0.16018$. For $\Om_H$, which is not star-shaped, we are not aware of any previously available explicit lower
bound for $\beta$ (the estimates of \cite{HorganPayne,CostabelDauge} require star-shapedness).

\begin{table}[t]
\caption{Square $(-1,1)^2$, manufactured solution $\boldsymbol u=\curl\big((1-x^2)^2(1-y^2)^2\big)$, $P=x^2y$; Taylor--Hood $\mathbb P_2/\mathbb P_1$, $h=1/s$.
Guaranteed bounds \eqref{apost} with the certified $\beta_0=0.426245$ and with the Horgan--Payne value $\beta_0=\sin(\pi/8)$ (columns HP); all bounds rounded up.
The true errors (floating point, for comparison only) are $|\boldsymbol u-\boldsymbol u_h|_1$ and $\|P-P_h\|_{L^2/\R}$; ``eff'' is bound/true error.}\label{tab:app-square}
\centering\footnotesize\setlength{\tabcolsep}{4pt}
\begin{tabular}{rllllllll}
\toprule
$s$ & $\eta_R$ & $\eta_D$ & $|\boldsymbol u-\boldsymbol u_h|_1\le$ & eff & HP & $\|P-P_h\|_{L^2/\R}\le$ & eff & HP\\
\midrule
4 & $3.42\cdot10^{-1}$ & $1.97\cdot10^{-1}$ & $5.75\cdot10^{-1}$ & 1.85 & $6.18\cdot10^{-1}$ & $1.89$ & 45.6 & $2.24$\\
8 & $8.59\cdot10^{-2}$ & $5.18\cdot10^{-2}$ & $1.49\cdot10^{-1}$ & 1.83 & $1.61\cdot10^{-1}$ & $4.87\cdot10^{-1}$ & 35.6 & $5.78\cdot10^{-1}$\\
16 & $2.13\cdot10^{-2}$ & $1.30\cdot10^{-2}$ & $3.72\cdot10^{-2}$ & 1.81 & $4.01\cdot10^{-2}$ & $1.22\cdot10^{-1}$ & 31.4 & $1.45\cdot10^{-1}$\\
32 & $5.29\cdot10^{-3}$ & $3.25\cdot10^{-3}$ & $9.27\cdot10^{-3}$ & 1.80 & $9.99\cdot10^{-3}$ & $3.03\cdot10^{-2}$ & 30.8 & $3.60\cdot10^{-2}$\\
64 & $1.32\cdot10^{-3}$ & $8.10\cdot10^{-4}$ & $2.32\cdot10^{-3}$ & 1.79 & $2.50\cdot10^{-3}$ & $7.55\cdot10^{-3}$ & 31.1 & $8.97\cdot10^{-3}$\\
128 & $3.28\cdot10^{-4}$ & $2.03\cdot10^{-4}$ & $5.78\cdot10^{-4}$ & 1.79 & $6.23\cdot10^{-4}$ & $1.89\cdot10^{-3}$ & 31.4 & $2.24\cdot10^{-3}$\\
\bottomrule
\end{tabular}
\end{table}

\begin{table}[t]
\caption{Guaranteed bounds \eqref{apost} on the finest mesh for each domain, with the certified $\beta_0$ of this paper and with the comparison value of $\beta$ from
\eqref{eq:cd-bounds} (last three columns). ``gap'': relative decrease of the pressure bound if $\beta_0$ were replaced by the certified \emph{upper} bound of $\beta$, i.e.\ the
largest possible improvement from a better value of $\beta$. Data: manufactured solutions on the square (Table~\ref{tab:app-square}) and on the rectangle
($\boldsymbol u=\curl(x^2(2-x)^2y^2(1-y)^2)$, $P=(x-1)^2(y-\frac12)$); $\boldsymbol f=(-y,x)$ on the L-shape and on $\Om_H$.}\label{tab:app-all}
\centering\scriptsize\setlength{\tabcolsep}{3pt}
\begin{tabular}{lrlllllll}
\toprule
domain & $1/h$ & $\beta_0$ & $|\boldsymbol u-\boldsymbol u_h|_1\le$ & $\|P-P_h\|_{L^2/\R}\le$ & gap & comp.\ $\beta$ & $|\boldsymbol u-\boldsymbol u_h|_1\le$ & $\|P-P_h\|_{L^2/\R}\le$\\
\midrule
$(-1,1)^2$ & 128 & $0.426245$ & $5.78\cdot10^{-4}$ & $1.89\cdot10^{-3}$ & $2.6\cdot10^{-5}$ & $0.38268$ & $6.23\cdot10^{-4}$ & $2.24\cdot10^{-3}$\\
$(0,2)\times(0,1)$ & 128 & $0.387220$ & $1.44\cdot10^{-4}$ & $5.19\cdot10^{-4}$ & $1.7\cdot10^{-4}$ & $0.22975$ & $2.09\cdot10^{-4}$ & $1.21\cdot10^{-3}$\\
L-shape & 128 & $0.303973$ & $9.23\cdot10^{-3}$ & $3.82\cdot10^{-2}$ & $6.1\cdot10^{-4}$ & $0.16018$ & $1.70\cdot10^{-2}$ & $1.22\cdot10^{-1}$\\
$\Om_H$ & 128 & $0.154596$ & $2.60\cdot10^{-2}$ & $1.93\cdot10^{-1}$ & $2.2\cdot10^{-3}$ & -- & -- & --\\
\bottomrule
\end{tabular}
\end{table}

\begin{table}[t]
\caption{L-shaped domain and $\Om_H$, $\boldsymbol f=(-y,x)$: refinement history of the guaranteed bounds \eqref{apost} with the certified $\beta_0$
($0.303973$ and $0.154596$); Taylor--Hood $\mathbb P_2/\mathbb P_1$, $h=1/s$; all bounds rounded up. ``order'': $\log_2$ of the ratio of the bounds on
consecutive meshes. No exact solution is known, so no true errors are given.}\label{tab:app-hist}
\centering\footnotesize\setlength{\tabcolsep}{4pt}
\begin{tabular}{lrrllllll}
\toprule
domain & $s$ & triangles & $\eta_R$ & $\eta_D$ & $|\boldsymbol u-\boldsymbol u_h|_1\le$ & order & $\|P-P_h\|_{L^2/\R}\le$ & order\\
\midrule
L-shape & 4 & 96 & $2.36\cdot10^{-2}$ & $2.31\cdot10^{-2}$ & $7.94\cdot10^{-2}$ & -- & $3.27\cdot10^{-1}$ & --\\
 & 8 & 384 & $1.32\cdot10^{-2}$ & $1.32\cdot10^{-2}$ & $4.54\cdot10^{-2}$ & 0.81 & $1.87\cdot10^{-1}$ & 0.81\\
 & 16 & 1\,536 & $8.71\cdot10^{-3}$ & $8.49\cdot10^{-3}$ & $2.93\cdot10^{-2}$ & 0.63 & $1.21\cdot10^{-1}$ & 0.63\\
 & 32 & 6\,144 & $5.94\cdot10^{-3}$ & $5.72\cdot10^{-3}$ & $1.98\cdot10^{-2}$ & 0.57 & $8.14\cdot10^{-2}$ & 0.57\\
 & 64 & 24\,576 & $4.08\cdot10^{-3}$ & $3.90\cdot10^{-3}$ & $1.35\cdot10^{-2}$ & 0.55 & $5.56\cdot10^{-2}$ & 0.55\\
 & 128 & 98\,304 & $2.80\cdot10^{-3}$ & $2.68\cdot10^{-3}$ & $9.23\cdot10^{-3}$ & 0.55 & $3.82\cdot10^{-2}$ & 0.54\\
\midrule
$\Om_H$ & 4 & 96 & $3.05\cdot10^{-2}$ & $2.66\cdot10^{-2}$ & $1.75\cdot10^{-1}$ & -- & $1.31$ & --\\
 & 8 & 384 & $1.88\cdot10^{-2}$ & $1.76\cdot10^{-2}$ & $1.16\cdot10^{-1}$ & 0.60 & $8.56\cdot10^{-1}$ & 0.61\\
 & 16 & 1\,536 & $1.28\cdot10^{-2}$ & $1.21\cdot10^{-2}$ & $7.90\cdot10^{-2}$ & 0.54 & $5.87\cdot10^{-1}$ & 0.55\\
 & 32 & 6\,144 & $8.77\cdot10^{-3}$ & $8.33\cdot10^{-3}$ & $5.46\cdot10^{-2}$ & 0.53 & $4.06\cdot10^{-1}$ & 0.53\\
 & 64 & 24\,576 & $6.04\cdot10^{-3}$ & $5.75\cdot10^{-3}$ & $3.77\cdot10^{-2}$ & 0.54 & $2.80\cdot10^{-1}$ & 0.54\\
 & 128 & 98\,304 & $4.15\cdot10^{-3}$ & $3.96\cdot10^{-3}$ & $2.60\cdot10^{-2}$ & 0.54 & $1.93\cdot10^{-1}$ & 0.54\\
\bottomrule
\end{tabular}
\end{table}

\subsection*{Results}
Tables~\ref{tab:app-square}, \ref{tab:app-all} and~\ref{tab:app-hist} show the bounds.

\emph{Cost.} All bounds were computed on one server with at most $24$ cores. The finest meshes ($h=1/128$, up to $1.3\cdot10^5$ triangles) took between $3$ and $13$ minutes, dominated
by the sparse direct Taylor--Hood solve, with at most $15$\,GB of memory.

\emph{Effectivity (manufactured solutions).} For the two manufactured solutions (square and rectangle) the velocity bounds overestimate the true error by factors between $1.66$ and
$1.85$, and the pressure bounds by factors $16$--$46$, on all meshes; the true errors are computed in floating point from the exact solution. The large pressure factors are intrinsic to
\eqref{apost}: the bound controls $\beta\|P-P_h\|_{L^2/\R}$ by the size of the velocity residual, whereas the Taylor--Hood pressure error is, for the manufactured solutions, much smaller
than the velocity error. On the square and the rectangle the bounds converge with order $2$.

\emph{Convergence on non-convex domains.} For $\Om_L$ and $\Om_H$, with $\boldsymbol f=(-y,x)$, no exact solution is known and no reference solution was computed, so no true errors and
no effectivity indices are available there; Table~\ref{tab:app-hist} gives the full refinement history of the guaranteed bounds. The observed order of the bounds between
consecutive meshes decreases to $0.55$ on $\Om_L$ and $0.54$ on $\Om_H$ on the finest meshes. This is close to the exponent $\lambda_{\rm S}\approx0.5445$ of the
leading corner singularity of the Stokes problem at the re-entrant corners (the smallest positive root of $\sin(\lambda_{\rm S}\omega)=-\lambda_{\rm S}\sin\omega$,
$\omega=\frac{3\pi}2$), which limits the convergence on uniform meshes; this statement concerns the bounds, not the (unknown) errors.

\emph{Influence of $\beta$.} The inf-sup constant matters quantitatively. On the finest meshes the divergence term $\eta_D/\beta_0$ is $1.3$--$6.2$ times larger than $\eta_R$, and
$\eta_D/\beta_0^2$ accounts for $56$--$86\%$ of the pressure bound. With the comparison values \eqref{eq:cd-bounds} the bounds increase as follows: on the square, with the Horgan--Payne
value $\sin(\pi/8)$, by $8\%$ (velocity) and $19\%$ (pressure); on the rectangle by factors $1.46$ and $2.33$; on the L-shape by factors $1.83$ and $3.19$. For $\Om_H$, without an
explicit lower bound for $\beta$, bounds of the form \eqref{apost} cannot be stated. Conversely, the certified enclosures of $\beta$ are tight enough that no better value of $\beta$
could reduce the pressure bounds by more than $0.22\%$ (column ``gap'').

\begin{remark}
Computer-assisted existence proofs for the stationary Navier--Stokes equations have so far avoided the inf-sup constant. The approach of
\cite{WatanabeYamamotoNakao1999} uses the a priori estimate of \cite{NakaoYamamotoWatanabe1998}, which is based on the constant of Korn's inequality, and
\cite{LiuNakaoOishi} (on 3D domains) uses exactly divergence-free Scott--Vogelius approximations together with the a priori estimate of \cite{LiuNakaoYouOishi}.
Certified values of $\beta$, as obtained here for planar polygons, could make it possible to base such proofs on standard inf-sup stable pairs, such as the
Taylor--Hood pairs used above, and to bound the pressure; we have not carried this out here.
\end{remark}

\section{Conclusions}
The multiplier identity separates the two difficulties of the inf-sup constant on polygons. Near the boundary the quadratic form $\Re(wf^2)$ is dominated by $\kappa_0|f|^2$ as soon as
$\kappa_0$ exceeds the sharp corner constant, and the compact remainder is reduced to a matrix inequality that can be verified. The resulting lower bounds are unconditional and
convergent, and together with trial functions they enclose $\beta$ to three or four significant digits on the examples. Explicit values of $\beta$ of this quality are what guaranteed pressure error bounds for the planar Stokes equations (Theorem~\ref{thm:apost}) require. They may
also allow computer-assisted proofs for the Navier--Stokes equations to use discretizations that are not exactly divergence-free, which
\cite{LiuNakaoYouOishi,LiuNakaoOishi} avoid by using Scott--Vogelius elements.

For rectangles, isolation of $\beta^2$ below the essential spectrum is proven for the five aspect ratios $r\in\{\frac14,\frac38,\frac12,\frac9{16},\frac58\}$, of which
$r=\frac9{16},\frac58$ were not covered by \cite[Sect.~5.1]{CCDL}; no statement is made for other ratios. Open mathematical points are a convergence rate and the question whether
$\kappa=\kess$ for the square and for the rectangles with $r\in\{\frac34,\frac78\}$. The restriction of the present computations to dyadic vertices and axis-parallel sides is a
limitation of the implementation (exact integer arithmetic in the verifier, and the data builder), not an open problem of the method, whose analysis in
Sections~\ref{sec:fried}--\ref{sec:conv} applies to every polygon.

As an outlook, we expect an analogue of the multiplier identity in three dimensions, based on a divergence-free quadratic tensor for the first-order system satisfied by $\Div u$ and
$\curl u$, in which the vertices of a polyhedron would contribute in addition to its edges; this is not proved here and is left to future work.

\section*{Data availability}
The verifier, the design data, the logs and certificate files of all certifying runs, the trial functions with their certificates, the Lean sources and the
code and data of the Stokes application are available at \texttt{https://github.com/xfliu/lbb2d-verified} \cite{Repo}; the file \texttt{MANIFEST.csv} there
maps each certified number in the tables and the text to the run that produced it. The corrector dumps and Lean exports (28.5\,GB) are listed there with their SHA-256 checksums and are
available from the author on request.

\section*{Acknowledgments}
{\sloppy This work was supported by JSPS KAKENHI Grant Numbers JP24K00538 and JP26K22276. A generative AI assistant was used in drafting parts of the text and of
the verification code; the author checked all statements, proofs and certificates and takes full responsibility for the content.\par}

\appendix
\section{Background from complex analysis}\label{app:complex}
This appendix collects, in the form needed here, the facts from complex analysis used in the paper. Standard references are \cite{Ahlfors,Rudin} for holomorphic
functions, \cite{Hormander} for distributions, and \cite[\S5.2]{Lebl} for Bergman spaces of general domains (\cite{HKZ} treats the unit disk).

\emph{Wirtinger derivatives.} We identify $z=x+iy\in\C$ with $(x,y)\in\R^2$. For a complex-valued $C^1$ function $u$ on an open set,
\begin{equation}\label{eq:wirtinger}
\dz u:=\tfrac12(\partial_xu-i\partial_yu),\qquad\dzb u:=\tfrac12(\partial_xu+i\partial_yu).
\end{equation}
For $u\in L^1_{\rm loc}(\Om)$, $\dzb u$ is the distribution defined by $\langle\dzb u,\varphi\rangle:=-\int_\Om u\,\dzb\varphi\,dA$ for $\varphi\in C_c^\infty(\Om)$.
In effect the operators \eqref{eq:wirtinger} differentiate a polynomial in $z$ and $\bar z$ as if $z$ and $\bar z$ were independent variables; e.g.\ $\dzb z^n=0$ and $\dzb|z|^2=z$.
\begin{lemma}[Wirtinger calculus]\label{lem:wirtinger}
Let $u$ be a complex-valued $C^1$ function on an open set.
\begin{enumerate}
\item[(a)] $\dz$ and $\dzb$ are linear and obey the product rule; $\dz z=\dzb\bar z=1$, $\dzb z=\dz\bar z=0$, $\overline{\dz u}=\dzb\bar u$, and $4\dz\dzb u=\Delta u$ if $u\in C^2$.
\item[(b)] Pointwise,
\begin{equation}\label{eq:A1}
|\dz u|^2+|\dzb u|^2=\tfrac12|\nabla u|^2 .
\end{equation}
\item[(c)] Every $\R$-linear map $L\colon\C\to\C$ has the form $L(z)=az+b\bar z$ with $a:=\frac12(L(1)-iL(i))$ and $b:=\frac12(L(1)+iL(i))$. Hence every real affine map
$\C\to\C$ has the form $c_0+az+b\bar z$ with $c_0,a,b\in\C$, and its $\dzb$ is the constant $b$.
\item[(d)] For $u=h+ig$ with $h$, $g$ real,
\begin{equation}\label{eq:A2}
\dzb u=\tfrac12\big[(\partial_xh-\partial_yg)+i(\partial_yh+\partial_xg)\big];
\end{equation}
so $\dzb u=0$ if and only if $\partial_xh=\partial_yg$ and $\partial_yh=-\partial_xg$ (the \emph{Cauchy--Riemann equations}), i.e.\ $\nabla h=\curl g$, with $\curl$ as in
Section~\ref{sec:fried}.
\end{enumerate}
\end{lemma}
\begin{proof}
(a) Linearity and the product rule are inherited from $\partial_x$ and $\partial_y$. Directly from \eqref{eq:wirtinger}, $\dz z=\frac12(1-i\cdot i)=1$, $\dzb z=\frac12(1+i\cdot i)=0$,
and similarly for $\bar z$; $\overline{\dz u}=\frac12(\partial_x\bar u+i\partial_y\bar u)=\dzb\bar u$; and $4\dz\dzb=(\partial_x-i\partial_y)(\partial_x+i\partial_y)=\partial_x^2+\partial_y^2$
on $C^2$ functions.
(b) By the parallelogram law, $4(|\dz u|^2+|\dzb u|^2)=|\partial_xu-i\partial_yu|^2+|\partial_xu+i\partial_yu|^2=2|\partial_xu|^2+2|\partial_yu|^2=2|\nabla u|^2$.
(c) $L(x+iy)=xL(1)+yL(i)$, and $az+b\bar z=(a+b)x+i(a-b)y=L(1)x+L(i)y$. An affine map is a constant plus an $\R$-linear map, and $\dzb(c_0+az+b\bar z)=b$ by (a).
(d) Insert $\partial_xu=\partial_xh+i\partial_xg$ and $\partial_yu=\partial_yh+i\partial_yg$ into \eqref{eq:wirtinger}.
\end{proof}
By (c), $w=1-\dzb m$ is piecewise constant for a piecewise linear $m$.

\emph{Holomorphic functions.}
\begin{definition}\label{def:holo}
A function $u\in C^1(\Om)$ on an open set $\Om\subset\C$ is \emph{holomorphic} if $\dzb u=0$ in $\Om$.
\end{definition}
For $u\in C^1(\Om)$ this is equivalent to complex differentiability at every point, and then $u'=\dz u$ \cite{Ahlfors}. The following facts are classical (see
\cite[Ch.~2--4]{Ahlfors} or \cite[Ch.~10]{Rudin}); we state them in the form used.
\begin{lemma}\label{lem:holo}
Let $f,g$ be holomorphic in an open set $\Om\subset\C$.
\begin{enumerate}
\item[(a)] $f+g$, $fg$ and $f^n$ are holomorphic; $f$ is $C^\infty$, and all its derivatives are holomorphic.
\item[(b)] \emph{(Taylor expansion)} If $D(a,R)\subset\Om$ and $c_n:=f^{(n)}(a)/n!$, then $f(z)=\sum_{n\ge0}c_n(z-a)^n$, absolutely and uniformly on every closed disk
$\overline{D(a,\varrho)}$, $\varrho<R$.
\item[(c)] \emph{(Identity theorem)} If $\Om$ is connected and $f=0$ on a nonempty open subset, then $f=0$ in $\Om$.
\item[(d)] \emph{(Mean value property)} If $D(a,\varrho)\subset\Om$ and $0<s<\varrho$, then $f(a)=\frac1{2\pi}\int_0^{2\pi}f(a+se^{i\phi})\,d\phi$. Consequently
$f(a)=\int f(a-\zeta)j(\zeta)\,dA(\zeta)$ for every radial $j\ge0$ supported in $D(0,\varrho)$ with $\int j=1$; in particular
$f(a)=(\pi\varrho^2)^{-1}\int_{D(a,\varrho)}f\,dA$.
\item[(e)] \emph{(Logarithms and powers)} $\operatorname{Log}\zeta:=\ln|\zeta|+i\operatorname{Arg}\zeta$ with $\operatorname{Arg}\zeta\in(-\pi,\pi)$ is holomorphic on
$\C\setminus(-\infty,0]$, and so is $\zeta^s:=\exp(s\operatorname{Log}\zeta)$ for $s\in\R$. A continuous function $L$ on a connected set $U$ with $e^L=F$ is called a continuous logarithm
of $F$; two continuous logarithms of $F$ differ by a constant in $2\pi i\mathbb Z$. A holomorphic function without zeros on a disk has a holomorphic logarithm there.
\end{enumerate}
\end{lemma}
\begin{proof}
For (a)--(c), (e) see \cite[Ch.~2--4]{Ahlfors} or \cite[Ch.~10]{Rudin}; the uniqueness statement in (e) holds because the difference of two continuous logarithms is continuous on the
connected set $U$ with values in the discrete set $2\pi i\mathbb Z$. (d) By (b) the series converges uniformly on the circle $|z-a|=s$; integrating it term by term and using
$\int_0^{2\pi}e^{in\phi}d\phi=0$ for $n\ge1$ gives the circle mean. Integrating the circle means (the mean over $f(a-se^{i\phi})$ is the same) against the radial profile of $j$, in
polar coordinates, gives the statement for $j$; the disk mean is the case $j=(\pi\varrho^2)^{-1}\mathbf 1_{D(0,\varrho)}$.
\end{proof}

The next lemma is the case of the operator $\dzb$ of the hypoellipticity of elliptic operators with constant coefficients \cite{Hormander}; for completeness we give the elementary
mollifier proof.
\begin{lemma}[Weyl's lemma for $\dzb$]\label{lem:weyl}
Let $f\in L^1_{\rm loc}(\Om)$ satisfy $\int_\Om f\,\dzb\varphi\,dA=0$ for all $\varphi\in C_c^\infty(\Om)$, i.e.\ $\dzb f=0$ in the sense of distributions. Then $f$
coincides almost everywhere with a holomorphic function.
\end{lemma}
\begin{proof}
Let $j_\epsilon$ be a radial mollifier supported in $D(0,\epsilon)$, $\Om_\epsilon:=\{z\in\Om:\dist(z,\partial\Om)>\epsilon\}$ and $f_\epsilon:=f*j_\epsilon$ on $\Om_\epsilon$. Then
$f_\epsilon$ is $C^\infty$, and
$\dzb f_\epsilon(z)=\int f(\zeta)\,\partial_{\bar z}[j_\epsilon(z-\zeta)]\,dA(\zeta)=-\int f(\zeta)\,\partial_{\bar\zeta}[j_\epsilon(z-\zeta)]\,dA(\zeta)=0$ by the
hypothesis, because $\zeta\mapsto j_\epsilon(z-\zeta)$ belongs to $C_c^\infty(\Om)$. So $f_\epsilon$ is holomorphic on $\Om_\epsilon$. By Lemma~\ref{lem:holo}(d) a
holomorphic function is invariant under convolution with a radial mollifier, so on $\Om_{\epsilon+\epsilon'}$
\[f_\epsilon=f_\epsilon*j_{\epsilon'}=(f*j_{\epsilon'})*j_\epsilon=f_{\epsilon'}*j_\epsilon=f_{\epsilon'} .\]
Thus on each $\Om_{2\epsilon_0}$ all $f_\epsilon$ with $\epsilon<\epsilon_0$ coincide with one holomorphic function $F$. Since $f_\epsilon\to f$ in $L^1_{\rm loc}$, $f=F$
almost everywhere.
\end{proof}
The name refers to Weyl's lemma for the Laplacian (weakly harmonic functions are smooth and harmonic), of which this is a variant: $\dzb f=0$ implies $\Delta f=4\dz\dzb f=0$.

The following facts are standard for arbitrary bounded open sets (see e.g.\ \cite[\S5.2]{Lebl}); since \cite{HKZ} treats the unit disk, we include the short proofs.
\begin{lemma}[Bergman space]\label{lem:bergman}
Let $\Om$ be a bounded open set.
\begin{enumerate}
\item[(a)] If $f\in\A$ and $D(a,\varrho)\subset\Om$, then
\[|f(a)|\le\frac{\|f\|_{L^2(D(a,\varrho))}}{\sqrt\pi\,\varrho},\qquad|f'(a)|\le\frac{\sqrt2\,\|f\|_{L^2(D(a,\varrho))}}{\sqrt\pi\,\varrho^2}.\]
\item[(b)] $\A$ is a closed subspace of $L^2(\Om)$, hence a Hilbert space, and $\Az$ is a closed subspace of $\A$. Constants and polynomials in $z$ belong to $\A$.
\item[(c)] If $K\subset\Om$ is compact, every bounded sequence in $\A$ has a subsequence that converges uniformly on $K$, and the restriction $\A\to L^2(K)$ is a compact operator.
Consequently, for every $v\in L^\infty(K)$ the quadratic form $f\mapsto\Re\int_Kvf^2\,dA$ maps weakly convergent sequences in $\A$ to convergent sequences.
\end{enumerate}
\end{lemma}
\begin{proof}
(a) The proof of Lemma~\ref{lem:taylor} gives, with the Taylor coefficients $c_n$ of $f$ at $a$ (so $c_0=f(a)$ and $c_1=f'(a)$),
\[\|f\|^2_{D(a,\varrho)}=\sum_n|c_n|^2\frac{\pi\varrho^{2n+2}}{n+1}\ge\max\Big\{|c_0|^2\pi\varrho^2,\ \frac{|c_1|^2\pi\varrho^4}{2}\Big\}.\]
(b) If $f_n\in\A$ and $f_n\to f$ in $L^2$, then $\int f\dzb\varphi=\lim\int f_n\dzb\varphi=-\lim\int\varphi\,\dzb f_n=0$ for $\varphi\in C_c^\infty(\Om)$, so $f\in\A$ by
Lemma~\ref{lem:weyl}. The mean $f\mapsto\int_\Om f$ is continuous on $L^2(\Om)$, so its kernel $\Az$ is closed. Polynomials are bounded on the bounded set $\Om$.

(c) Let $(f_n)$ be bounded in $\A$, $d:=\dist(K,\partial\Om)/2$ and $K':=\{z:\dist(z,K)\le d\}\subset\Om$. For $z\in K'$, $D(z,d)\subset\Om$, so by (a) $|f_n|$ and $|f_n'|$ are bounded on $K'$ by a
constant times $\sup_n\|f_n\|$. For $z,z'\in K$ with $|z-z'|<d$ the segment $[z,z']$ lies in $K'$, so $|f_n(z)-f_n(z')|\le C|z-z'|$: the sequence is equicontinuous on $K$.
The Arzel\`a--Ascoli theorem gives a uniformly convergent subsequence. Now let $f_n\rightharpoonup f$ weakly in $\A$. Weakly convergent sequences are bounded, and by (a) the point
evaluations are continuous linear functionals on $\A$, so $f_n(z)\to f(z)$ for every $z\in K$. Hence every subsequence of $(f_n)$ has a further subsequence converging uniformly on $K$,
necessarily to $f$; as $K$ has finite measure, $f_n\to f$ in $L^2(K)$. This is the compactness of the restriction. Finally,
$|\int_Kv(f_n^2-f^2)|\le\|v\|_{L^\infty}\|f_n-f\|_{L^2(K)}\|f_n+f\|_{L^2(K)}\to0$.
\end{proof}

The following decomposition is classical and goes back to Friedrichs \cite{Friedrichs}; 
for completeness we include the short proof.
\begin{lemma}[orthogonal decomposition]\label{lem:decomp}
Let $\Om$ be a bounded open set. For $\varphi\in H^1_0(\Om;\C)$, $\|\dz\varphi\|=\|\dzb\varphi\|=\frac12|\varphi|_1$. The space $\dz H^1_0(\Om)$ is closed in
$L^2(\Om;\C)$, and
\[L^2(\Om;\C)=\A\oplus\dz H^1_0(\Om),\qquad(\Az)^\perp=\dz H^1_0(\Om)\oplus\C,\]
both sums being orthogonal with respect to $\langle u,v\rangle:=\int u\bar v$.
\end{lemma}
\begin{proof}
For $\varphi\in C_c^\infty(\Om;\C)$, integrating by parts twice (no boundary terms), and using $\overline{\dz\varphi}=\dzb\bar\varphi$ and $\dz\dzb=\dzb\dz$,
\[\int|\dz\varphi|^2=\int\dz\varphi\,\dzb\bar\varphi=-\int\dzb\dz\varphi\,\bar\varphi=-\int\dz\dzb\varphi\,\bar\varphi=\int\dzb\varphi\,\dz\bar\varphi=\int|\dzb\varphi|^2 .\]
Together with \eqref{eq:A1} this gives $\|\dz\varphi\|^2=\|\dzb\varphi\|^2=\frac14|\varphi|_1^2$, which extends to $H^1_0$ by density. So $2\dz$ is an
isometry from $H^1_0$ into $L^2$, and its range is closed. Next, $k\perp\dz C_c^\infty(\Om)$ means $0=\int k\,\overline{\dz\psi}=\int k\,\dzb\bar\psi$ for all $\psi$. As
$\bar\psi$ ranges over all of $C_c^\infty(\Om;\C)$, this says $\dzb k=0$ in $\mathcal D'(\Om)$, i.e.\ $k\in\A$ by Lemma~\ref{lem:weyl}; conversely every $k\in\A$ satisfies
it by Lemma~\ref{lem:holo}(a) and integration by parts. By density of $C_c^\infty$ in $H^1_0$, $(\dz H^1_0)^\perp=\A$, and since $\dz H^1_0$ is closed,
$L^2=\A\oplus\dz H^1_0$. Finally $\A=\Az\oplus\C$ orthogonally, because $\langle f,1\rangle=\int f$. Hence $(\Az)^\perp=(\A)^\perp\oplus\C=\dz H^1_0\oplus\C$; the last sum
is orthogonal because $\int\dz\varphi=0$ for $\varphi\in H^1_0$.
\end{proof}

\section{Auxiliary results}\label{app:aux}
Lemma~\ref{lem:findim} is a reformulation of Albert's nonnegativity criterion \cite{Albert1969}: for $t>0$ its hypothesis says
$\left(\begin{smallmatrix}S&x\\x^{\mathsf T}& t\end{smallmatrix}\right)\succeq0$, which holds if and only if $x\in\operatorname{ran}S$ and $x^{\mathsf T}S^+x\le t$. We include the short direct proof.
\begin{lemma}[finite-dimensional duality]\label{lem:findim}
Let $S$ be a real symmetric positive semidefinite $m\times m$ matrix, $x\in\R^m$ and $t\ge0$. Then
\[(y\cdot x)^2\le t\,y^{\mathsf T}Sy\quad\text{for all }y\in\R^m\qquad\Longleftrightarrow\qquad x=S\alpha\ \text{for some }\alpha\in\R^m\text{ with }\alpha^{\mathsf T}S\alpha\le t .\]
In this case every $\alpha$ with $S\alpha=x$ has $\alpha^{\mathsf T}S\alpha=x^{\mathsf T}S^+x\le t$, where $S^+$ is the Moore--Penrose pseudoinverse.
\end{lemma}
\begin{proof}
($\Leftarrow$) The Cauchy--Schwarz inequality for the semi-inner product $(y,y')\mapsto y^{\mathsf T}Sy'$ gives
$(y\cdot S\alpha)^2\le(y^{\mathsf T}Sy)(\alpha^{\mathsf T}S\alpha)\le t\,y^{\mathsf T}Sy$.
($\Rightarrow$) If $Sy=0$, then $y^{\mathsf T}Sy=0$ and hence $y\cdot x=0$; so $x\perp\ker S$. As $S$ is symmetric, $\ker S=(\operatorname{ran}S)^\perp$, so
$x\in(\ker S)^\perp=\operatorname{ran}S$. Let $\alpha_0:=S^+x$; since $SS^+$ is the orthogonal projection onto $\operatorname{ran}S$, $S\alpha_0=x$. If also $S\alpha=x$, then
$\alpha-\alpha_0\in\ker S$ and $\alpha^{\mathsf T}S\alpha=\alpha_0^{\mathsf T}S\alpha_0=x^{\mathsf T}S^+x$. Taking $y=\alpha$ in the hypothesis,
$(\alpha^{\mathsf T}S\alpha)^2=(\alpha\cdot x)^2\le t\,\alpha^{\mathsf T}S\alpha$, hence $\alpha^{\mathsf T}S\alpha\le t$ (trivially if $\alpha^{\mathsf T}S\alpha=0$).
\end{proof}
The lemma identifies the set of vectors $x$ satisfying the hypothesis with the ellipsoid $\{S\alpha:\alpha^{\mathsf T}S\alpha\le\mu\}$; it remains valid, with the same
constant, when $S$ is singular or ill conditioned. This is why Section~\ref{sec:impl} never inverts $\mathcal R_{\rm up}$.

\begin{lemma}\label{lem:nonvanish}
Let $\Om$ be a polygonal domain and $a\in\Om$. Then there is $f\in\Az$ with $f(a)\ne0$.
\end{lemma}
\begin{proof}
Suppose $f(a)=0$ for all $f\in\Az$. For an entire function $g$, $g-|\Om|^{-1}\int g\in\Az$, so $g(a)=|\Om|^{-1}\int_\Om g\,dA$. Choose
$\xi\in\C$, $|\xi|=1$, such that $z\mapsto\Re(\xi(z-a))$ is constant on no side and takes different values at different vertices; this excludes finitely many $\xi$.

\emph{The extreme vertex.} The maximum $M$ of this linear function over $\overline\Om$ is attained on $\partial\Om$, hence at a vertex, and by the choice of $\xi$ at exactly one vertex $v$;
moreover $M>0$, since $a$ is an interior point. The opening $\omega_v$ is less than $\pi$, because near a vertex with $\omega_v\ge\pi$ the set $\Om$ contains points in which the
linear function is larger than at $v$. Let $K:=\{v+re^{i\theta}:r\ge0,\ \theta_v\le\theta\le\theta_v+\omega_v\}$ be the tangent cone at $v$, and $R_0>0$ with
$\Om\cap D(v,R_0)=\operatorname{int}K\cap D(v,R_0)$. The function $\theta\mapsto\Re(\xi e^{i\theta})$ is negative at $\theta_v$ and at $\theta_v+\omega_v$ (the sides at $v$), and the
interval has length $<\pi$, so it is negative on the whole interval. Hence $\Re(\xi(z-v))\le-c_1|z-v|$ on $K$ for some $c_1>0$. On the compact set
$\overline\Om\setminus D(v,R_0)$ we have $\Re(\xi(z-v))\le-c_0<0$.

\emph{Growth of exponential moments.} For $s>0$ let
\[\begin{gathered}E(s):=\int_\Om e^{s\xi(z-a)}dA,\\ J_1:=\int_Ke^{s\xi(z-v)}dA,\qquad J_2:=\int_{\Om\setminus K}e^{s\xi(z-v)}dA,\qquad J_3:=\int_{K\setminus\Om}e^{s\xi(z-v)}dA .\end{gathered}\]
By assumption $E(s)=|\Om|$, and $e^{-s\xi(v-a)}E(s)=J_1+J_2-J_3$. Using $\int_0^\infty re^{-\lambda r}dr=\lambda^{-2}$ for $\Re\lambda>0$,
\[J_1=\int_{\theta_v}^{\theta_v+\omega_v}\frac{d\theta}{(s\xi e^{i\theta})^2}=\frac{e^{-2i\theta_v}-e^{-2i(\theta_v+\omega_v)}}{2i\,s^2\xi^2},\qquad|J_1|=\frac{\sin\omega_v}{s^2}.\]
The sets $\Om\setminus K$ and $K\setminus\Om$ lie outside $D(v,R_0)$ (up to a null set), so $|J_2|\le|\Om|e^{-sc_0}$ and
$|J_3|\le2\pi\int_{R_0}^\infty re^{-sc_1r}dr$, which decays exponentially in $s$. Since $|e^{-s\xi(v-a)}|=e^{-sM}$,
$|\Om|=|E(s)|\ge e^{sM}\big(s^{-2}\sin\omega_v-|J_2|-|J_3|\big)\to\infty$ as $s\to\infty$, a contradiction.
\end{proof}

\begin{proposition}[complex-linear representers]\label{prop:complexlin}
Let the triangulation, the weights $w_T$ and the passive and active sets $\mathcal P$, $\mathcal A$ be as in Section~\ref{sec:cert}, let $\kappa_0>\kappa_{\mathcal P}$, and let
$\mathsf A$ be the real-linear operator on $L^2(\Om;\C)$ with $\mathsf Ae:=\kappa_0e-\bar w_T\bar e$ on $T\in\mathcal P$ and $\mathsf Ae:=\kappa_0e$ on $T\in\mathcal A$, whose inverse is the
map $k\mapsto k^\sharp$ of Section~\ref{sec:conv}. Assume that $\mathcal A\ne\emptyset$, that some passive triangle has $w_T\ne0$, and that some Taylor functional $\ell_\mu$ does not vanish
on $\Az$. Let $k_1,\dots,k_N\in L^2(\Om)$ satisfy $\langle f,k_\nu\rangle=\ell_\nu(f)$ for all $f\in\Az$, $\nu=1,\dots,N$, and put $K_y:=\sum_\nu y_\nu k_\nu$ for $y\in\C^N$. Then there is $y\in\C^N$
with $K_y^\sharp\notin\Az$.
\end{proposition}
By the proposition, for the complex-linear combinations of \eqref{gram} the corrector part $e_y:=K_y^\sharp-P_0K_y^\sharp$ of the splitting in Lemma~\ref{lem:split}
($P_0$ the $\langle\cdot,\cdot\rangle_q$-orthogonal projection onto $\Az$) cannot vanish for all $y$, whatever the representers. This is why Lemma~\ref{lem:optcorr} uses separate
correctors for the real and imaginary parts of each Taylor functional.
\begin{proof}
Suppose $K_y^\sharp\in\Az$, i.e.\ $K_y\in\mathsf A(\Az)$, for all $y\in\C^N$. Taking $y$ equal to the $\mu$-th unit vector and to $i$ times it gives $k_\mu\in\mathsf A(\Az)\cap i\mathsf A(\Az)$.
Let $k_\mu=\mathsf Ag_1=i\mathsf Ag_2$ with $g_1,g_2\in\Az$. On the active triangles $\kappa_0g_1=i\kappa_0g_2$, so $g_1=ig_2$ in $\Om$ by the identity theorem
(Lemma~\ref{lem:holo}(c)). On a passive triangle with $w_T\ne0$, $0=\mathsf A(ig_2)-i\mathsf Ag_2=2i\bar w_T\bar g_2$, so $g_2=0$ there, hence everywhere, and $k_\mu=0$.
This contradicts $\ell_\mu\not\equiv0$ on $\Az$.
\end{proof}
The reason is that $\sharp$, and with it the optimal correction, is only real-linear. Theorem~\ref{thm:rel} holds for both versions; the implementation (Section~\ref{sec:impl})
uses the real-linear representers $\tilde k_j$.

\section{Proofs of results from the literature}\label{app:lit}
This appendix contains self-contained proofs of Theorem~\ref{thm:beta}, due to Costabel and Dauge \cite[Thm.~2.1]{CostabelDauge} (and, for smooth
domains, to Horgan and Payne \cite{HorganPayne}), and of Lemma~\ref{lem:kess}, which by Theorem~\ref{thm:beta} is \cite[Thm.~3.3]{CCDL} 
(see the attributions before the statements in Section~\ref{sec:fried}). They are included for completeness. The notation is that of Section~\ref{sec:fried};
in particular all functions in \eqref{beta} are real valued.

\subsection{The Friedrichs identity}\label{app:lit-beta}
We write $\langle\cdot,\cdot\rangle_{\mathcal D}$ for the pairing of distributions on $\Om$ with test functions.
By the definition of distributional derivatives (e.g.\ $\langle\partial_y\varphi,v_1\rangle_{\mathcal D}+\langle-\partial_x\varphi,v_2\rangle_{\mathcal D}
=-(\varphi,\partial_yv_1)+(\varphi,\partial_xv_2)$), for real $\varphi\in L^2(\Om)$ and $v\in C_c^\infty(\Om)^2$,
\begin{equation}\label{eq:distr}
(\varphi,\Div v)=-\langle\nabla\varphi,v\rangle_{\mathcal D},\qquad(\varphi,\rot v)=\langle\curl\varphi,v\rangle_{\mathcal D}.
\end{equation}

\begin{lemma}\label{lem:fq}
Let $q\in L^2_0(\Om)$, let $u_q$ be its Riesz representer \eqref{eq:uq}, and put $d:=\Div u_q$, $g:=\rot u_q$ and $h:=q-d$. Then
\begin{equation}\label{eq:beta-c}
(h,\Div v)=(g,\rot v)\qquad(v\in H^1_0(\Om)^2),
\end{equation}
the function $f:=h-ig$ belongs to $\Az$, and $(h,d)=\|g\|^2$.
\end{lemma}
\begin{proof}
By \eqref{eq:uq} and \eqref{eq:beta-a}, $(q,\Div v)=(\nabla u_q,\nabla v)=(d,\Div v)+(g,\rot v)$ for every $v\in H^1_0(\Om)^2$; this is \eqref{eq:beta-c}.

\emph{$f$ is holomorphic.} For $v\in C_c^\infty$, \eqref{eq:distr} turns \eqref{eq:beta-c} into $-\langle\nabla h,v\rangle_{\mathcal D}=\langle\curl g,v\rangle_{\mathcal D}$.
Hence $\nabla h=-\curl g$ in the sense of distributions, i.e.\ $\partial_xh=-\partial_yg$ and $\partial_yh=\partial_xg$. These are the Cauchy--Riemann
equations for the pair $(h,-g)$:
\[\dzb f=\tfrac12(\partial_x+i\partial_y)(h-ig)=\tfrac12\big[(\partial_xh+\partial_yg)+i(\partial_yh-\partial_xg)\big]=0\qquad\text{in }\mathcal D'(\Om).\]
Since $f\in L^2(\Om)$, Weyl's lemma (Lemma~\ref{lem:weyl}) shows that $f$ coincides almost everywhere with a holomorphic function, so $f\in\A$.

\emph{$f$ has mean zero.} $\int h=\int q-\int\Div u_q=0$ because $q\in L^2_0$ and by Lemma~\ref{lem:divrot}, and $\int g=\int\rot u_q=0$ by
Lemma~\ref{lem:divrot}. Hence $\int f=0$ and $f\in\Az$.

\emph{Orthogonality.} Testing \eqref{eq:beta-c} with $v=u_q$ gives $(h,d)=\|g\|^2$.
\end{proof}

\begin{lemma}\label{lem:imf}
For every $f\in\Az$,
\begin{equation}\label{eq:beta-e}
\|\Im f\|^2\ge\beta(\Om)^2\|f\|^2 .
\end{equation}
\end{lemma}
\begin{proof}
For $f=0$ there is nothing to prove. Let $0\ne f\in\Az$, $h:=\Re f$ and $g:=\Im f$, so $f=h+ig$ with $h,g$ real; let $a:=\|h\|^2$, $b:=\|g\|^2$. The
Cauchy--Riemann equations $\partial_xh=\partial_yg$, $\partial_yh=-\partial_xg$ (Lemma~\ref{lem:wirtinger}(d)) say $\nabla h=\curl g$ in $\Om$. Since
$\int f=0$, $h\in L^2_0$ and $\int g=0$. Both $a$ and $b$ are positive: if $h=0$, then $\nabla g=0$, so $g$ is constant on the connected set $\Om$ and
$\int g=0$ forces $g=0$, i.e.\ $f=0$; the case $g=0$ is the same with the roles exchanged.

By \eqref{eq:distr}, $(h,\Div v)=-\langle\nabla h,v\rangle_{\mathcal D}=-\langle\curl g,v\rangle_{\mathcal D}=-(g,\rot v)$ for $v\in C_c^\infty$, and by continuity
\begin{equation}\label{eq:beta-d}
(h,\Div v)=-(g,\rot v)\qquad(v\in H^1_0(\Om)^2).
\end{equation}
Let $v^*:=u_h$ be the Riesz representer of $h$ (Lemma~\ref{lem:riesz}), and $s:=|v^*|_1^2=(h,\Div v^*)$. By \eqref{eq:beta-b}, $|v^*|_1\ge\beta\|h\|$,
i.e.\ $s\ge\beta^2a$. The Cauchy--Schwarz inequality applied to $s=(h,\Div v^*)$ and, via \eqref{eq:beta-d}, to $s=-(g,\rot v^*)$ gives
\[s^2\le a\,\|\Div v^*\|^2,\qquad s^2\le b\,\|\rot v^*\|^2 .\]
Dividing by $a$ and by $b$, adding, and using Lemma~\ref{lem:divrot}: $s^2(\frac1a+\frac1b)\le\|\Div v^*\|^2+\|\rot v^*\|^2=|v^*|_1^2=s$. Hence
$s\le ab/(a+b)$ (trivially if $s=0$) and
\[\beta^2a\le s\le\frac{ab}{a+b},\qquad\text{i.e.}\qquad\beta^2(a+b)\le b .\]
Since $a+b=\|f\|^2$ and $b=\|\Im f\|^2$, this is \eqref{eq:beta-e}.
\end{proof}

\begin{proof}[Proof of Theorem~\ref{thm:beta}]
\emph{($\ge$)} Let $0\ne q\in L^2_0$, $u:=u_q$, and let $d$, $g$, $h$, $f$ be as in Lemma~\ref{lem:fq}. Let $a:=\|h\|^2$, $b:=\|g\|^2$ and
$D:=\|d\|^2$. By Lemma~\ref{lem:divrot}, $|u|_1^2=D+b$; since $q=h+d$ and $(h,d)=b$ (Lemma~\ref{lem:fq}), $\|q\|^2=a+2(h,d)+D=a+2b+D$.
Moreover $\|f\|^2=\int(h^2+g^2)=a+b$, and since $f^2=h^2-g^2-2ihg$, $\Re\int f^2=a-b$. By the Cauchy--Schwarz inequality, $b^2=(h,d)^2\le aD$.

If $a=0$, then $b=0$ by $b^2\le aD$, so $|u|_1^2/\|q\|^2=D/D=1\ge(1-\kappa)/2$ (here $D=\|q\|^2>0$, and $\kappa\ge0$).

If $a>0$, then $D\ge b^2/a$. The function $\Phi(t):=(t+b)/(a+2b+t)$ has $\Phi'(t)=(a+b)/(a+2b+t)^2>0$, so it is increasing on $[0,\infty)$ and
\[\frac{|u|_1^2}{\|q\|^2}=\Phi(D)\ge\Phi\Big(\frac{b^2}{a}\Big)=\frac{b^2+ab}{a^2+2ab+b^2}=\frac{b}{a+b}.\]
As $f\in\Az$ and $f\ne0$ ($\|f\|^2=a+b>0$), \eqref{eq:kappa-re} gives $a-b=\Re\int f^2\le\kappa(a+b)$, i.e.\ $2b\ge(1-\kappa)(a+b)$, i.e.\
$b/(a+b)\ge(1-\kappa)/2$.

In both cases $|u_q|_1^2/\|q\|^2\ge(1-\kappa)/2$, and taking the infimum over $q$ in \eqref{eq:beta-b} gives $\beta^2\ge(1-\kappa)/2$.

\emph{($\le$)} Let $0\ne f\in\Az$, and let $F:=e^{i\phi}f$ be as in Lemma~\ref{lem:phase}, so that $F\in\Az$, $\|F\|=\|f\|$ and $\int F^2=|\int f^2|$.
Writing $F=H+iG$ with $H,G$ real, $\int F^2=\|H\|^2-\|G\|^2+2i\int HG$; comparing real parts and applying \eqref{eq:beta-e} to $F$,
\[\Big|\int f^2\Big|=\|H\|^2-\|G\|^2=\|F\|^2-2\|G\|^2\le(1-2\beta^2)\|f\|^2 .\]
Taking the supremum over $f$ in \eqref{kappa} gives $\kappa\le1-2\beta^2$, i.e.\ $\beta^2\le(1-\kappa)/2$.
\end{proof}

\subsection{The corner bound}\label{app:lit-kess}
\begin{proof}[Proof of Lemma~\ref{lem:kess}]
Fix a corner $c$ with opening $\omega$, and let $\theta_0$, $R_0$, $S$, $t$, $c'$, $M$ and $g_\delta$ be as in \eqref{eq:gdelta}; put $d:=|c-c'|=t$. By
Lemma~\ref{lem:gdelta}, $g_\delta\in\A$ for $0<\delta<1$.

\emph{Local factorisation.} Fix $\rho$ with $0<\rho\le\min\{R_0,d/2\}$ and let $S_\rho:=S\cap D(c,\rho)$. On $D(c,\rho)$, $d/2\le|z-c'|\le3d/2$, so $z-c'$ has a holomorphic
logarithm $L_2$ on this disk; on $S_\rho$, $L_1(c+re^{i\theta}):=\ln r+i\theta$ ($\theta_0<\theta<\theta_0+\omega$) is a continuous logarithm of $z-c$. Both
$\operatorname{Log}M$ and $L_1-L_2$ are continuous logarithms of $M$ on the connected set $S_\rho$, so they differ by a constant $2\pi ik$, $k\in\mathbb Z$
(Lemma~\ref{lem:holo}(e)). Hence on $S_\rho$
\[g_\delta(z)=(z-c)^{-1+\delta}\varphi_\delta(z),\qquad (z-c)^{-1+\delta}:=r^{-1+\delta}e^{i(-1+\delta)\theta},\qquad\varphi_\delta:=e^{2\pi ik(-1+\delta)}e^{(1-\delta)L_2},\]
where $\varphi_\delta$ is holomorphic on $D(c,\rho)$, $|\varphi_\delta(z)|=|z-c'|^{1-\delta}$, and in particular $|\varphi_\delta(c)|^2=d^{2-2\delta}$. The gradients of
$\varphi_\delta^2$ and of $|\varphi_\delta|^2=|z-c'|^{2-2\delta}$ are bounded in modulus by $2|z-c'|^{1-2\delta}\le C_0:=2\max\{3d/2,2/d\}$ on $D(c,\rho)$, uniformly in
$\delta\in(0,1)$, so that
\[|\varphi_\delta(z)^2-\varphi_\delta(c)^2|\le C_0|z-c|,\qquad\big||\varphi_\delta(z)|^2-d^{2-2\delta}\big|\le C_0|z-c|\qquad(z\in D(c,\rho)).\]

\emph{Polar-coordinate integrals.} With $z=c+re^{i\theta}$ and $dA=r\,dr\,d\theta$,
\begin{gather*}
\int_{S_\rho}|z-c|^{-2+2\delta}dA=\int_{\theta_0}^{\theta_0+\omega}\!\!\int_0^\rho r^{-1+2\delta}\,dr\,d\theta=\frac{\omega\rho^{2\delta}}{2\delta},\\
I_\delta:=\int_{S_\rho}(z-c)^{-2+2\delta}dA=\int_0^\rho r^{-1+2\delta}dr\int_{\theta_0}^{\theta_0+\omega}e^{-2i(1-\delta)\theta}d\theta
=\frac{\rho^{2\delta}}{2\delta}\,e^{-2i(1-\delta)\theta_0}\,\frac{1-e^{-2i(1-\delta)\omega}}{2i(1-\delta)},
\end{gather*}
and, as $|1-e^{-2ix}|=2|\sin x|$, $|I_\delta|=\rho^{2\delta}|\sin((1-\delta)\omega)|/(2\delta(1-\delta))$. Finally
$\int_{S_\rho}|z-c|^{-2+2\delta}|z-c|\,dA=\omega\rho^{1+2\delta}/(1+2\delta)$ is bounded uniformly in $\delta\in(0,1)$.

\emph{The singular part.} Let $O(1)$ denote quantities bounded independently of $\delta\in(0,1)$. On $S_\rho$, $g_\delta^2=(z-c)^{-2+2\delta}\varphi_\delta^2$ and
$|g_\delta|^2=|z-c|^{-2+2\delta}|\varphi_\delta|^2$; inserting $\varphi_\delta^2=\varphi_\delta(c)^2+O(|z-c|)$ and $|\varphi_\delta|^2=d^{2-2\delta}+O(|z-c|)$,
\[\int_{S_\rho}g_\delta^2\,dA=\varphi_\delta(c)^2I_\delta+O(1),\qquad\int_{S_\rho}|g_\delta|^2dA=d^{2-2\delta}\frac{\omega\rho^{2\delta}}{2\delta}+O(1).\]

\emph{All other contributions are bounded.} On $\Om\setminus D(c,\rho)$, which contains the neighbourhoods of all other corners, $|g_\delta|\le\max\{1,\rho^{-1}\}\max\{1,D'\}$
with $D':=\sup_{z\in\Om}|z-c'|$ (the same bound holds with $\rho$ replaced by any $r>0$); so the integrals of $g_\delta^2$, $|g_\delta|^2$ and $g_\delta$ over
$\Om\setminus D(c,\rho)$ are $O(1)$. On $S_\rho$,
$\int_{S_\rho}|g_\delta|\le\max\{1,3d/2\}\int_{S_\rho}|z-c|^{-1+\delta}dA=\max\{1,3d/2\}\,\omega\rho^{1+\delta}/(1+\delta)=O(1)$. Hence $I(\delta):=\int_\Om g_\delta\,dA=O(1)$.

\emph{Mean correction and limit.} By Lemma~\ref{lem:mean}, $f_\delta:=g_\delta-\mu_\delta$ with $\mu_\delta:=I(\delta)/|\Om|$ lies in $\Az$, and
\[\int f_\delta^2=\int g_\delta^2-\frac{I(\delta)^2}{|\Om|},\qquad\|f_\delta\|^2=\|g_\delta\|^2-\frac{|I(\delta)|^2}{|\Om|}.\]
Therefore $2\delta\int f_\delta^2=2\delta\varphi_\delta(c)^2I_\delta+O(\delta)$ and $2\delta\|f_\delta\|^2=d^{2-2\delta}\omega\rho^{2\delta}+O(\delta)$. As $\delta\to0$,
$|2\delta\varphi_\delta(c)^2I_\delta|=d^{2-2\delta}\rho^{2\delta}|\sin((1-\delta)\omega)|/(1-\delta)\to d^2|\sin\omega|$ and $d^{2-2\delta}\omega\rho^{2\delta}\to d^2\omega>0$. Thus
$f_\delta\ne0$ for small $\delta$, and
\[\frac{|\int f_\delta^2|}{\|f_\delta\|^2}=\frac{|2\delta\int f_\delta^2|}{2\delta\|f_\delta\|^2}\longrightarrow\frac{|\sin\omega|}{\omega}.\]
By \eqref{kappa}, $\kappa(\Om)\ge|\sin\omega|/\omega$; maximising over the corners gives the claim.\end{proof}


\begin{thebibliography}{99}
\bibitem{Ahlfors} L. V. Ahlfors, \emph{Complex analysis}, 3rd ed., McGraw-Hill, New York, 1979.
\bibitem{Albert1969} A. Albert, \emph{Conditions for positive and nonnegative definiteness in terms of pseudoinverses}, SIAM J. Appl. Math. \textbf{17} (1969), no.~2, 434--440.
\bibitem{Arb} F. Johansson, \emph{Arb: efficient arbitrary-precision midpoint-radius interval arithmetic}, IEEE Trans. Comput. \textbf{66} (2017), no.~8, 1281--1292.
\bibitem{BabuskaAziz} I. Babu\v{s}ka and A. K. Aziz, \emph{Survey lectures on the mathematical foundations of the finite element method}, in: The mathematical foundations of the finite element method with applications to partial differential equations, Academic Press, New York, 1972, pp.~1--359.
\bibitem{BCDG} C. Bernardi, M. Costabel, M. Dauge, and V. Girault, \emph{Continuity properties of the inf-sup constant for the divergence}, SIAM J. Math. Anal. \textbf{48} (2016), no.~2, 1250--1271.
\bibitem{ChaumontFrelet} T. Chaumont-Frelet, \emph{Guaranteed stability bounds for second-order PDE problems satisfying a G\aa rding inequality}, arXiv:2601.00404, 2026.
\bibitem{CCDL} M. Costabel, M. Crouzeix, M. Dauge, and Y. Lafranche, \emph{The inf-sup constant for the divergence on corner domains}, Numer. Methods Partial Differential Equations \textbf{31} (2015), 439--458.
\bibitem{CostabelDauge} M. Costabel and M. Dauge, \emph{On the inequalities of Babu\v{s}ka--Aziz, Friedrichs and Horgan--Payne}, Arch. Ration. Mech. Anal. \textbf{217} (2015), no.~3, 873--898.
\bibitem{EvansGariepy} L. C. Evans and R. F. Gariepy, \emph{Measure theory and fine properties of functions}, revised ed., CRC Press, Boca Raton, FL, 2015.
\bibitem{Friedrichs} K. Friedrichs, \emph{On certain inequalities and characteristic value problems for analytic functions and for functions of two variables}, Trans. Amer. Math. Soc. \textbf{41} (1937), no.~3, 321--364.
\bibitem{Gallistl19} D. Gallistl, \emph{Rayleigh--Ritz approximation of the inf-sup constant for the divergence}, Math. Comp. \textbf{88} (2019), 73--89.
\bibitem{Gallistl21} D. Gallistl, \emph{A posteriori error analysis of the inf-sup constant for the divergence}, SIAM J. Numer. Anal. \textbf{59} (2021), no.~1, 249--264.
\bibitem{Grassle} B. Gr\"a\ss le, \emph{Guaranteed inf-sup bounds and existence verification for semilinear elliptic problems via nonconforming finite elements}, arXiv:2604.21887, 2026.
\bibitem{HSV2012} A. Hannukainen, R. Stenberg, and M. Vohral\'{\i}k, \emph{A unified framework for a posteriori error estimation for the Stokes problem}, Numer. Math. \textbf{122} (2012), no.~4, 725--769.
\bibitem{HKZ} H. Hedenmalm, B. Korenblum, and K. Zhu, \emph{Theory of Bergman spaces}, Graduate Texts in Mathematics, vol.~199, Springer, New York, 2000.
\bibitem{HorganPayne} C. O. Horgan and L. E. Payne, \emph{On inequalities of Korn, Friedrichs and Babu\v{s}ka--Aziz}, Arch. Rational Mech. Anal. \textbf{82} (1983), no.~2, 165--179.
\bibitem{Hormander} L. H\"ormander, \emph{The analysis of linear partial differential operators I}, 2nd ed., Springer, Berlin, 1990.
\bibitem{Lebl} J. Lebl, \emph{Tasty bits of several complex variables}, open textbook, \texttt{https://www.jirka.org/scv/}. 
\bibitem{Mathlib} The mathlib Community, \emph{The Lean mathematical library}, in: Proceedings of the 9th ACM SIGPLAN International Conference on Certified Programs and Proofs (CPP 2020), ACM, 2020, 367--381.
\bibitem{Liu2015} X. Liu, \emph{A framework of verified eigenvalue bounds for self-adjoint differential operators}, Appl. Math. Comput. \textbf{267} (2015), 341--355.
\bibitem{Repo} X. Liu, \emph{lbb2d-verified: certificates, verifier, Lean sources and data for this paper}, \texttt{https://github.com/xfliu/lbb2d-verified}, 2026.
\bibitem{LiuNakaoOishi} X. Liu, M. T. Nakao, and S. Oishi, \emph{Computer-assisted proof for the stationary solution existence of the Navier--Stokes equation over 3D domains}, Commun. Nonlinear Sci. Numer. Simul. \textbf{108} (2022), 106223.
\bibitem{NakaoYamamotoWatanabe1998} M. T. Nakao, N. Yamamoto, and Y. Watanabe, \emph{A posteriori and constructive a priori error bounds for finite element solutions of the Stokes equations}, J. Comput. Appl. Math. \textbf{91} (1998), no.~1, 137--158.
\bibitem{WatanabeYamamotoNakao1999} Y. Watanabe, N. Yamamoto, and M. T. Nakao, \emph{A numerical verification method of solutions for the Navier--Stokes equations}, Reliab. Comput. \textbf{5} (1999), no.~3, 347--357.
\bibitem{LiuNakaoYouOishi} X. Liu, M. T. Nakao, C. You, and S. Oishi, \emph{Explicit a posteriori and a priori error estimation for the finite element solution of Stokes equations}, Jpn. J. Ind. Appl. Math. \textbf{38} (2021), no.~2, 545--559.
\bibitem{PragerSynge} W. Prager and J. L. Synge, \emph{Approximations in elasticity based on the concept of function space}, Quart. Appl. Math. \textbf{5} (1947), 241--269.
\bibitem{ReedSimon} M. Reed and B. Simon, \emph{Methods of modern mathematical physics I: Functional analysis}, revised and enlarged ed., Academic Press, New York, 1980.
\bibitem{Repin2002} S. I. Repin, \emph{A posteriori estimates for the Stokes problem}, J. Math. Sci. (N.Y.) \textbf{109} (2002), no.~5, 1950--1964.
\bibitem{Repin2008} S. Repin, \emph{A posteriori estimates for partial differential equations}, Radon Series on Computational and Applied Mathematics, vol.~4, Walter de Gruyter, Berlin, 2008.
\bibitem{RepinStenberg} S. Repin and R. Stenberg, \emph{A posteriori error estimates for the generalized Stokes problem}, J. Math. Sci. (N.Y.) \textbf{142} (2007), 1828--1843.
\bibitem{Rudin} W. Rudin, \emph{Real and complex analysis}, 3rd ed., McGraw-Hill, New York, 1987.
\bibitem{Storn} J. Storn, \emph{Computation of the LBB constant for the Stokes equation with a least-squares finite element method}, SIAM J. Numer. Anal. \textbf{58} (2020), no.~1, 86--108.
\bibitem{Stoyan1999} G. Stoyan, \emph{Towards discrete Velte decompositions and narrow bounds for inf-sup constants}, Comput. Math. Appl. \textbf{38} (1999), no.~7--8, 243--261.
\end{thebibliography}
\end{document}